\documentclass[3p,times]{elsarticle}

\usepackage{array}
\usepackage{color} 
\usepackage{tabularx}
\usepackage{graphicx} 
\usepackage{amsmath}
\usepackage{amssymb}
\usepackage{amsfonts}	
\usepackage{moreverb}
\usepackage{dsfont}
\usepackage{bm}
\usepackage{multirow}
\usepackage{soul}

\newcommand\BibTeX{{\rmfamily B\kern-.05em \textsc{i\kern-.025em b}\kern-.08em
T\kern-.1667em\lower.7ex\hbox{E}\kern-.125emX}}

\usepackage{hyperref} 
\hypersetup{
    colorlinks=true,                          
    linkcolor=blue, 
    citecolor=red, 
    urlcolor=blue  } 

\newcommand{\x}{\mbf{x}}

\newcommand{\mbf}[1]{\mathbf{#1}}			%

\newcommand{\Q}{\mathbf{Q}}
\renewcommand{\S}{\mathbf{S}}
\renewcommand{\u}{\mathbf{u}}
\newcommand{\w}{\mathbf{w}}
\newcommand{\q}{\mathbf{q}}
\newcommand{\F}{\mathbf{F}}
\newcommand{\f}{\mathbf{f}}
\newcommand{\g}{\mathbf{g}}

\renewcommand{\v}{\mathbf{v}}

\newcommand{\G}{\mathbf{G}}
\newcommand{\A}{\mathbf{A}}

\newcommand{\halb}{\frac{1}{2}}

\newcommand{\be}{\begin{equation}}
\newcommand{\ee}{\end{equation}}
\newcommand{\bdm}{\begin{displaymath}}
\newcommand{\edm}{\end{displaymath}}

\newcommand{\bea}{\begin{eqnarray} }
\newcommand{\eea}{\end{eqnarray} }

\newcommand{\aposteriori}{\textit{a posteriori} }

\newcommand{\dev}{\textnormal{dev}} 

\newcommand{\AAA}{{\boldsymbol{A}}}

\newcommand{\GG}{{\mathbf{G}}}

\newcommand{\vv}{{\mathbf{v}}}
\newcommand{\xx}{{\mathbf{x}}}

\newcommand{\II}{{\mathbf{I}}}
\newcommand{\JJ}{{\mathbf{J}}}
\newcommand{\QQ}{{\mathbf{Q}}}
\newcommand{\PP}{{\mathbf{P}}}
\newcommand{\WW}{{\mathbf{W}}}

\newcommand{\Id}{{\mathbf{I}}}
\newcommand{\tr}{\textnormal{tr}}
\newcommand{\BS}{{\boldsymbol{\sigma}}}
\renewcommand{\Re}{\textnormal{Re}}

\newfont{\numerikEleven}{ecrm1000}
\newfont{\numerikTen}{cmss10}
\newfont{\numerikNine}{cmss9}
\newfont{\numerikEight}{cmss8}

\journal{Journal of Computational Physics}

\begin{document} 
\begin{frontmatter}
\title{High order ADER schemes for a unified first order hyperbolic formulation of continuum mechanics: 
viscous heat-conducting fluids and elastic solids} 
\author[UniTN]{Michael Dumbser$^{*}$}
\ead{michael.dumbser@unitn.it}
\cortext[cor1]{Corresponding author}

\author[UniPau]{Ilya Peshkov$^{**}$}
\ead{peshenator@gmail.com}
\cortext[cor1]{{Ilya Peshkov is on leave from Sobolev Institute of Mathematics , 4 Acad. Koptyug Avenue, 630090 Novosibirsk, Russia}}

\author[NSC]{Evgeniy Romenski}
\ead{evrom@math.nsc.ru}

\author[UniTN]{Olindo Zanotti}
\ead{olindo.zanotti@unitn.it}

\address[UniTN]{Department of Civil, Environmental and Mechanical Engineering, 
University of Trento, Via Mesiano 77, 38123 Trento, Italy.} 
\address[UniPau]{{Open and Experimental Center for Heavy Oil, Universit\'e de Pau et des Pays de l'Adour, Avenue de l'Universit\'e, 64012 Pau, France.}}
\address[NSC]{{Sobolev Institute of Mathematics, 4 Acad. Koptyug Avenue, 630090 Novosibirsk, Russia \& Novosibirsk State University, 2 Pirogova Str., 630090 Novosibirsk, Russia}}

\begin{abstract}
This paper is concerned with the numerical solution of the \textit{unified} first order hyperbolic formulation of continuum 
mechanics recently proposed by Peshkov \& Romenski \cite{PeshRom2014}, further denoted as \textit{HPR model}. In that 
framework, the viscous stresses are computed from the so-called \textit{distortion tensor} $\mathbf{A}$, which is one of  
the primary state variables in the proposed first order system. A very important key feature of the HPR model is its 
ability to describe \textit{at the same time} the behavior of inviscid and viscous compressible Newtonian and non-Newtonian 
\textit{fluids} with heat conduction, as well as the behavior of elastic and visco-plastic \textit{solids}. Actually, the 
model treats viscous and inviscid fluids as generalized visco-plastic solids. 
This is achieved via a stiff source term that accounts for strain relaxation in the evolution equations of $\mathbf{A}$. 
Also heat conduction is included via a first order hyperbolic evolution equation of the thermal impulse, from which the 
heat flux is computed. 
The governing PDE system is hyperbolic and fully consistent with the first and the second principle of thermodynamics. 
It is also fundamentally \textit{different} from first order Maxwell-Cattaneo-type relaxation models based on extended 
irreversible thermodynamics. The HPR model represents therefore a \textit{novel} and \textit{unified} description of 
continuum mechanics, which applies at the same time to \textit{fluid mechanics} and \textit{solid mechanics}. 
In this paper, the direct connection between the HPR model and the classical hyperbolic-parabolic Navier-Stokes-Fourier 
theory is established for the first time via a formal asymptotic analysis in the stiff relaxation limit. 

From a numerical point of view, the governing partial differential equations are very challenging, since they 
form a large nonlinear hyperbolic PDE system that includes stiff source terms and non-conservative products. We apply 
the successful family of one-step ADER-WENO finite volume (FV) and ADER discontinuous Galerkin (DG) finite element schemes
to the HPR model in the stiff relaxation limit, and compare the numerical results with exact or 
numerical reference solutions obtained for the Euler and Navier-Stokes equations. Numerical convergence 
results are also provided. To show the universality of the HPR model, the paper is rounded-off with an application 
to wave propagation in elastic solids, for which one only needs to switch off the strain relaxation source term in the 
governing PDE system.  

We provide various examples showing that for the purpose of \textit{flow visualization}, the distortion  
tensor $\mathbf{A}$ seems to be particularly useful. 

\end{abstract}

\end{frontmatter}


%
\section{Introduction} \label{sec:introduction}

\subsection{A unified first order hyperbolic approach to continuum mechanics} 

An attempt to build a unified and overarching formulation of continuum mechanics in first order hyperbolic form that includes fluid mechanics as well as solid mechanics 
has been very recently described by Peshkov and Romenski in~\cite{PeshRom2014}.  
The proposed model, hereafter the Hyperbolic Peshkov-Romenski (HPR) model,
can potentially cover the entire spectrum of viscous flows ranging from non-equilibrium gas dynamics to Newtonian and non-Newtonian fluids, and even elastic and plastic deformation in solids, provided that the continuum description is applicable. 
In order to make this possible, the \textit{material element}\footnote{In fluid mechanics, the terms fluid elements, fluid particles and fluid parcels are also used.} view point is employed and the very essence of any macroscopic flow, \textit{i.e.} the process of material element 
rearrangements, is explicitly described in the mathematical model. 
We note that the term \textit{material element} should be understood in
the conventional meaning of continuum mechanics, \textit{i.e.} as an \textit{ensemble} of a \textit{sufficiently large number} of molecules or atoms. 

An important difference between the HPR model and the classical continuum  models is that the material elements not only have a finite size, but they 
also have an internal structure,  which is subject to rearrangements, and which can be macroscopically described after introducing suitable quantities. 
Thus, in order to describe the deformability of material elements, a tensorial field\footnote{Rigorously speaking, $ \AAA $ is not a tensor field 
of rank 2, since it 
it transforms like a tensor of rank 1 with respect to a change of  coordinates.} 
$ \AAA(\xx,t)=[A_{ij}] $ is used. It maps the material elements from a current deformed state to the undeformed state, and it contains the information about deformation and rotation of material elements. While this approach is standard in the framework of solid mechanics, it is much less obvious for gas dynamics.
Because of the rearrangements of material elements, the field $ \AAA $ is not integrable in the sense that it does not relate Eulerian and Lagrangian coordinates of the continuum. As a result, the field $ \AAA $ is \textit{local}, see~\cite{PeshRom2014,God1978,GodRom1998,GodRom2003}. 
This is also the reason why we \textit{cannot} call $ \AAA $ the \textit{deformation gradient}, and thus, 
following \cite{God1978,GodRom1998,GodRom2003}, we shall instead refer to it as 
the \textit{material distortion field}, or simply the \textit{distortion tensor}.

In addition to the distortion field $ \AAA $, another important information 
is required to describe rearrangements in a system of material elements of finite size. This information should 
characterize how easy or how hard it is for material elements to rearrange (fluidity). 
In the kinetic theory of liquids,  Frenkel~\cite{Frenkel1955} proposed to use the average time $ \tau_{\rm F} $ 
between two solid-like vibration states of an atom to describe the ability of a liquid to flow.\footnote{Frenkel's ideas have been discussed, used and extended during the last 20 years to compute the thermodynamic and dynamic properties of liquids, see \cite{wallace1997statistical,chisolm2001dynamics,brazhkin2012two,bolmatov2013thermodynamic,bolmatov2015revealing} and references therein.}. 
Following this idea of Frenkel, it was proposed in~\cite{PeshRom2014} to use a continuum analog $ \tau $ of Frenkel's time $ \tau_{\rm F} $. 
Thus, in our continuum approach, the time $ \tau $ is the time taken by  a given material element to "escape" from the \textit{cage} composed of its neighbor elements, 
i.e. the time taken to   \textit{rearrange} with one of its neighbors. 
The more viscous a fluid is, the larger the 
time $\tau$, \textit{i.e.} the longer the fluid elements stay in contact with each other. 
The limiting cases, inviscid fluids and elastic solids, are recovered when $ \tau=0 $ and $ \tau=\infty $, respectively, 
while for viscous fluids, the time $ \tau$ is finite with $ 0<\tau<\infty $ (see the discussion in~\cite{PeshRom2014}). 
We shall call $ \tau $ the \textit{strain dissipation time}, because, in the 
mathematical formulation of the HPR model the inverse time $ \tau^{-1} $ defines the rate at which shear strains dissipate during the rearrangement process. 

Our material element point of view allows to formulate the system of governing partial differential equations (PDE) with rather convenient mathematical properties:
\begin{itemize}
\item
First, the model is described by a system of first order PDEs. 
We recall that first order systems are less sensitive to the quality of the computational mesh and in general they allow to get a numerical scheme 
of higher order of accuracy than for a second order model on the same discrete stencil. 

\item
Second, the model is \textit{hyperbolic} if the total energy potential is a \textit{convex} function of the state variables, see~\cite{PeshRom2014}. 
In other words, the model  is based on a \textit{wave} formulation. 
Indeed, from the point of view of the physics of wave propagation and because of the causality principle, any macroscopic transport phenomenon should be considered as a  wave propagation process. 
In particular, the momentum transfer in a viscous fluid in the transverse direction to the mean flow is nothing but a wave propagation process. These waves are known 
as the shear waves, which are very dissipative waves propagating over a distance that equals just a few wave lengths. Nevertheless, such waves give rise to very important phenomena known as 
boundary layers. Thus, one may expect that a physically based boundary layer theory has to be based on such a transverse wave dynamics. In full agreement with the above discussion, there are 
two types of waves in our hyperbolic model, \textit{longitudinal waves} and \textit{shear waves}, which transfer momentum in the transverse flow directions. 

\item
Third, the dissipative process of material element rearrangements is modeled by a stiff algebraic source term, \textit{i.e.} this term does not depend on the space derivatives, which automatically 
implies that the characteristic speeds of the corresponding hyperbolic system are always \textit{finite} (as they should), 
whatever the time $ \tau $ is. One may recall that in hyperbolic 
Maxwell-Cattaneo-type  models some characteristic speeds tend to infinity if the relaxation parameter tends to zero.  

\end{itemize}

We also note that the system of the governing equations discussed in~\cite{PeshRom2014} has already been derived by Godunov \& Romenski in the 
1970ies ~\cite{GodRom1972,God1978} in the context of
elasto-plastic deformation of metals, for which
it has been used  by several  
authors over the years~\cite{Rom1989,GavrFavr2008,BartonRom2010,Pesh2010,Barton2013}.
On the contrary, the idea that the same model could also describe the dynamics of any 
continuum, including inviscid fluids, viscous Newtonian and non-Newtonian fluids, elastic and visco-plastic solids was discussed in~\cite{PeshRom2014} for the first time. 
In order to allow a quantitative comparison also with the Fourier heat conduction theory, in this paper we extend the model proposed by Peshkov and Romenski in \cite{PeshRom2014}  
by including also hyperbolic heat conduction equations, as proposed by Romenski in~\cite{Rom1986heat,Rom1989,RomToro2007,RomDrikToro2010}. The essential difference of our 
hyperbolic heat conduction model from that proposed by Cattaneo~\cite{Cattaneo1948} is that the speed of the heat propagation front is always \textit{finite},
whatever the heat flux relaxation parameter is. 

We emphasize that it is not our aim to provide a link with kinetic theory, although this could be very illuminating, but rather to verify the capabilities of the HPR model 
to account for a wide variety of dynamical systems.

\subsection{High order ADER-WENO finite volume and ADER discontinuous Galerkin finite element schemes} 

The resulting governing partial differential equations of the HPR model, introduced in \cite{PeshRom2014} and presented later in Section \ref{sec:model},  
are rather challenging from a numerical point of view, since they constitute a \textit{large system} of nonlinear hyperbolic conservation laws that also 
includes \textit{non-conservative products} and \textit{stiff source terms}. To the best knowledge of the authors, the complete first order HPR model presented in  
\cite{PeshRom2014} has never been solved so far by any numerical method in multiple space dimensions and including all terms, hence one of the main goals of 
this paper is to thoroughly investigate the behavior of the HPR model in a large number of different standard benchmark problems of computational fluid mechanics 
and computational solid mechanics.  

It is important to mention that exactly for such a general class of nonlinear time-dependent hyperbolic PDEs, the families of ADER finite volume (FV) and ADER 
discontinuous Galerkin (DG) finite element methods have been developed in the past decade. The starting point of the original ADER (arbitrary high order derivatives) 
schemes of Toro \& Titarev et al. for hyperbolic conservation laws \cite{toro1,schwartzkopff,toro3,toro4,titarevtoro,Toro:2006a,dumbser_jsc,CastroToro,toro-book} was the 
approximate solution of the generalized Riemann problem (GRP) \cite{LeFloch:1991a,BenArtzi:2006a} that arises naturally in the 
context of high order finite volume and DG schemes, due to their piecewise high order polynomial data representation, for which the vector of conserved variables 
and all its spatial derivatives are known at a given time level.  
The ADER approach has been successfully extended also to hyperbolic PDEs with stiff source terms \cite{DumbserEnauxToro,DumbserZanotti,HidalgoDumbser,Montecinos2014d}, 
to hyperbolic PDEs with non-conservative products \cite{ADERNC,USFORCE2} and to parabolic problems \cite{MunzDiffusionFlux,Hidalgo2009,ADERNSE}. Recent developments 
include space-time adaptive meshes \cite{AMR3DCL,AMR3DNC,Zanotti2015}, moving meshes \cite{LagrangeNC,Lagrange3D}, ADER-WENO finite volume schemes for 
divergence-free magnetohydrodynamics \cite{Balsara2009,Balsara2013,ADERdivB} and \textit{a posteriori} limiting of high order ADER-DG and ADER-FV schemes 
\cite{ADER_MOOD_14,Dumbser2014,Zanotti2015a,Zanotti2015b}. In the context of ADER schemes, first order hyperbolic reformulations of parabolic viscous problems 
have been tackled by Toro and Montecinos in \cite{Montecinos2014a,Montecinos2014b,Montecinos2014c}, while a series of interesting previous work on first order 
hyperbolic reformulations of advection-diffusion equations was proposed by Nishikawa in \cite{Nishikawa1,Nishikawa2}. Although not directly related to viscous
problems, we also would like to refer to the well-known relaxation system of Jin and Xin \cite{JinXin}, which allows to reformulate any nonlinear hyperbolic 
conservation law as an augmented linear first order system with stiff relaxation source terms.  

In this paper, we concentrate our attention on compressible viscous Newtonian fluids, which in the classical continuum theory can be described by the 
hyperbolic-parabolic Navier-Stokes-Fourier (NSF) theory, as well as on elastic solids. 
It should also be noted that there are several advantages of a first order hyperbolic formulation of viscous fluids: first, the use of explicit Godunov-type 
shock-capturing finite volume schemes and, 
even more, the use of high order discontinuous Galerkin finite element methods is - at least in principle - \textit{straightforward} for first order 
systems, while DG schemes need some special care in the presence of parabolic and higher order derivative terms, see the very interesting discussions in 
the well-known papers of Bassi \& Rebay \cite{BassiRebay}, Baumann \& Oden \cite{Baumann1,Baumann2}, Cockburn and Shu \cite{CBS-convection-diffusion,CBS-convection-dominated},  
Yan and Shu \cite{YanShu,Yan2002,LevyShuYan} and others \cite{ArnoldBrezzi,HartmannHouston,HartmannHouston2,KlaijVanDerVegt,Feistauer4,MunzDiffusionFlux,DumbserFacchini}.
Second, the use of a parabolic theory can lead to a severe time step size restriction, if explicit time stepping schemes are used, since the infinite 
propagation speed of perturbations that is intrinsically inherent in parabolic PDEs is reflected in explicit numerical methods by a stability condition 
on the time step that scales with the square of the mesh size, while it scales only linearly with the mesh size for first order hyperbolic systems due to 
the classical CFL condition \cite{CFL}.  
The situation is even worse for high order discontinuous Galerkin finite element schemes, where the explicit time step size scales not only quadratically 
with the mesh size, but where it decreases even quadratically with the order of the method. In Section \ref{sec:results} we will show one numerical example 
with an explicit time stepping scheme, where the use of the first order HPR model is clearly more convenient in terms of time step size and CPU time compared 
to the classical parabolic Navier-Stokes theory. 
As a third and last advantage of a first order hyperbolic model, we would like to emphasize that, by avoiding 
the presence of \textit{infinite wave speeds} even in the Newtonian framework, the new formulation 
suggests that its extension to \textit{relativistic continuum mechanics} should also be possible.

\subsection{Outline of the paper} 

The rest of this paper is organized as follows: in Section \ref{sec:model} we recall and discuss the extended hyperbolic Peshkov-Romenski model, denoted by HPR model in the 
following, including also a hyperbolic formulation of heat conduction. In particular, we show that the system is thermodynamically consistent and symmetric hyperbolic. 
A sketch of the analysis of the characteristics of the model is provided, together with a dispersion analysis of the wave speeds for relaxation times ranging from zero to 
infinity. We also carry out a formal asymptotic analysis of the system in the stiff relaxation limit, which reveals the direct connection of the first order HPR model 
with the well-established hyperbolic-parabolic Navier-Stokes-Fourier equations of viscous heat conducting fluids. 
In Section \ref{sec:ader} we briefly summarize the numerical methods used to solve the HPR model in this paper, 
namely ADER-WENO finite volume schemes and ADER discontinuous Galerkin finite element methods, making use of the unified $P_NP_M$ framework established in 
\cite{Dumbser2008}, which contains FV schemes and DG methods as two special cases of a more general class of numerical methods. 
In Section \ref{sec:results} we present computational results for a large set of different multi-dimensional test problems from computational fluid mechanics  
and also one example from computational solid mechanics, ranging from viscous low Mach number flows over viscous and inviscid compressible flows to the simulation of 
wave propagation in elastic solids. 
The paper is rounded-off by some concluding remarks and an outlook to future research in Section \ref{sec:conclusion}.  

\newpage

\section{Presentation and discussion of the mathematical model} 
\label{sec:model} 

\subsection{Formulation of the model}

The unified first order hyperbolic model for continuum mechanics proposed by Peshkov \& Romenski in \cite{PeshRom2014}, including a hyperbolic formulation of 
heat conduction, reads: 
\begin{subequations}\label{eqn.HPR}
	\begin{align}
	& \frac{\partial \rho}{\partial t}+\frac{\partial \rho v_k}{\partial x_k}=0,\label{eqn.conti}\\[2mm]
	&\displaystyle\frac{\partial \rho v_i}{\partial t}+\frac{\partial \left(\rho v_i v_k + p \delta_{ik} - \sigma_{ik} \right)}{\partial x_k}=0, \label{eqn.momentum}\\[2mm]
	&\displaystyle\frac{\partial A_{i k}}{\partial t}+\frac{\partial A_{im} v_m}{\partial x_k}+v_j\left(\frac{\partial A_{ik}}{\partial x_j}-\frac{\partial A_{ij}}{\partial x_k}\right)
	=-\dfrac{ \psi_{ik} }{\theta_1(\tau_1)},\label{eqn.deformation}\\[2mm]
	&\displaystyle\frac{\partial \rho J_i}{\partial t}+\frac{\partial \left(\rho J_i v_k+ T \delta_{ik}\right)}{\partial x_k}=-\dfrac{\rho H_i}{\theta_2(\tau_2)}, \label{eqn.heatflux}\\[2mm]
	&\displaystyle\frac{\partial \rho s}{\partial t}+\frac{\partial \left(\rho s v_k + H_k \right)}{\partial x_k}=\dfrac{\rho}{\theta_1(\tau_1) T} \psi_{ik} \psi_{ik} + \dfrac{\rho}{\theta_2(\tau_2) T} H_i H_i \geq0, \label{eqn.entropy}
	\end{align}
\end{subequations}
The solutions of the above PDE system fulfill also the additional conservation law 
\begin{equation}\label{eqn.energy}
\frac{\partial \rho  E}{\partial t}+\frac{\partial \left(v_k \rho  E + v_i (p \delta_{ik} - \sigma_{ik}) + q_k \right)}{\partial x_k}=0, 
\end{equation}
which is the conservation of total energy. Actually, in the numerical computations shown later in Section \ref{sec:results} of this paper, we solve the energy equation 
\eqref{eqn.energy} instead of the entropy equation ~\eqref{eqn.entropy}, but from the point of view of the model formulation, the entropy should be considered among the 
vector of unknowns (see Section~\ref{sec.therm.hyperb} for a discussion). 

Here we use the following notation:  $\rho$ is the mass density, $ [v_i]=\vv=(u,v,w) $ is the velocity vector, $ [A_{ik}]=\AAA $ is the distortion tensor, $ [J_i]=\JJ $ is the 
thermal impulse vector, $ s $ is the entropy, $ E =E(\rho,s,\vv,\AAA,\JJ)$ is the total energy, $ p = \rho^2E_{\rho} $ is the pressure, $ \delta_{ik} $ is the Kronecker delta, 
$ [\sigma_{ik}]=\boldsymbol{\sigma} = - [\rho A_{mi}E_{A_{mk}}] $ is the symmetric viscous shear stress tensor, $ T=E_s $ is the temperature, $[q_k] = \mathbf{q} = [E_s E_{J_k}]$  
is the heat flux vector and $ \theta_1 = \theta_1(\tau_1) > 0$ and $ \theta_2 = \theta_2(\tau_2) > 0$ are positive scalar functions, which will be specified below, 
depending on the strain dissipation time $ \tau_1 > 0$ and the thermal impulse relaxation time $ \tau_2 > 0$, respectively. The dissipative terms $\psi_{ik}$ and $H_i$ on the right hand side  
of the evolution equations for $\AAA$, $\JJ$ and $s$ are defined as $[\psi_{ik}] = \boldsymbol{\psi} = [E_{A_{ik}}]$ and $[H_i] = \mathbf{H} = [E_{J_i}]$, respectively. 
Hence, the viscous stress tensor and the heat flux vector are directly related to the dissipative terms on the right hand side via 
$ \boldsymbol{\sigma} = - \rho \AAA^T \boldsymbol{\psi} $ and $ \mathbf{q} = T \, \mathbf{H}$.     
Note that $ E_\rho $, $ E_s $, $ E_{A_{ik}} $ and $ E_{J_i} $ should be understood as the partial derivatives $ \partial E/\partial \rho$, $ \partial E/\partial s $,  
$ \partial E/\partial A_{ik}$ and $ \partial E/\partial J_i$; they are the so-called \textit{energy gradients in the state space} or the \textit{thermodynamic forces}. 
The Einstein summation convention over repeated indices is implied.

These equations are the mass conservation (\ref{eqn.conti}), the momentum conservation~(\ref{eqn.momentum}), the time evolution for the distortion~(\ref{eqn.deformation}), the time evolution for the thermal impulse~(\ref{eqn.heatflux}), the entropy time evolution~(\ref{eqn.entropy}), and the total energy conservation~(\ref{eqn.energy}). 
The PDE governing the time evolution of the thermal impulse~\eqref{eqn.heatflux} looks formally very similar to the momentum equation \eqref{eqn.momentum}, where the temperature $T$ takes the role of the pressure $p$.  Due to this similarity, it will also be called the \textit{thermal momentum equation} in the following.

One can clearly see that in order to close the system, it is necessary to specify the total energy potential $ E(\rho,s,\vv,\AAA,\JJ) $. This potential then generates all the constitutive fluxes (\textit{i.e.} non advective fluxes)  and source terms by means of its partial derivatives with respect to the state variables. Hence, the energy specification is one of the key steps in the model formulation.

In order to specify $ E $, we note that there are three scales involved in the continuum model formulation described in the introduction. Namely, the molecular scale, or
the  \textit{microscale}; the scale of the material elements, called here \textit{mesoscale}; and the flow scale, or the \textit{macroscale}. It is therefore assumed that the total energy $ E $ is the sum of three terms, each of which represents the energy distributed in its corresponding scale. Thus, we assume that
\begin{equation}\label{eq:total_energy}
E(\rho,s,\vv,\AAA,\JJ)=E_{1}(\rho,s)+E_{2}(\AAA,\JJ) + E_{3}(\vv).
\end{equation}
The terms $ E_3 $ and $ E_1 $ are conventional. They are the specific kinetic energy per unit mass $ E_3(\vv) =\dfrac{1}{2} v_iv_i$, which represents the 
macroscale part of the total energy, and the \textit{internal energy}  $ E_1(\rho,s) $, which is related to the kinetic energy of the molecular motion. 
$ E_1(\rho,s) $ is the only energy which does not disappear in the thermodynamic equilibrium where any meso- and macroscopic dynamics are absent, and only molecular dynamics is present. For this reason, it  is sometimes referred to as the \textit{equilibrium energy}. In this paper, for $ E_1 $, we shall use either the \textit{ideal gas equation of state} 
\begin{equation}\label{eq:ideal_gas_eos}
E_1(\rho,s)=\frac{c_0^2}{\gamma(\gamma-1)},\ \ c_0^2=\gamma\rho^{\gamma-1}e^{s/c_V},
\end{equation}
or the \textit{stiffened gas equation of state}
\begin{equation}\label{eq:stiff_gas_eos}
E_1(\rho,s)=\dfrac{c^2_0}{\gamma(\gamma-1)}\left ( \dfrac{\rho}{\rho_0}\right )^{\gamma-1}e^{s/c_V}+\dfrac{\rho_0 c_0^2-\gamma p_0}{\gamma\rho},\ \ c_0^2=const.
\end{equation}
In both cases, $c_0$ has the meaning of the adiabatic sound speed; $ c_V $ and $c_p$ are the specific heat capacities at constant volume and at constant pressure, respectively, 
which are related by the ratio of specific heats $ \gamma=c_p/c_V $. In~\eqref{eq:stiff_gas_eos}, $ \rho_0 $ is the reference mass density and $ p_0 $ is the reference 
(atmospheric) pressure. 

For the mesoscopic, or \textit{non-equilibrium}, part of the total energy, we shall use a quadratic form
\begin{equation}\label{eq:e_2}
E_2(\AAA,\JJ)=\dfrac{c_s^2}{4}G^{\rm TF}_{ij}G^{\rm TF}_{ij}+\frac{\alpha^2}{2}J_i J_i,
\end{equation}
with 
\begin{equation}
 [G_{ij}^{\rm TF}] = \dev(\GG) = \GG-\frac{1}{3} {\rm tr}(\GG) \II,  \qquad \textnormal{ and } \qquad \GG=\AAA^\mathsf{T}\AAA. 
\end{equation} 
Here, $[G_{ij}^{\rm TF}] = \dev(\GG)$ is the deviator, or the \textit{trace-free} part, of the tensor $\GG=\AAA^\mathsf{T}\AAA$ and ${\rm tr}(\GG)=G_{ii}$ is its trace, 
$ \II $ is the unit tensor and $ c_s $ is the characteristic velocity of propagation of transverse perturbations. 
In the following we shall refer to it as the \textit{shear sound velocity}. The characteristic velocity of heat wave propagation $c_h$ is related to 
$\alpha$\footnote{The physical units of $\alpha$ are $ \rm kg/(K\cdot m \cdot s^2) $.}, as discussed later in Section~\ref{sec.charcteristics}.  
We stress that $E_2(\AAA,\JJ)$ is a simple \textit{quadratic form} in terms of $G_{ij}^{\rm TF}$ and $\JJ$.

We also note that, because of the frame invariance principle, or \textit{objectivity principle}, the total energy can depend on vectors and tensors by means of their invariants only. By a direct calculation, one can see that
\[G_{ij}^{\rm TF} G_{ij}^{\rm TF}\equiv I_2-I_1^2/3,\]
where $ I_1={\rm tr}(\GG) $ and $ I_2={\rm tr}(\GG^2) $, and therefore $ E_2 $, as well as the total energy $ E $, are a function of invariants of $ \AAA $ and $ \JJ $.


In general, the mesoscopic energy $E_2(\AAA,\JJ)$ can also be a function of $ \rho $ and $ s $ in addition to $ \AAA $ and $ \JJ $. This would
correspond to a coupling between the  molecular scale and the scale of material elements. 
Such a  dependence on $ \rho $ and $ s $ should be introduced in the velocities $ c_s $ and $ \alpha $, \textit{i.e.} $ c_s=c_s(\rho,s) $, $ \alpha=\alpha(\rho,s) $.  The dependencies  
$ c_s(\rho,s) $ and $ \alpha(\rho,s) $ should be taken into account when strongly non-equilibrium flows are considered. This would affect the computation of the pressure and 
of the temperature through the partial derivatives $ E_\rho $ and $ E_s $ and give rise to a so-called \textit{non-equilibrium pressure} and a \textit{non-equilibrium temperature}. 
For simplicity, however, in this paper we do \textit{not} consider such a possibility, and $ c_s $ and $ \alpha $ are assumed to be 
\textit{constant}.

The algebraic source term on the right-hand side of equation~(\ref{eqn.deformation}) describes the shear strain dissipation due to material element rearrangements, and the source term 
on the right-hand side of (\ref{eqn.heatflux}) describes the relaxation of the thermal impulse due to heat exchange between material elements.  

After the total energy potential has been specified, one can write all fluxes and source terms in an explicit form. Thus, for the energy $ E_2(\AAA,\JJ) $ given by~(\ref{eq:e_2}), 
we have $\boldsymbol{\psi}= E_{\AAA}= c_s^2 \AAA \dev(\GG)$, hence the shear stresses are 
\begin{equation}
\label{eqn.stress} 
\boldsymbol{\sigma}= -\rho\AAA^\mathsf{T} \boldsymbol{\psi} = -\rho\AAA^\mathsf{T} E_{\AAA} = -\rho c_s^2 \GG \dev(\GG),\ \ \qquad {\rm tr}(\boldsymbol{\sigma})=0,
\end{equation}
and the strain dissipation source term is 
\begin{equation}
\label{eqn.psi} 
-\dfrac{\boldsymbol{\psi}}{\theta_1(\tau_1)} = -\dfrac{E_{\AAA}}{\theta_1(\tau_1)}=-\dfrac{3}{\tau_1 } \left| \AAA \right|^{\frac{5}{3}} \AAA \dev(\GG),
\end{equation}
where we have chosen $ \theta_1(\tau_1) = \tau_1 c_s^2 / 3 \, |\AAA|^{-\frac{5}{3}} $, with $|\AAA|=\det(\AAA) > 0$ the determinant of $\AAA$ and $\tau_1$ being 
the strain relaxation time, or, in other words, the time scale that characterizes how long a material element is connected with its neighbor elements before 
rearrangement.\footnote{Following Frenkel~\cite{Frenkel1955},  this relaxation time was called particle-settled-life (PSL) time in \cite{PeshRom2014}.
}
Note, that the determinant of $\AAA$ must 
satisfy the \textit{constraint}    
\begin{equation} 
\label{eqn.compatibility} 
 |\AAA| = \frac{\rho}{\rho_0},    
\end{equation} 
where $\rho_0$ is the density at a reference configuration, see \cite{PeshRom2014}. 
Furthermore, from the energy potential $E_2(\AAA,\JJ)$ the heat flux vector follows with $E_{\JJ} = \alpha^2 \JJ$ directly as 
\begin{equation}
\label{eqn.hyp.heatflux}
 \mathbf{q} = T \, \mathbf{H} = E_s E_{\JJ} = \alpha^2 T \JJ.    
\end{equation} 
For the thermal impulse relaxation source term, we choose $ \theta_2=\tau_2\alpha^2 \frac{\rho}{\rho_0} \frac{T_0}{T}$, and hence
\begin{equation}
-\dfrac{\rho \mathbf{H}}{\theta_2(\tau_2)} = -\dfrac{\rho E_{\JJ}}{\theta_2(\tau_2)}= - \frac{T}{T_0} \frac{\rho_0}{\rho} \dfrac{\rho\JJ}{\tau_2}.
\end{equation}
It contains another characteristic relaxation time $\tau_2$ that is associated to heat conduction.  

The motivation for this particular choice of $\theta_1$ and $\theta_2$ can be found later in Section \ref{sec:asy}, where a formal asymptotic analysis of the model is presented, 
and where the connection with classical Navier-Stokes-Fourier theory is established in the stiff limit $\tau_1 \to 0$ and $\tau_2 \to 0$. 

\subsection{Discussion} 
\label{sec:discussion}

In this section, we discuss a few additional important properties of the HPR model. 
We first illustrate the relation of the HPR model to the laws of thermodynamics and the important role  played by the 
total energy potential. In particular, we demonstrate that the HPR model is compatible with the first and second law of thermodynamics, and that this 
automatically implies that the HPR model is a \textit{hyperbolic} system of PDEs, \textit{i.e.} the Cauchy problem for the system~\eqref{eqn.HPR} is well-posed. 
We complete this section by unveiling the characteristic structure of the HPR model. 
%

\subsubsection{Thermodynamically compatible systems of hyperbolic conservation laws and well-posedness}\label{sec.therm.hyperb}

\paragraph{\bf Overdetermined system of PDEs and the first law of thermodynamics} 
As many  other models of continuum mechanics, the system~\eqref{eqn.conti}-\eqref{eqn.energy}
is an \textit{overdetermined} system of PDEs. It consists of 18 PDEs for just 17
unknowns, and hence the natural question arises of whether it is \textit{consistent}, \textit{i.e.} whether it has at least one solution satisfying all the PDEs.
This is in general not guaranteed and one needs to provide evidences 
that a solution satisfying all the PDEs of the system does exist. 

In 1961, after discovering the mutual relations between thermodynamics, well-posedness of the initial value problem for systems of conservation laws and stability of numerical schemes, Godunov~\cite{God1961,God1962} concluded  that  an overdetermined system of conservation laws representing a continuum
mechanics model is consistent if it is compatible with the first law of thermodynamics, i.e. with the total energy conservation. 
In order to illustrate Godunov's idea, let us consider equations~\eqref{eqn.conti}-\eqref{eqn.energy} and let's also assume that it is an abstract system of PDEs, not 
necessarily related to the subject of this paper.  Following 
Godunov~\cite{God1961,God1962,God1972MHD}, we now show that if the unknown function $ E(t,\xx) $ is in fact not an unknown but a \textit{potential}, depending on all other 
unknowns, \textit{i.e.} $ E=E(\rho,\vv,\AAA,\JJ,s) $, then, if a solution of system~\eqref{eqn.HPR} exists, it also satisfies equation~\eqref{eqn.energy}, 
\textit{i.e.} the system \eqref{eqn.conti}-\eqref{eqn.energy} is consistent. 
In fact, we have to use the so-called \textit{conservative variables}, \textit{i.e.} we should consider the potential $ \rho E $ as a function of $ \rho $, $
\rho\vv $, $ \AAA $, $ \rho\JJ $, $ \rho s $. After this remark, one can see that equation~\eqref{eqn.energy} can be obtained as a linear combination of
equations~\eqref{eqn.HPR} multiplied by the factors\footnote{We recall that these factors should be understood as the partial derivatives, \textit{e.g.} $ (\rho E)_{\rho v_i} =\partial(\rho E)/\partial(\rho v_i)$.} $ E-V E_V-sE_s-v_iE_{v_i}-J_i E_{J_i} $, $ (\rho E)_{\rho v_i} $, $ (\rho E)_{A_{ik}} $, $ (\rho E)_{\rho J_i} $, and $ (\rho E)_{\rho s} $, \textit{i.e.}
\begin{equation}\label{eqn.summation}
(E-V E_V-sE_s-v_iE_{v_i}-J_i E_{J_i})\cdot\eqref{eqn.conti}+(\rho E)_{\rho v_i}\cdot\eqref{eqn.momentum} +(\rho E)_{A_{ik}}\cdot\eqref{eqn.deformation}+(\rho E)_{\rho J_i}\cdot\eqref{eqn.heatflux}+ (\rho E)_{\rho s}\cdot\eqref{eqn.entropy}\equiv\eqref{eqn.energy}.
\end{equation}
Here, the notation $ V=\rho^{-1} $ was used.  Because of the Gibbs identity
\begin{equation}
{\rm d}(\rho E)\equiv (E-V E_V-sE_s-v_iE_{v_i}-J_i E_{J_i}){\rm d}\rho+(\rho E)_{\rho v_i}{\rm d}\rho v_i +(\rho E)_{A_{ik}}{\rm d}A_{ik}+(\rho E)_{\rho J_i}{\rm d}\rho J_i+ (\rho E)_{\rho s}{\rm d}\rho s,
\end{equation}
it is obvious that~\eqref{eqn.summation} indeed holds for the time derivatives, as well as it holds for the right-hand sides, but it is less obvious that it is true for the space derivatives. In fact, the constitutive terms in the fluxes, \textit{i.e.} $ \rho^2E_\rho $, $ \rho A_{mi}E_{mk} $, $ E_{s} $, and $ E_{J_k} $ are chosen in these forms on purpose, because otherwise  it is impossible to get fully conservative fluxes in the energy conservation, but some non-conservative products would appear (details can be found in~\cite{God1972MHD,God1978,GodRom2003}, see also appendix in~\cite{Pesh2015}), which apparently violates  energy conservation.

Thus, identity~\eqref{eqn.summation} shows that if equations~\eqref{eqn.HPR} are fulfilled, then equation \eqref{eqn.energy}
 is also automatically  fulfilled. We stress once more that 
in order to have the property that the \textit{overdetermined} system~\eqref{eqn.conti}-\eqref{eqn.energy} of 18 PDEs for 17 unknowns is a \textit{consistent} system, 
the following constraints should hold 
\begin{itemize}
\item the function $ E(t,\xx) $ is \textit{not} an unknown but rather  a \textit{potential}, depending on the remaining unknowns, i.e. $ E=E(\rho,\vv,\AAA,\JJ,s)$;

\item all the constitutive terms in the fluxes and the dissipative source terms of the HPR model \eqref{eqn.HPR} are directly generated by the total energy potential 
by means of its gradients $ E_\rho $, $ E_{A_{ij}} $, $ E_{J_i} $, $ E_s $ in the state space and they must have this particular form in order to guarantee total 
energy conservation. 
\end{itemize}

In other words, these two requirements form the \textit{closure} for the overdetermined system \eqref{eqn.conti}-\eqref{eqn.energy},  making it consistent.

\paragraph{\bf Well-posedness of the Cauchy problem} 
It is not sufficient to propose a new continuum model that respects only some fundamental physical 
principles, but it is also required that the Cauchy problem for the proposed system of governing PDEs be \textit{well-posed}, \textit{i.e.} the solution of the system 
with initial data at time $t=0$ exists, at least locally, it is unique and stable. Otherwise, the practical value of the model would be questionable. 
In this context, hyperbolic conservation laws are very desirable for modeling dynamical phenomena, because hyperbolicity implies that the model 
is \textit{causal} (finite speed of perturbation propagation) and that the Cauchy problem for the nonlinear PDE system under consideration is 
\textit{well-posed} (hence, suitable for numerical treatment), see \textit{e.g.} see~\cite{GodUrMathPhys,Godlewski:1996a,Dafermos2005,KulikPogorSemen}.

From the discussion of the previous paragraph, 
it is obvious that the total energy potential plays a central role in the formulation of the HPR model. Moreover, we shall demonstrate that the convexity of the energy potential also guarantees that system~\eqref{eqn.HPR} is \textit{symmetric hyperbolic}, \textit{i.e.} the initial value problem for~\eqref{eqn.HPR} is well-posed. 

As noted by Godunov~\cite{God1961,God1962,God1972MHD},  an interesting parametrization of ovedetermined systems of conservation laws is possible. 
This parametrization allows to rewrite the original system in a symmetric quasilinear form. 
If, in addition, the total energy $ E $ is a \textit{convex function} of the state variables, 
then the system is symmetric hyperbolic. 
After a careful analysis of a large number of models in continuum mechanics, the original observation of Godunov was later extended to a wide class of thermodynamically consistent systems of hyperbolic conservation laws in a series of papers~\cite{GodRom1995,GodRom1996,godunov1996systems,GodRom1998,Rom1998,Rom2001} by Godunov and Romenski. All models belonging to this class of conservation laws are automatically symmetric hyperbolic. In particular, the system~\eqref{eqn.conti}-\eqref{eqn.energy} belongs to this class, see~\cite{godunov1996systems,Rom1998,Rom2001}. Therefore, in order to demonstrate that system~\eqref{eqn.HPR} is symmetric hyperbolic,  
we introduce the so-called \textit{thermodynamically conjugate}, or dual, state variables, which are in fact the factors in~\eqref{eqn.summation}:
\begin{equation}\label{eqn.conjstate}
r=E-V E_V-sE_s-v_iE_{v_i}-J_i E_{J_i},\ \ \nu_i=(\rho E)_{\rho v_i},\ \ \alpha_{ik}=(\rho E)_{A_{ik}},\ \ \Theta_i=(\rho E)_{\rho J_i},\ \ \sigma=(\rho E)_{\rho s},
\end{equation}
and the new thermodynamic potential $ L $ as the Legendre transform of $ \rho E $, \textit{i.e.}
\begin{equation}
L(r,\nu_i,\alpha_{ik},\Theta_i,\sigma)=r\rho+\nu_i \rho v_i+\alpha_{ik}A_{ik}+\Theta_i \rho J_i+\sigma\rho s-\rho E=\rho^2E_\rho + \rho A_{ij}E_{A_{ij}}.
\end{equation}
Now, the left hand side of~\eqref{eqn.HPR} can be rewritten as follows\footnote{We restrict the demonstration by considering only the left-hand side of~\eqref{eqn.HPR} because the type of a system of PDEs is defined by the leading terms.} (details can be found in~\cite{GodRom1995,GodRom1996,GodRom1998,Rom1998,Rom2001,Pesh2015})
\begin{subequations}\label{eqn.HPR.symmetric}
	\begin{align}
	& \frac{\partial L_r}{\partial t}+\frac{\partial (\nu_k L)_r}{\partial x_k}=0,\label{eqn.conti.sym}\\[2mm]
	&\displaystyle\frac{\partial L_{\nu_i}}{\partial t}+\frac{\partial (\nu_k L)_{\nu_i}}{\partial x_k}+L_{\alpha_{im}}\dfrac{\partial\alpha_{km}}{\partial x_k}-L_{\alpha_{mk}}\dfrac{\partial\alpha_{mk}}{\partial x_i}=0, \label{eqn.momentum.sym}\\[2mm]
	&\displaystyle\frac{\partial L_{\alpha_{il}}}{\partial t}+\frac{\partial (\nu_k L)_{\alpha_{il}}}{\partial x_k}+L_{\alpha_{ml}}\frac{\partial \nu_m}{\partial x_i}-L_{\alpha_{il}}\frac{\partial \nu_k}{\partial x_k}=0,\label{eqn.deformation.sym}\\[2mm]
	&\displaystyle\frac{\partial L_{\Theta_i}}{\partial t}+\frac{\partial ( \nu_k L)_{\Theta_i}}{\partial x_k}
	+\frac{\partial \sigma\delta_{ik}}{\partial x_k}=0, \label{eqn.heatflux.sym}\\[2mm]
	&\displaystyle\frac{\partial L_\sigma}{\partial t}+\frac{\partial (\nu_k L)_\sigma}{\partial x_k} + \frac{\partial \Theta_k}{\partial x_k}=0, \label{eqn.entropy.sym}
	\end{align}
\end{subequations}
and then in the quasilinear form
\begin{equation}\label{eqn.symm}
\mathsf{A}(\PP)\dfrac{\partial\PP}{\partial t} + \mathsf{B}_k(\PP)\dfrac{\partial \PP}{\partial x_k}=0, 
\end{equation}
where $ \PP=(r,\nu_i,\alpha_{ik},\Theta_i,\sigma) $, and matrices $ \mathsf{A}^\mathsf{T}=\mathsf{A}$ and $ \mathsf{B}^\mathsf{T} =\mathsf{B}$ are symmetric, and moreover $ \mathsf{A}>0 $ if the potential $ L(r,\nu_i,\alpha_{ik},\Theta_i,\sigma) $ is a convex function. We recall, that because of the properties of the Legendre transformation, the convexity of $ L(r,\nu_i,\alpha_{ik},\Theta_i,\sigma) $ is equivalent to the convexity of $ \rho E $ with respect to the conservative variable. In other words, the system~\eqref{eqn.symm}, as well as \eqref{eqn.HPR}, is \textit{symmetric hyperbolic} if $ \rho E $ is a convex potential of the conservative state variables, and the solution to the initial value problem exists locally. In turn, we note that via a direct calculation one can verify that the convexity of $ \rho E $ with respect to the conservative state variables is equivalent to the convexity of $ E $ with respect to the primitive state variables $ \rho $, $ v_i $, $ A_{ij} $, $ J_i $ and $ s $.

\paragraph{\bf Energy transformation and the second law} 
The energy is the only quantity that is allowed to be transferred among all the  three scales involved, namely the micro-, meso-, and macroscales. 
Therefore, the scales can interact only through an  energy exchange, and the total energy potential has to be involved in some way in the mathematical formulation of this interaction. Indeed, the energy transfer from meso- to macroscale, $ E_2(\AAA,\JJ)\rightarrow E_3(\vv) $, is known as \textit{reversible} energy transformation, and is controlled by the momentum fluxes, and as we have seen in the previous paragraph, these fluxes are given by the gradients $ E_\rho $ and $ E_{A_{ij}} $. The energy transfer from meso- to microscale, $ E_1(\rho,s)\leftarrow E_2(\AAA,\JJ) $, is an  \textit{irreversible} transformation, which is controlled by the dissipative source terms in the governing equations for the distortion tensor, the thermal impulse, and the entropy. 
Thus, it is natural to expect that these dissipative source terms in the HPR model are also generated by the energy potential, via its partial derivatives with respect to the state variables.  
Indeed, the total energy conservation principle holds regardless of whether dissipation is present, or not. Thus, even if the dissipative source terms are present, we anyway have to 
have zero on the right-hand side of the total energy conservation law. 
Since the HPR model is an overdetermined system of PDEs, we require that the summation identity \eqref{eqn.summation} holds. Hence, each dissipative source term is multiplied 
by the corresponding factor (conjugate state variables~\eqref{eqn.conjstate}) and the sum must vanish. Let us denote the source terms in equations \eqref{eqn.deformation} and \eqref{eqn.heatflux} 
by $ S^{A}_{\! ik} = -\psi_{ik} / \theta_1 $ and $ S^{\! \rho J}_{\! i} = - \rho H_i/\theta_2 $, respectively, while the source term in the entropy equation \eqref{eqn.entropy} 
is denoted by $S^{\! \rho s}$. In the summation identity~\eqref{eqn.summation}, they are 
multiplied by $ (\rho E)_{A_{ik}} =\rho E_{A_{ik}}$, $ (\rho E)_{\rho J_i} =E_{J_i}$ and $ (\rho E)_{\rho s}=E_s =T$, respectively. 
Total energy conservation requires that the right hand side of \eqref{eqn.summation} vanishes, i.e. 
\begin{equation} 
   \label{eqn.summation.rhs} 
   \rho E_{A_{ik}} S^{A}_{\! ik} + E_{J_i} S^{\! \rho J}_{\! i} + E_s \, S^{\! \rho s} = 0. 
\end{equation}  
The only freedom we have to satisfy the total energy conservation law is to set  
\begin{equation}\label{eqn.entropyprod}
 S^{\! \rho s} = -\dfrac{1}{E_s} \left( \rho E_{A_{ik}} S^{A}_{\! ik} + E_{J_i} S^{\! \rho J}_{\! i}  \right) 
 = \dfrac{1}{E_s} \left( \rho E_{A_{ik}} \frac{\psi_{ik}}{\theta_1} + E_{J_i} \frac{\rho H_i}{\theta_2} \right).
\end{equation}

At that point, we recall that the thermodynamics of dissipative processes requires that the entropy cannot decrease, and hence the entropy production \eqref{eqn.entropyprod} has to be  nonnegative. A simple possibility to guarantee this is to assume that the terms $ \psi_{ik} $ and $ H_{i} $ are proportional to the gradients $ E_{A_{ik}} $ 
and $ E_{J_i} $, respectively, with some positive coefficients. This makes the entropy source term a positive definite quadratic form, which guarantees that the entropy does not 
decrease. Since the functions $\theta_1 >0$ and $\theta_2 > 0$ are positive, in this paper we have simply chosen 
\begin{equation}
\psi_{ik} = E_{A_{ik}} \qquad \textnormal{ and } \qquad H_i = E_{J_i}.  
\end{equation}


\subsubsection{Characteristic speeds and sound speeds}\label{sec.charcteristics}

Understanding the \textit{characteristic} structure of a hyperbolic system is an important step in studying the solution properties, because the solution of a hyperbolic model 
is a combination of waves propagating along the \textit{characteristic lines}, \textit{e.g.} see~\cite{GodUrMathPhys,Dafermos2005}. In this section we study the characteristic 
structure of the HPR model.  First we shall present the characteristic structure of the viscous part of the HPR model (equations~(\ref{eqn.conti}), (\ref{eqn.momentum}), (\ref{eqn.deformation}) and (\ref{eqn.energy})), then we discuss the characteristic structure of the heat conducting part (equations \eqref{eqn.heatflux} and \eqref{eqn.entropy}), and eventually we close this section by presenting the structure of the entire model~(\ref{eqn.HPR}). It is also important to recall that the characteristic speeds of a hyperbolic model with  stiff dissipative source terms are not the true sound speeds in the media, because these apparent sound speeds are strongly influenced by the dissipative processes giving rise to the  phenomena called sound dispersion. In fact, the characteristic speeds of a hyperbolic system with stiff dissipative source terms are the high frequency limits for sound speeds~\cite{Rug1992}.

We also note that the HPR model is fundamentally different from the classical parabolic NSF theory in the way it treats viscous and heat conducting phenomena. In the classical NSF theory, the transport phenomena are treated by means of phenomenological transport relations such as, for example, Newton's law of viscosity and Fourier's law of heat conduction, while in the HPR model all transport phenomena are treated from the wave propagation point of view. Thus, as it will be shown below, there are four types of sound waves in the HPR model, one for the 
transport of longitudinal (or pressure) perturbations, two for shear perturbations, and one for heat transfer, in contrast to only one pressure wave in the NSF model.


\paragraph{\bf Viscous subsystem} Let us consider system \eqref{eqn.HPR} in the one-dimensional case.
If a direction $ x_k $ is chosen, then in terms of the state variables
\begin{equation}
\QQ=(\rho,p,v_1,v_2,v_3,A_{1k},A_{2k},A_{3k})^\mathsf{T},
\end{equation}
PDEs \eqref{eqn.conti}, \eqref{eqn.energy}, \eqref{eqn.momentum} and \eqref{eqn.deformation}  can be written in the quasilinear form
\begin{equation}\label{eq:quasilinear}
\dfrac{\partial \QQ}{\partial t} + \mathsf{A}_k(\QQ)\dfrac{\partial\QQ}{\partial x_k}=\mathsf{S}_k(\QQ)
\end{equation}
with the source vector 
\[\mathsf{S}_k(\QQ)=(0,0,0,0,0,-E_{A_{1k}}/\theta_1,-E_{A_{2k}}/\theta_1,-E_{A_{3k}}/\theta_1)^\mathsf{T},\]
and matrix $ \mathsf{A}_k(\QQ) $ given in \ref{ap.matrix} as well as the formulas of the eigenvalues for $ \mathsf{A}_k $.
The full basis consisting of eigenvectors of the matrix $ \mathsf{A}_k $ can be also obtained in the same way as it is done in \cite{Barton2009,BartonRom2010}. 
We note that in~\cite{Barton2009,BartonRom2010}, the time evolution equation for $ \AAA^{-1} $ was used instead of equation~\eqref{eqn.deformation}.

In order to illustrate the characteristic structure of \eqref{eq:quasilinear} we restrict ourselves to the consideration of a fluid (or solid) at the rest state $ \QQ_0 $, \textit{i.e.} $  \vv=0 $, $ \AAA=\II $, $ \rho=\rho_0 $. If the internal energy $ E_1(\rho,s) $ is considered in the form~\eqref{eq:ideal_gas_eos}, or \eqref{eq:stiff_gas_eos}, and $ E_2(\AAA) $ in the form \eqref{eq:e_2} then  matrix $ \mathsf{A}_k $, $ k=1 $, looks as follows (we omit the subscript 1 in what follows)
\[
\mathsf{A}=\left (\begin{array}{cccccccc}
	0 &     0     & \rho  & 0 & 0 &        0         &   0   &   0   \\
	0 &     0     & c_0^2 & 0 & 0 &        0         &   0   &   0   \\
	0 & \rho^{-1} &   0   & 0 & 0 & \frac{4}{3}c_s^2 &   0   &   0   \\
	0 &     0     &   0   & 0 & 0 &        0         & c_s^2 &   0   \\
	0 &     0     &   0   & 0 & 0 &        0         &   0   & c_s^2 \\
	0 &     0     &   1   & 0 & 0 &        0         &   0   &   0   \\
	0 &     0     &   0   & 1 & 0 &        0         &   0   &   0   \\
	0 &     0     &   0   & 0 & 1 &        0         &   0   &   0
\end{array}\right ),
\]
Its non-zero eigenvalues are $ \lambda_{1,2,3,4}=\pm c_s $, $ \lambda_{5,6}=\pm \sqrt{c_0^2+\frac{4}{3}c_s^2} $. Thus, the eigenvalues $ \lambda_{1,2,3,4} $ are (in general, they are distinct if shear flow is present) 
transverse, or shear, characteristic speeds, while $ \lambda_{5,6} $ is the longitudinal characteristic speed. In the framework of solid mechanics,
the existence of  two types of waves comes with no surprise, 
but this is much less obvious for fluid mechanics. However, this fact is in full agreement with the causality principle and the wave propagation 
point of view on the transport phenomena, as it was already mentioned above.

Yet another point has to be explained. One may note that the characteristic velocity corresponding to the propagation of pressure perturbations in fluids modeled by EOS (\ref{eq:ideal_gas_eos}) or (\ref{eq:stiff_gas_eos}) is $ c_0 $, while we get $ \sqrt{c_0^2+\frac{4}{3}c_s^2} \neq c_0$ for $ c_s\neq 0 $. In fact, this is not a paradox. It is necessary to recall that for the hyperbolic PDEs with
stiff dissipative source terms like system~\eqref{eq:quasilinear}, the characteristic speeds are not the \textit{ true sound speeds}, but the true sound speeds are the result of a 
coupling of the non-dissipative waves modeled by the left hand side of~\eqref{eq:quasilinear} and the dissipative processes modeled by the algebraic source terms, and therefore the true sound speeds can be obtained only via a \textit{dispersion analysis}. However, such an analysis is outside the scopes of this paper, some details of it are given just to demonstrate that there is no controversy between the sound speeds predicted by the HPR model and experimental observations on sound propagation in fluids.

The dispersion relation for a hyperbolic system of PDEs of the form~(\ref{eq:quasilinear}) with algebraic dissipative source terms is~\cite{Rug1992}
\begin{equation}\label{eq:dispers_rel}
\det \left ( \II -z \mathsf{A} + \dfrac{i}{\omega} \mathsf{B}\right )=0,
\end{equation}
where $\omega=2\pi f$ is the angular frequency, $ f $ is the wave frequency, $z=k/\omega$, $ k $ is the complex wave number, $i$ is the imaginary unit, matrices $ \mathsf{A}=\mathsf{A}(\QQ_0) $ and $ \mathsf{B}(\QQ_0)=\partial\mathsf{S}/\partial \QQ $ are taken at the rest state $ \QQ_0 $, and $ \II $ is the identity matrix of the same size as $\mathsf{A}$ and $\mathsf{B}$. 
Once the solution $z$ to~(\ref{eq:dispers_rel}) is found, the phase velocity, $V$, and the attenuation factor, $\alpha$, of a harmonic sound wave of frequency $ \omega $ are given by
\begin{equation}\label{eq:vel_attenuation}
V=\dfrac{1}{{\rm Re}(z)}, \ \ \ \mu=-\omega\, {\rm Im}(z).
\end{equation}

Equation~(\ref{eq:dispers_rel}) has six nontrivial solutions, four corresponding to the transverse waves
\begin{equation}
z_{1,2}=-\sqrt{\dfrac{\Omega - 3 i}{\Omega\,c_s^2}},\ \  z_{3,4}=\sqrt{\dfrac{\Omega - 3 i}{\Omega\,c_s^2}},
\end{equation}
and two corresponding to the longitudinal waves
\begin{equation}
z_{5,6}=\pm\sqrt{\dfrac{3(\Omega-2 i)}{3 c_0^2 (\Omega - 2 i) + 4\Omega\, c_s^2}}\ ,
\end{equation}
where $\Omega=\tau_1\omega$. Thus, for the \textit{longitudinal sound waves} the phase velocity and the attenuation factor are
\begin{align}\label{eq:phase_vel_l}
V^{long}=\frac{c_\infty}{2 X}\left(\sqrt{(X-\Omega)(Y-\Omega)}+\sqrt{(X+\Omega)(Y+\Omega)}\right), \\
\mu^{long}=\frac{\omega}{2 c_\infty Y}\left(\sqrt{(X-\Omega)(Y+\Omega)}-\sqrt{(X+\Omega)(Y-\Omega)}\right),\\
c_\infty=\sqrt{c_0^2+\frac{4 }{3}c_s^2},\ \ X=\sqrt{\Omega^2+16},\ \ Y=\sqrt{\Omega^2+16\left(\frac{c_0}{c_\infty}\right) ^4}\nonumber,
\end{align}
and for \textit{shear sound waves} they are
\begin{equation}\label{eq:phase_vel_t}
V^{shear}=c_t \sqrt{\dfrac{\Omega(Z+\Omega)}{2 Z^2}}\ ,\ \ 
\mu^{shear}=\dfrac{\omega}{c_s\sqrt{2\Omega}}\sqrt{Z-\Omega},\ \ Z=\sqrt{\Omega^2+9}.
\end{equation}

\begin{figure}[!htbp]
	\begin{center}
		\begin{tabular}{cc} 
			\includegraphics[trim = 0mm 0mm 0mm 0mm, clip, width=0.405\textwidth]{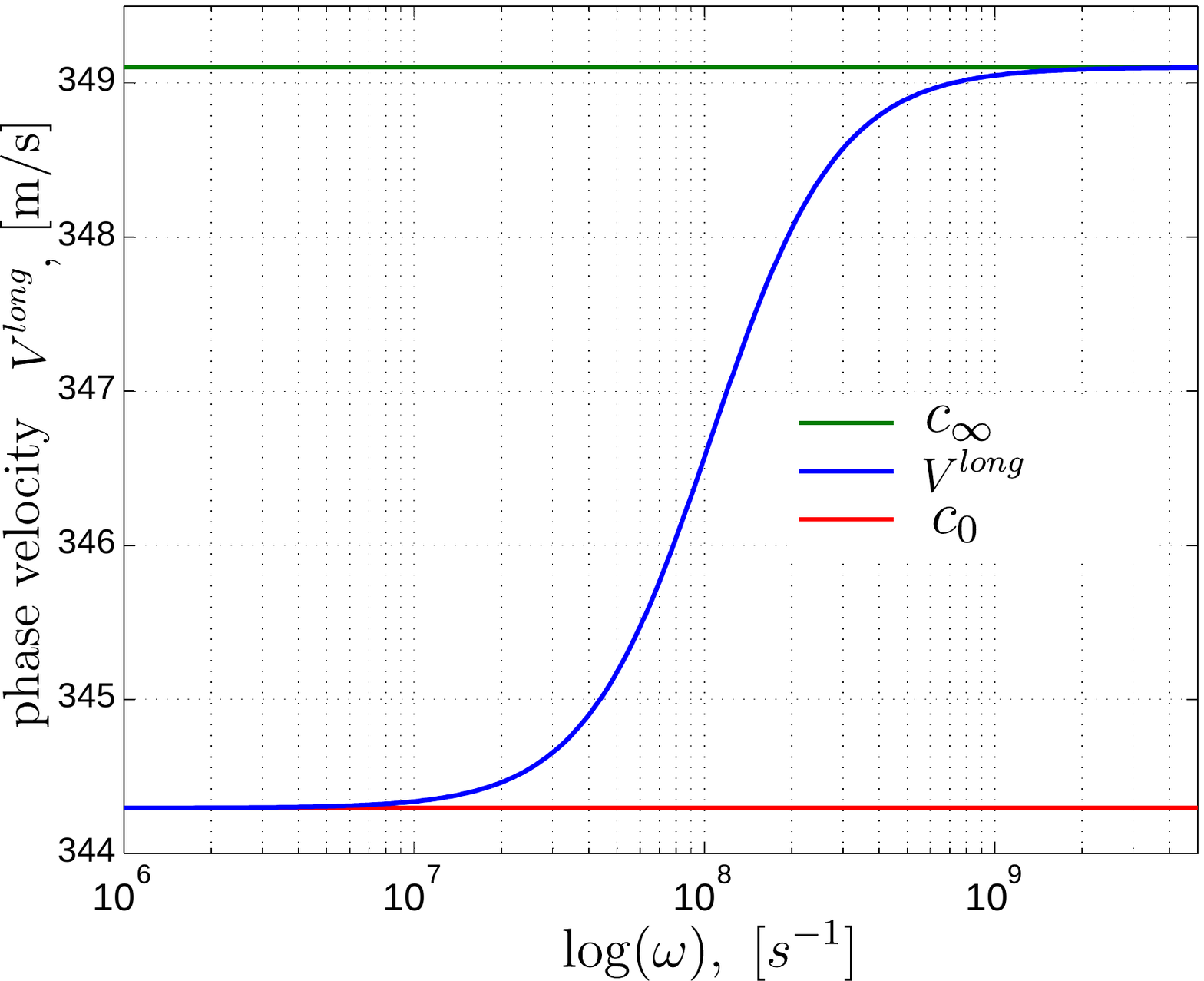} & 
			\includegraphics[trim = 0mm 0mm 0mm 0mm, clip, width=0.40\textwidth]{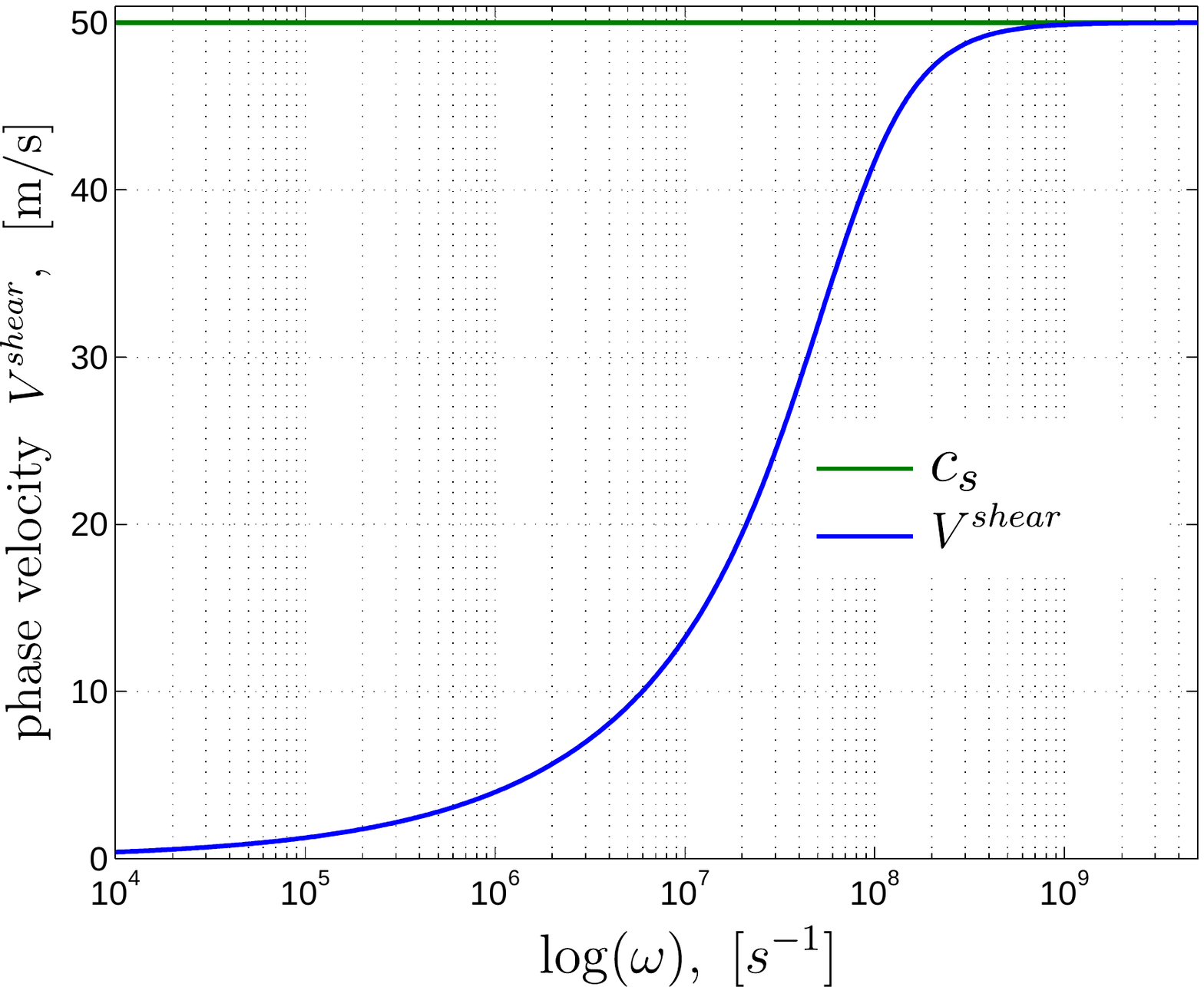} 
		\end{tabular} 
		\caption{Phase velocity of the longitudinal wave (left) and shear wave (right) versus $ \log(\omega) $ propagating in a viscous gas with parameters $ \rho=1.177\ {\rm kg/m^3}$, $ \gamma=1.4 $, $ c_v=718\ {\rm J/(kg\,K)} $, $ s=8100 $, $ c_0=344.3$ m/s, $ c_s =50$ m/s, $ \mu= 1.846\cdot10^{-5}$ Pa$ \cdot $s, $ \tau_1= 3.76\cdot10^{-8}$~s. } 
		\label{fig.dispers}
	\end{center}
\end{figure}

By a direct verification, one can see that the low ($ \Omega\rightarrow 0 $) and high ($ \Omega\rightarrow\infty $) frequency limits of $ V^{long} $ are $ c_0 $ and $ c_\infty $, accordingly, while for $ V^{shear} $ they are $ 0 $ and $ c_s $. This clearly indicates that (i) perturbations of any frequency propagate at finite speeds in contrast to the classical NSF theory and that (ii) the low frequency sound waves propagate at velocities $ \approx c_0 $ what we in fact use to call the sound speed in fluids. Fig.~\ref{fig.dispers} shows the longitudinal and shear sound speeds as a function of the angular frequency $ \omega $. 

%

\paragraph{\bf Heat conducting subsystem}

For convenience, we rewrite the heat conduction equations~(\ref{eqn.heatflux}) and \eqref{eqn.entropy} in the form
\begin{subequations}\label{eqn.heatsubsys}
	\begin{align}
	&\rho\displaystyle\frac{{\rm d} J_k}{{\rm d} t}+\frac{\partial E_s}{\partial x_k}=-\dfrac{\rho E_{J_k}}{\theta_2}, \label{eqn.heatflux2}\\[2mm]
	&\rho\displaystyle\frac{{\rm d }s}{{\rm d} t}+\frac{\partial E_{J_k}}{\partial x_k}=\dfrac{\rho}{\theta_2 E_s} E_{J_i} E_{J_i}\geq0, \label{eqn.entropy2}
	\end{align}
\end{subequations}
where $ {\rm d}/{\rm d}t =\partial/\partial t+ v_k\partial/\partial x_k$ is the material time derivative, and the energy potential $ E $ is taken to be $ E=E_1(\rho,s) + \alpha^2 J_k J_k/2$, while the ideal gas EOS~\eqref{eq:ideal_gas_eos} is used for $ E_1(\rho,s) $. Let us consider this system in the direction $ x_1 $, then it can be rewritten in a quasilinear form
\begin{equation}\label{eqn.heatquasi}
\dfrac{{\rm d}}{{\rm d} t}\left (\begin{array}{c}
J_1\\
s
\end{array}\right )+\left ( \begin{array}{cc}
0 & \frac{T}{\rho c_V}\\
\frac{\alpha^2}{\rho} & 0
\end{array}\right )
\dfrac{\partial}{\partial x_1}\left (\begin{array}{c}
J_1\\
s
\end{array}\right )
=
\frac{\rho_0}{\rho T_0 \tau_2} \left (\begin{array}{c}
- T  J_1  \\
+ \alpha^2 J_1^2   
\end{array}\right ).
\end{equation}
The eigenvalues of the homogeneous part of the system \eqref{eqn.heatquasi} are $ \lambda_{1,2}=\mp \alpha\sqrt{T}/(\rho\sqrt{c_V})$, and in the following we shall use the notation
\begin{equation}
c_h=\dfrac{\alpha}{\rho}\sqrt{\dfrac{T}{c_V}}
\end{equation}
for the velocity of the heat characteristic. The dispersion relation \eqref{eq:dispers_rel} for the heat conducting subsystem (\ref{eqn.heatquasi}) can also be treated analytically. Thus, the phase velocity for heat harmonic wave of angular frequency $ \omega $ is
\begin{equation}
V^{heat}=2\dfrac{\alpha\sqrt{T\Omega}}{\rho\sqrt{c_V}}\left(\dfrac{\sqrt{X+\rho}+\sqrt{X-\rho}}{\sqrt{X+\rho}-\sqrt{X-\rho}+2X+2\Omega}\right),\ \ \ X=\sqrt{\rho^2+\Omega^2},\ \ \ \Omega=\tau_2\omega.
\end{equation}
In particular, one can see that the low frequency limit ($ \Omega=\tau_2\omega \rightarrow 0$) and the high frequency limit ($ \Omega=\tau_2\omega \rightarrow \infty$) of the phase velocity are 0 and $ c_h $, respectively.

\paragraph{\bf The full system} We shall now consider the full HPR system and study its characteristic structure assuming that
 the space coordinate $ x_1 $ is chosen as the direction of wave propagation. We chose the following vector of state variables
\begin{equation}
\QQ=(\rho,p,J_1,v_1,v_2,v_3,A_{11},A_{21},A_{31})^\mathsf{T}\,.
\end{equation}
To discuss the characteristic structure, it is again sufficient to consider wave propagation near the rest state $ \QQ_0 $ characterized by $ \vv=0 $, $ \AAA=\II $, $ \JJ=0 $. If the ideal gas EOS is used for the $ E_1(\rho,s) $, the matrix $ \mathsf{A}(\QQ_0) $ reads as
\begin{equation}
 \mathsf{A}(\QQ_0)=\left(
\begin{array}{ccccccccc}
0 & 0 & 0 & \rho  & 0 & 0 & 0 & 0 & 0 \\
0 & 0 & \beta  c_h^2 & \rho  c_0^2 & 0 & 0 & 0 & 0 & 0 \\
-\frac{T}{\rho } & \beta^{-1} & 0 & 0 & 0 & 0 & 0 & 0 & 0 \\
0 & \rho^{-1} & 0 & 0 & 0 & 0 & \frac{4 }{3}c_s^2 & 0 & 0 \\
0 & 0 & 0 & 0 & 0 & 0 & 0 & c_s^2 & 0 \\
0 & 0 & 0 & 0 & 0 & 0 & 0 & 0 & c_s^2 \\
0 & 0 & 0 & 1 & 0 & 0 & 0 & 0 & 0 \\
0 & 0 & 0 & 0 & 1 & 0 & 0 & 0 & 0 \\
0 & 0 & 0 & 0 & 0 & 1 & 0 & 0 & 0
\end{array}
\right),
\end{equation}
where $ \beta= c_V(\gamma-1)\rho^2$, $ T=E_s $. This matrix has eight non-zero eigenvalues; four eigenvalues $ \lambda_{1,2,3,4} =\mp c_s$ corresponding to two shear waves, two eigenvalues
\begin{equation}
\lambda_{5,6}=\mp\dfrac{1}{\sqrt{2}}\sqrt{C-\dfrac{\sqrt{C^2-4\rho c_h^2(\beta T+\frac{4}{3}\rho c_s^2)}}{\rho}},\ \ \ C=c_0^2+c_h^2+\frac{4}{3}c_s^2
\end{equation}
corresponding to heat waves,
and two eigenvalues
\begin{equation}
\lambda_{7,8}=\mp\dfrac{1}{\sqrt{2}}\sqrt{C+\dfrac{\sqrt{C^2-4\rho c_h^2(\beta T+\frac{4}{3}\rho c_s^2)}}{\rho}},\ \ \ C=c_0^2+c_h^2+\frac{4}{3}c_s^2
\end{equation}
corresponding to longitudinal pressure waves. The same dispersion analysis as above can be performed for the viscous heat conducting case, but we 
do not enter such details here, as they would distract us from the main purpose of the present work.

\subsection{Formal asymptotic analysis, Newton's viscous law and Fourier's law of heat conduction} 
\label{sec:asy}

In this section we show how to establish a link between the HPR model \eqref{eqn.conti}-\eqref{eqn.energy} and 
the classical Navier-Stokes-Fourier (NSF) theory in the stiff relaxation limit $\tau_1 \ll 1$ and $\tau_2 \ll 1$. 

\subsubsection{Asymptotic limit of the viscous stress tensor}

We first concentrate on the relaxation limit of the viscous stress tensor $\boldsymbol{\sigma}$. For that purpose, we can ignore the rotational degree of freedom contained in the 
distortion tensor $ \A $, since $\BS$ is only a function of the symmetric tensor $ \G=\A^\mathsf{T}\A $, which contains only the information about the deformation of the material elements. The temporal evolution equation of $\GG$ can be obtained from Eqn. \eqref{eqn.deformation} as \footnote{To obtain this PDE, it is necessary to sum up equation \eqref{eqn.deformation} multiplied by $ \A^\mathsf{T} $  from the left and transpose equation \eqref{eqn.deformation} multiplied by $ \A $ from the right, since $\dot{\G} = \A^\mathsf{T} \dot{\A} + \dot{\A}^\mathsf{T} \A$. We also use here that $ \BS=-\rho\A^\mathsf{T}E_\A=-\rho(E_\A)^\mathsf{T} \A =\BS^\mathsf{T}.$} 
\begin{equation}\label{eqn.GODE}
	\dot{\G} = - \left( \G \nabla \v + \nabla \v^\mathsf{T} \G \right) +\dfrac{2}{\rho\,\theta_1} \, \BS,
\end{equation} 
where $ \dot{\G} =\partial\G/\partial t + \v\cdot\nabla\G$ is the material time derivative of $\GG$ and $\nabla \v$ is the velocity gradient. 

We now proceed with a \textit{formal} asymptotic expansion\footnote{the so-called Chapman-Enskog expansion} of the tensor $\GG$ in a series of the small relaxation 
parameter $\tau_1$, 
\begin{equation}
  \GG = \GG_0 + \tau_1 \GG_1 + \tau_1^2 \GG_2 + ... 
	\label{eqn.Gseries} 
\end{equation} 
Furthermore, we have to specify $ \theta_1 $. In this paper, we choose $ \theta_1 = \tau_1 |\A|^{-\frac{5}{3}} c_s^2/3 = \tau_1 |\G|^{-\frac{5}{6}} c_s^2/3  $. Now, inserting \eqref{eqn.Gseries} into \eqref{eqn.GODE} and collecting terms of the same power in $\tau_1$ yields: 
\begin{eqnarray}
 \frac{ {\rm d} }{{\rm d} t}(\G_0+\tau_1 \G_1 + ...) = -\left ( (\G_0+\tau_1 \G_1 + ...)\nabla\v+\nabla\v^\mathsf{T}
 (\G_0+\tau_1 \G_1 + ...)\right )-\nonumber\\
 \dfrac{6}{\tau_1}|\G_0+\tau_1\G_1+...|^{\frac{5}{6}}(\G_0+\tau_1 \G_1 + ...)\dev(\G_0+\tau_1 \G_1 + ...),  
 \label{eqn.aux2}
\end{eqnarray}  
and then
\begin{equation}
 \tau_1^{-1} \, \underbrace{\left(6|\G_0|^{\frac{5}{6}}\G_0\dev(\G_0)\right )}_{0} + 
  \tau_1^0  \, \underbrace{\left( \frac{ {\rm d} \G_0}{{\rm d} t}  + ...\right )}_{0} + ... = 0.  
 \label{eqn.aux3}
\end{equation}

\paragraph{\bf Inviscid fluid as a zeroth order approximation}
Since \eqref{eqn.aux3} is valid for any $\tau_1$, the coefficients multiplying powers of $\tau_1$ must all vanish. Furthermore, density positivity $\rho =\rho_0|\A|=\rho_0|\G|^{\frac{1}{2}}> 0$ implies that $|\G| > 0$, hence $\G$ is invertible. Therefore, from the first term in \eqref{eqn.aux3}, it follows that
\begin{equation}
  \dev \left( \G_0 \right) = 0,  \quad \Rightarrow  \quad \G_0 - \frac{1}{3} \tr \left( \GG_0 \right) \Id = 0, \quad \Rightarrow 
	\G_0 = \frac{1}{3} \tr \left( \GG_0 \right) \Id. 
	\label{eqn.devg0} 
\end{equation} 
Setting $g := \frac{1}{3} \tr \left( \GG_0 \right)$ and neglecting higher order terms in $\tau_1$, we get the first important result:
\begin{equation}
  \G = g \Id + \tau_1 \G_1 = \AAA^{\mathsf{T}} \AAA, 
	\label{eqn.gtemp} 
\end{equation} 
which means that the distortion tensor $\AAA$ tends towards an \textit{orthogonal matrix} in the stiff limit $\tau_1 \ll 1$. 

To obtain the unknown coefficient $g$ from known quantities, we compute the determinant of $\G$ from \eqref{eqn.gtemp}, neglecting small terms of the order $\mathcal{O}(\tau_1)$. Hence,  
\begin{equation}
\label{eqn.detg} 
 | \G | = g^3 = |\AAA|^2, \qquad \Rightarrow \qquad g = | \G |^{\frac{1}{3}} = | \AAA |^{\frac{2}{3}} = \left( \frac{\rho}{\rho_0} \right)^{\frac{2}{3}}.
\end{equation} 
If we retain only the leading zeroth order term $\G_0$ of the expansion \eqref{eqn.Gseries}, then $ \G=\G_0=g\II $ and thus $ \sigma=-\rho c_s^2\G_0 \dev(\G_0)=0 $, \textit{i.e.} viscous stresses 
vanish and we retrieve the \textit{inviscid} case (compressible Euler equations) as a zeroth order approximation of the HPR model in the stiff limit $\tau_1 \ll 1$. 
The relations \eqref{eqn.gtemp} and \eqref{eqn.detg} above also imply that in the inviscid limit, the \textit{shape} of the material elements does not change, but only their \textit{volume}.

\paragraph{\bf Newton's viscous law as a first order approximation}
We are now interested in a first order approximation of the viscous stress tensor $\BS$ that results in the stiff relaxation limit $\tau_1 \ll 1$. For that purpose, we expand the stress tensor \eqref{eqn.stress} in a series of $ \tau_1 $. 
Here, we will use that $ \G=g\II+\tau_1\G_1 $, which results in $ \rho=\rho_0|\A|=\rho_0|g\II+\tau_1\G_1|^\frac{1}{2}=\rho_0 (g^{3/2}+\frac{\tau_1}{2} g^{1/2} \tr(\G_1)+\mathcal{O}(\tau_1^2))$ and 
$ \dev(\G)=\dev(g\II+\tau_1\G_1)=\tau_1\dev(\G_1)$: 
\begin{equation}
\BS = - \rho c_s^2 \G \dev(\G)=-\rho_0 c_s^2\left  (g^{3/2}+\frac{\tau_1}{2} g^{1/2} \tr(\G_1)\right )\left (g\II+\tau_1\G_1\right )\tau_1\dev(\G_1).
\label{eqn.sdot} 
\end{equation} 
Thus, ignoring higher order terms in $ \tau_1 $ yields   
\begin{equation}\label{eqn.sigma.expan}
\BS = - \tau_1  \rho_0 c_s^2 g^{5/2}\dev(\G_1).
\end{equation}

The evolution equation for $ \dev(\G) $, which is obtained by applying the ``$ \dev $'' operator to \eqref{eqn.GODE}, reads:  
\begin{equation}\label{eqn.devG}
\dfrac{{\rm d}}{{\rm d}t}\dev(\G)+\G\nabla \v+\nabla \v^\mathsf{T}\G-\dfrac{1}{3}\tr(\G\nabla \v+\nabla \v^\mathsf{T}\G)\II=-\dfrac{6}{\tau_1}|\G|^{5/6}\dev(\G\dev(\G)).
\end{equation}
Inserting the expansion \eqref{eqn.Gseries} into \eqref{eqn.devG}, collecting terms of the same power of $ \tau_1 $ and recalling from \eqref{eqn.devg0} that $\dev{\G_0}=0$, 
one gets for the leading order terms ($ \tau_1^0 $): 
\begin{equation*}
\G_0\nabla \v+\nabla \v^\mathsf{T}\G_0-\dfrac{2}{3}\tr(\G_0\nabla \v)\II=-6|\G_0|^{7/6}\dev(\G_1).
\end{equation*}
Since $ \G_0=g\II $, the last equality can be rewritten as
\begin{equation}\label{eqn.devG1}
g\left  (\nabla \v+\nabla \v^\mathsf{T}-\dfrac{2}{3}\tr(\nabla \v)\II\right  )=-6\, g^{7/2}\dev(\G_1).
\end{equation}
By inserting \eqref{eqn.devG1} into \eqref{eqn.sigma.expan}, we conclude that 
\begin{equation}
\BS = \frac{1}{6} \tau_1 \rho_0 c_s^2 \left( \nabla \v + \nabla \v^T  - \frac{2}{3} \tr ( \nabla \v ) \Id  \right) := \mu \left( \nabla \v + \nabla \v^T  - \frac{2}{3} ( \nabla \cdot \v ) \Id  \right),  
\label{eqn.stresslimit} 
\end{equation} 
which is the classical stress tensor known from the compressible Navier-Stokes equations based on Stokes' hypothesis, with the dynamic viscosity coefficient 
\begin{equation}
	\mu = \frac{1}{6} \tau_1 \rho_0 c_s^2,  
	\label{eqn.mu} 
\end{equation} 
as already given in \cite{PeshRom2014} \footnote{Please note that for obtaining a constant viscosity coefficient $\mu$ as given in \eqref{eqn.mu}, in \cite{PeshRom2014} there was a 
factor $|A|^{\frac{8}{3}}$ missing in front of the relaxation source term $\boldsymbol{\psi}$.}. This completes the formal asymptotic analysis, establishing a direct  
connection of the stress tensor $\BS$ of the HPR model with the known viscous stress tensor of the compressible Navier-Stokes equations, which is automatically recovered by the 
HPR model for small relaxation times $\tau_1 \ll 1$.  
At this point it is very important to highlight that in the HPR model, the viscous stress tensor $\BS$ obtained in the stiff relaxation limit \eqref{eqn.stresslimit} is a 
\textit{result} of the choice of a simple quadratic form for the contribution $E_2(\AAA,\JJ)$ to the total energy potential $E$, see Eqn. \eqref{eq:e_2}. 
In contrast, in classical Navier-Stokes theory, the stress tensor $\BS$ is postulated as a constitutive relation right from the beginning. 

Our particular choice of $\theta_1$ has been made only in order to obtain a \textit{constant} viscosity coefficient $\mu$. In order to obtain a variable viscosity coefficient that depends, 
for example, on the temperature, like in Sutherland's law, it is sufficient to modify the function $\theta_1(\tau_1)$ accordingly. 

\subsubsection{Asymptotic limit of the heat flux} 

Next, we proceed with a similar formal asymptotic analysis of the heat flux $\mathbf{q} = \alpha^2 T \JJ$, which is, however, much simpler than the 
previous analysis of the stress tensor $\BS$. We recall the governing PDE \eqref{eqn.heatflux} for the vector $\JJ$ with the choice
$ \theta_2=\tau_2\alpha^2 \frac{\rho}{\rho_0} \frac{T_0}{T}$:  
\begin{equation}
  \frac{\partial \rho \JJ}{\partial t} + \nabla \cdot \left( \rho \JJ \otimes \u \right) + \nabla T = -\frac{1}{\tau_2} \, \frac{T}{T_0} \frac{\rho_0}{\rho} \, \rho \JJ. 
 \label{eqn.j} 	
\end{equation} 
The Chapman-Enskog expansion of $\JJ$ in terms of the small parameter $\tau_2 \ll 1$ reads 
\begin{equation}
  \JJ = \JJ_0 + \tau_2 \JJ_1 + \tau_2^2 \JJ_2 + ...,  
 \label{eqn.jseries} 
\end{equation} 
which can be directly inserted into \eqref{eqn.j}. Collecting terms of equal powers in $\tau_2$ and setting all the individual coefficients to zero, like in the previous 
section, yields: 
\begin{equation}
 \tau_2^{-1} \underbrace{\left( \frac{T}{T_0} \frac{\rho_0}{\rho}  \, \rho \JJ_0 \right)}_{0} 
+ \tau_2^0 \underbrace{\left( \frac{\partial \rho \JJ_0}{\partial t} + \nabla \cdot \left( \rho \JJ_0 \otimes \u \right) + \nabla T +  \frac{T}{T_0} \frac{\rho_0}{\rho} \, \rho \JJ_1 \right)}_{0} + ... = 0,  
\end{equation} 
hence, the leading zeroth order term and the first order term of the expansion are given by 
\begin{equation} 
  \JJ_0 = 0, \qquad \textnormal{ and } \qquad \JJ_1 = - \frac{T_0}{T \rho_0} \nabla T, 
\end{equation} 
so that the expansion \eqref{eqn.jseries} up to first order terms becomes
\begin{equation} 
 \JJ = - \tau_2 \frac{T_0}{\rho_0} \frac{\nabla T}{T}. 
\label{eqn.asyj} 
\end{equation} 
Inserting \eqref{eqn.asyj} into the heat flux $\mathbf{q} = \alpha^2 T \JJ$ present in the energy equation \eqref{eqn.energy} yields 
\begin{equation}
 \mathbf{q} = \alpha^2 T \JJ = - \alpha^2 \tau_2 \frac{T_0}{\rho_0} \nabla T := -\kappa \nabla T,
\end{equation} 
which is the familiar form of the Fourier heat flux with heat conduction coefficient 
\begin{equation}
\label{eqn.kappa}
\kappa = \alpha^2 \tau_2 \frac{T_0}{\rho_0} 
\end{equation} 
that is recovered for small relaxation times $\tau_2 \ll 1$. 


\subsubsection{On the experimental measurement of the model parameters}

As we have seen in the previous section, the relation between the conventional transport coefficients, the viscosity coefficient $ \mu $ and heat conductivity $ \kappa $, are given by \eqref{eqn.mu} and \eqref{eqn.kappa}, respectively. From these relations, however, it is impossible to recover both parameters, $ \tau_1 $ and $ c_s $ or $ \tau_2 $ and $ \alpha $, of the HPR model, and an experimental way to do that has to be pointed out. So far, we see a possibility to make measurements of the HPR model parameters using experiments on high frequency sound propagation. For example, experimental results from \cite{Sette1955} show that for the vapor of methyl-chloride (CH$ _3 $Cl) at temperature $ 30^{\rm o} $C the low (adiabatic sound speed) and high frequency limits of the longitudinal phase velocity $ V^{long} $ are (for simplicity, we ignore here the heat conducting effect)
\[ c_0\approx250\ {\rm m/s},\ \ \ c_\infty\approx258 {\rm m/s}.\]
Using these data and that $ c_\infty=\sqrt{c_0^2+\frac{4}{3}c_s^2} $ in the HPR model, we get that $ c_s\approx 55.21$ m/s, and the dissipation time is $ \tau_1=6\mu/(\rho_0 c_s^2) \approx 1.545\cdot 10^{-7}$ s, where $ \mu= 1.57\cdot10^{-4}$ Pa$ \cdot $s and $ \rho_0= 2$ kg/m$ ^3 $. In general, the heat conducting effect cannot be ignored, and some extra high frequency heat wave propagation experiments should be conducted in order to measure the characteristic heat wave velocity $ c_h $.




\section{ADER finite volume and ADER discontinuous Galerkin finite element schemes} 
\label{sec:ader}
The equations \eqref{eqn.conti}-\eqref{eqn.energy} of the HPR model described above can be written in the following 
general form of a nonlinear system of hyperbolic PDEs with non-conservative products and stiff source terms: 
\begin{equation}	
\label{eqn.pde.nc}
	\frac{\partial \Q}{\partial t} + \nabla \cdot \bf F(\Q) + \mathbf{B}(\Q) \cdot \nabla \Q = \S(\Q), 
\end{equation}
where $\Q=\Q(\x,t)$ is the state vector; $\x=(x,y) \in \Omega$ is the vector of spatial coordinates and 
$\Omega$ denotes the computational domain; ${\bf F}(\Q) = (\f, \g)$ is the nonlinear flux tensor that contains 
the conservative part of the PDE system and $\mathbf{B}(\Q) \cdot \nabla \Q$ is a genuinely non-conservative 
term. When written in quasilinear form, the system (\ref{eqn.pde.nc}) becomes 
\begin{equation}	\label{eq:Csyst}
\frac{\partial \Q}{\partial t}+ \mathsf{A} (\Q) \cdot\nabla \Q = \S(\Q)\,,
\end{equation}
where the matrix ${\mathsf{A}}(\Q)=\partial {\bf F}(\Q)/\partial \Q + \mathbf{B}(\Q)$ includes both 
the Jacobian of the conservative flux, as well as the non-conservative product. The hyperbolicity of system \eqref{eq:Csyst} 
has been discussed in \cite{PeshRom2014}. However, for the practical implementation of the numerical schemes used in this paper, 
the eigenvectors $\mathbf{R}_n$ of the matrix $\mathsf{A}_n = \mathsf{A}(\Q) \cdot \mathbf{n}$ ($\mathbf{n}$ is a unit-normal vector)
will not be needed, even if they were in principle available. 

The PDE system \eqref{eqn.pde.nc} is solved by resorting to a high order one-step ADER-FV and ADER-DG method 
\cite{QiuDumbserShu,Dumbser2008,ADERNC}, which provides at the same time high order of accuracy in 
both space and time in one single step, hence completely avoiding the Runge-Kutta sub-stages that are typically
used in Runge-Kutta DG and Runge-Kutta WENO schemes. The method will be presented in the \textit{unified} framework
of $P_NP_M$ methods introduced in \cite{Dumbser2008}, which contains both, DG schemes and FV schemes as special cases
of a more general class of methods. For related work on $P_NP_M$ schemes, the reader is referred to 
\cite{luo1,luo2}. The construction of fully-discrete high order one-step schemes is typical of the ADER approach \cite{titarevtoro}.  
In the following we only summarize the main steps, while for more details the reader is referred to 
\cite{Dumbser2008,DumbserZanotti,HidalgoDumbser,GassnerDumbserMunz,Balsara2013,Dumbser2014,Zanotti2015a,Zanotti2015b}.

\subsection{Data representation and reconstruction}

The computational domain $\Omega$ is discretized by a computational mesh (structured or unstructured), 
composed of conforming elements denoted by $T_i$, where the index $i$ ranges from 1 to the 
total number of elements $N_E$. We will further denote the volume (area) of an individual cell by 
$|T_i| = \int_{T_i} d\x$. The discrete solution of PDE \eqref{eqn.pde.nc} is denoted by 
$\u_h(\x,t^n)$ and is represented by piecewise polynomials of maximum degree 
$N \geq 0$. Within each cell $T_i$ we have 
\begin{equation}
\label{eqn.ansatz.uh}
  \u_h(\x,t^n) = \sum_l^{\mathcal{N}} \Phi_l(\x) \hat{\u}^n_{l,i} := \Phi_l(\x) \, \hat{\u}^n_{l,i},  \quad \x \in T_i, 
\end{equation}
where we have introduced the classical Einstein summation convention over two repeated indices. 
The discrete solution $\u_h(\x,t^n)$ is defined in the space of piecewise polynomials up 
to degree $N$, spanned by a set of basis functions $\Phi_l=\Phi_l(\x)$. Throughout this paper we use the 
orthogonal Dubiner-type basis for simplex elements, which is a so-called \textit{modal basis}, 
detailed in \cite{Dubiner,orth-basis}, while we use a tensor-product-type \textit{nodal basis} 
for quadrilateral elements \cite{Dumbser2014}. The nodal basis is given by the Lagrange interpolation polynomials 
passing through the Gauss-Legendre quadrature nodes on the unit square \cite{stroud}. The symbol $\mathcal{N}$ denotes 
the number of degrees of freedom per element and is given by $\mathcal{N}=(N+1)(N+2)/2$ for simplex elements 
and by $\mathcal{N}=(N+1)^2$ for quadrilateral elements in two space dimensions. In the framework of $P_NP_N$ methods, the 
discrete solution $\u_h$ is now \textit{reconstructed} in order to obtain for each element a piecewise polynomial $\w_h(\x,t)$ 
of degree $M\geq N$, with a total number of $\mathcal{M}$ degrees of freedom. Details on the nonlinear WENO reconstruction and on the 
$P_NP_M$ reconstruction can be found in \cite{DumbserKaeser06b,DumbserKaeser07,Dumbser2008} and are not repeated here. The number of degrees 
of freedom $\mathcal{M}$ is again $\mathcal{M}=(M+1)(M+2)/2$ for simplex elements and $\mathcal{M}=(M+1)^2$ for quadrilateral elements in 2D, 
respectively. The reconstruction step is simply abbreviated by $\w_h(\x,t) = \mathcal{R}( \u_h(\x,t) )$, and the reconstruction 
polynomial $\w_h(\x,t)$ is written as 
\begin{equation}
  \w_h(\x,t^n) = \sum_l^{\mathcal{M}} \Psi_l(\x) \hat{\w}^n_{l,i} := \Psi_l(\x) \, \hat{\w}^n_{l,i},  \quad \x \in T_i.  
 \label{eqn.recsol} 
\end{equation} 
Note that for $N=M$ the $P_NP_M$ method reduces to a classical discontinuous Galerkin finite element scheme, with the reconstruction operator
equal to the identity operator, $\mathcal{R}=\mathcal{I}$, or, equivalently, $\w_h(\x,t^n) = \u_h(\x,t^n)$, while for the case $N=0$ the method reduces to a standard high order WENO finite 
volume scheme if a WENO reconstruction operator is adopted. 

For WENO schemes on structured meshes we have found that it is particularly convenient to adopt one-dimensional stencils, each composed by $n_e=M+1$ cells, 
which are subsequently oriented along each spatial direction. The resulting reconstruction is still multidimensional, but implemented in a dimension-by-dimension 
strategy. A complete description of this approach can be found in \cite{AMR3DCL,Zanotti2015}. For unstructured meshes, on the contrary, intrinsically 
multidimensional stencils are built, with $n_e = 2 \mathcal{M}$, where $\mathcal{M} = (M+1)(M+2)/2$.
Moreover, the total number of stencils is seven, i.e. one central stencil, three primary sector stencils and  three reverse sector stencils.  
Further details can be found in \cite{DumbserKaeser07,DumbserKaeser06b,MixedWENO2D,MixedWENO3D}.

In this paper, however, we will only use these two special limits of the general $P_NP_M$ approach, i.e. either $N=0$ (pure FV) or $N=M$ (pure DG). 

\subsection{Local space-time predictor}
\label{sec.predictor}

The discrete solution $\w_h(\x,t^n)$ is now evolved in time according to an element-local weak formulation of the governing PDE in space-time, 
see \cite{DumbserEnauxToro,Dumbser2008,HidalgoDumbser,DumbserZanotti,GassnerDumbserMunz,Balsara2013,Dumbser2014,Zanotti2015a,Zanotti2015b}. 
The local space-time Galerkin method is only used for the construction of an element-local predictor solution of the PDE 
\textit{in the small}, hence neglecting the influence of neighbor elements. This predictor will subsequently be inserted into the corrector 
step described in the next section, which then provides the appropriate coupling between neighbor elements via a numerical flux function 
(Riemann solver) and a path-conservative jump term for the discretization of the non-conservative product. To simplify notation, we define  
\begin{equation}
 \label{eqn.operators1}
  \left<f,g\right> =
      \int \limits_{t^n}^{t^{n+1}} \int \limits_{T_i}  f(\x, t)  g(\x, t)  \, d \x \, d t,
\qquad 
  \left[f,g\right]^{t} =
      \int \limits_{T_i}f(\x, t) g(\x, t) \,  d \x,
\end{equation}
which denote the scalar products of two functions $f$ and $g$ over the space-time element $T_i \times \left[t^n;t^{n+1}\right]$ and  
over the spatial element $T_i$ at time $t$, respectively. Within the local space-time predictor, the discrete solution of 
equation \eqref{eqn.pde.nc} is denoted by $\q_h=\q_h(\x,t)$. 
We then multiply \eqref{eqn.pde.nc} with a space-time test function $\theta_k=\theta_k(\x,t)$ and subsequently integrate over 
the space-time control volume $T_i \times \left[t^n;t^{n+1}\right]$. Inserting $\q_h$, the following weak formulation of 
the PDE is obtained: 
\begin{equation}
\label{eqn.pde.nc.weak1}
 \left< \theta_k, \frac{\partial \q_h}{\partial t}  \right>
    + \left< \theta_k, \nabla \cdot \F \left(\q_h\right) + \mathbf{B}(\q_h) \cdot \nabla \q_h \right> = \left< \theta_k, \S \left( \q_h \right)  \right>.
\end{equation}
The discrete representation of $\q_h$ in element $T_i \times [t^n,t^{n+1}]$ is assumed to have the following form 
\begin{equation}
\label{eqn.st.state}
 \q_h = \q_h(\x,t) =
 \sum \limits_l \theta_l(\x,t) \hat{\q}^n_{l,i} := \theta_l \hat{\q}^n_{l,i},
\end{equation}
where $\theta_l(\x,t)$ is a space-time basis function of maximum degree $M$. 
For the basis functions $\theta_l$ we use the nodal basis given in \cite{Dumbser2008} on simplex elements, 
while we use a tensor-product of 1D nodal basis functions given by the Lagrange interpolation polynomials of the 
Gauss-Legendre quadrature points for quadrilateral elements. 
After integration by parts in time of the first term, eqn. \eqref{eqn.pde.nc.weak1} reads 
\begin{equation}
\label{eqn.pde.nc.dg1}
 \left[ \theta_k, \q_h \right]^{t^{n+1}} - \left[ \theta_k, \w_h(\x,t^n) \right]^{t^n} - \left< \frac{\partial}{\partial t} \theta_k, \q_h \right> 
    + \left< \theta_k, \nabla \cdot \F \left(\q_h\right) + \mathbf{B}(\q_h) \cdot \nabla \q_h \right> = \left< \theta_k, \S \left( \q_h \right)  \right>. 
\end{equation}
Note that the high order polynomial reconstruction of the $P_NP_M$ scheme $\w_h(\x,t^n)$ is taken into account 
in \eqref{eqn.pde.nc.dg1} in a \textit{weak sense} by the term $\left[ \theta_k, \w_h(\x,t^n) \right]^{t^n}$. 
This corresponds to the choice of a numerical flux in time direction, which is nothing else than \textit{upwinding in time}, 
according to the causality principle. 

Note further that due to the DG approximation in space-time, we may have $\q_h(\x,t^n) \neq \w_h(\x,t^n)$ in general, hence
the choice of a numerical flux in time direction is necessary. Note further that in \eqref{eqn.pde.nc.dg1} we have
\textit{not} used integration by parts in space, nor any other coupling to spatial neighbor elements. The integrals
appearing in the weak form \eqref{eqn.pde.nc.dg1}, as well as the space-time test and basis functions involved are 
conveniently written by making use of a space-time reference element $T_e \times [0;1]$. \\ 
The solution of \eqref{eqn.pde.nc.dg1} yields the unknown space-time degrees of freedom $\hat{\q}^n_{l,i}$ for each 
space-time element $T_i \times [t^n; t^{n+1}]$ and is easily achieved with a fast converging iterative scheme, see 
\cite{Dumbser2008,HidalgoDumbser,DumbserZanotti} for more details. 
The above space-time Galerkin predictor has replaced the cumbersome Cauchy-Kovalewski procedure that has been initially 
employed in the original version of ADER finite volume and ADER discontinuous Galerkin schemes 
\cite{schwartzkopff,toro3,toro4,titarevtoro,dumbser_jsc,taube_jsc,DumbserKaeser07}. 

\subsection{Fully discrete one-step finite volume and discontinuous Galerkin schemes}
\label{sec.ADERNC}

At the aid of the local space-time predictor $\q_h$, a fully discrete one-step $P_NP_M$ scheme can now be 
simply obtained by multiplication of the governing PDE system \eqref{eqn.pde.nc} by test functions  
$\Phi_k$, which are identical with the spatial basis functions of the original data representation
before reconstruction, and subsequent integration over the space-time control volume $T_i \times [t^n;t^{n+1}]$. Due 
to the presence of non-conservative products, the jumps of $\q_h$ across element boundaries are taken into account in 
the framework of path-conservative schemes put forward by Castro and Par\'es in the finite volume context 
\cite{Castro2006,Pares2006} and subsequently extended to DG schemes in \cite{Rhebergen2008} and \cite{ADERNC,USFORCE2}, 
where also a generalization to the unified $P_NP_M$ framework has been provided. All these approaches are based on the 
theory of Dal Maso, Le Floch and Murat \cite{DLMtheory}, which gives a definition of weak solutions in the context of 
non-conservative hyperbolic PDE. 
For open problems concerning path-conservative schemes, the reader is referred to \cite{NCproblems}. 

If $\mathbf{n}$ is the outward pointing unit normal vector on the surface $\partial T_i$ of element $T_i$ and  
the path-conservative jump term in normal direction is denoted by $\mathcal{D}^-\left(\q_h^-, \q_h^+ \right) \cdot\mathbf{n}$, 
which is a function of the left and right boundary-extrapolated data, 
$\q_h^-$ and $\q_h^+$, respectively, then we obtain the following path-conservative one-step 
$P_NP_M$ scheme, see \cite{ADERNC}: 
\begin{equation}
\label{eqn.pde.nc.gw2}
\begin{split}
\left( \int \limits_{T_i} \Phi_k \Phi_l d\x \right) \left( \hat{\u}_l^{n+1} -  \hat{\u}_l^{n} \right) +
\int \limits_{t^n}^{t^{n+1}} \int \limits_{\partial T_i} \Phi_k \, \mathcal{D}^-\left(\q_h^-, \q_h^+ \right)\cdot\mathbf{n} \, dS dt 
\\ 
+\int\limits_{t^n}^{t^{n+1}} \int \limits_{T_i \backslash \partial T_i} \Phi_k \left( \nabla \cdot \F\left(\q_h \right) + \mathbf{B}(\q_h) \cdot \nabla \q_h \right) d\x dt  
= \int \limits_{t^n}^{t^{n+1}} \int \limits_{T_i} \Phi_k \S(\q_h) d\x dt. 
\end{split} 
\end{equation}
The element mass matrix appears in the first integral of \eqref{eqn.pde.nc.gw2}, the second term accounts for the jump in the discrete 
solution at element boundaries and the third term takes into account the smooth part of the non-conservative product. 
For general complex nonlinear hyperbolic PDE systems we use the simple Rusanov method \cite{Rusanov:1961a} (also called the local Lax Friedrichs
method), although any other kind of Riemann solver could be also used, see \cite{toro-book} for an overview of state-of-the-art Riemann solvers. 
At that point we would also like to point out the new general reformulation of the HLLEM Riemann solver of Einfeldt and Munz 
\cite{Einfeldt88,munz91}, within the setting of path-conservative schemes recently forwarded in \cite{NewHLLEM}, as well as the family of 
MUSTA schemes, which has been applied to the equations of nonlinear elasticity in \cite{TitarevRomenskiToro}.
\\

The path-conservative Rusanov jump term reads  
\begin{equation}
  \mathcal{D}^-\left(\q_h^-, \q_h^+ \right)\cdot\mathbf{n} = \halb \left( \F(\q_h^+) - \F(\q_h^-) \right) \cdot \mathbf{n} + 
	\halb \left( \tilde{\mathbf{B}} \cdot \mathbf{n} - s_{\max} \mathbf{I} \right) \left( \q_h^+ - \q_h^- \right), 
	\label{eqn.rusanov} 
\end{equation} 
with the maximum signal speed at the element interface $s_{\max} = \max\left( \left|\boldsymbol{\Lambda}(\q_h^+) \right|, \left|\boldsymbol{\Lambda}(\q_h^-) \right| \right)$ and
the matrix $\tilde{\mathbf{B}} \cdot \mathbf{n}$ given by the following path-integral along a straight line segment path $\psi$: 
\begin{equation}
 \tilde{\mathbf{B}} \cdot \mathbf{n} = \int \limits_0^1 \mathbf{B}\left( \psi(\q_h^-, \q_h^+, s \right) \cdot \mathbf{n} \, ds, 
\qquad 
\psi \left( \q_h^-, \q_h^+, s \right) = \q_h^- + s \left( \q_h^+ - \q_h^- \right).  
\end{equation} 
According to the suggestions made in \cite{ADERNC,USFORCE2,OsherUniversal,ApproxOsher,OsherNC}, the path-integrals can be conveniently 
evaluated numerically by the use of a classical Gauss-Legendre quadrature formula on the unit interval $[0;1]$. For an alternative choice 
of the path, see \cite{MuellerToro1,MuellerToro2}. 

This completes the brief description of the $P_NP_M$ scheme used for the discretization of the governing PDE system 
\eqref{eqn.pde.nc}. In the case of ADER-WENO finite volume schemes, we simply have $N=0$, $\mathcal{N}=1$, $\Phi_k=1$, and 
the limiter is directly incorporated in the \textit{nonlinear} reconstruction operator $\w_h(\x,t^n) = \mathcal{R}\left(\u_h(\x,t^n) \right)$, 
while for ADER-DG schemes ($N=M$, $\Phi_k = \Psi_k$) a new family of \aposteriori 
sub-cell finite volume limiters has been forwarded in \cite{Dumbser2014,Zanotti2015a,Zanotti2015b}. For alternative finite volume 
subcell limiters in the context of DG schemes, see the work of Sonntag \& Munz \cite{Sonntag} and Meister \& Ortleb \cite{MeisterOrtleb}.


%
\section{Numerical results}       
\label{sec:results}

\subsection{Numerical convergence studies in the stiff inviscid limit}  
\label{sec:shuvortex}

We first present a numerical convergence study on a smooth unsteady flow, for which an exact analytical solution is known for the compressible Euler equations, i.e. 
in the inviscid limit $\tau_1 \to 0$ and $\tau_2 \to 0$ of the HPR model. 

The computational setup is the classical one of a convected isentropic vortex, see \cite{balsarashu,HuShuTri}. 
The initial condition is given in terms of primitive variables and it consists in a linear superposition of a homogeneous background field and some perturbations $\delta$: 
\begin{equation}
\label{ShuVortIC}
(\rho, u, v, p) = (1+\delta \rho, 1+\delta u, 1+\delta v, 1+\delta p).
\end{equation}  
We furthermore set the distortion tensor initially to $\mathbf{A} = \sqrt[3]{\rho} \, \mathbf{I}$, while the heat flux vector is initialized with $\mathbf{J}=0$. The radial 
coordinate is related to the Cartesian coordinates $x$ and $y$ by the relation $r^2=(x-5)^2+(y-5)^2$. The vortex strength is chosen as $\epsilon=5$ and the perturbation of entropy 
$S=\frac{p}{\rho^\gamma}$ is assumed to be zero, while the perturbations of temperature $T$ and velocity $\mathbf{v}$ are given by 
\begin{equation}
\label{ShuVortDelta}
\left(\begin{array}{c} \delta u \\ \delta v \end{array}\right) = \frac{\epsilon}{2\pi}e^{\frac{1-r^2}{2}} \left(\begin{array}{c} -(y-5) \\ \phantom{-}(x-5) \end{array}\right), \quad \delta S = 0, \quad \delta T = -\frac{(\gamma-1)\epsilon^2}{8\gamma\pi^2}e^{1-r^2}. 
\end{equation} 
From \eqref{ShuVortDelta} it follows that the perturbations for density and pressure are given by 
\begin{equation}
\label{rhopressDelta}
\delta \rho = (1+\delta T)^{\frac{1}{\gamma-1}}-1, \quad \delta p = (1+\delta T)^{\frac{\gamma}{\gamma-1}}-1. 
\end{equation} 

The computational domain is the square $\Omega=[0;10]\times[0;10]$ and periodic boundary conditions are applied everywhere. The reference solution $\Q_e$ is given by the exact solution 
of the compressible Euler equatons. In the inviscid case $\Q_e$ is simply the time--shifted initial condition $\Q_e(\x,t)=\Q(\x-\v_c t,0)$, where the convective mean velocity is 
$\v_c=(1,1)$. The test problem is run on a sequence of successively refined meshes until a final time of $t=1.0$. The chosen physical parameters are 
$\gamma = 1.4$, $c_v = 2.5$, $\rho_0=1$, $c_s=0.5$ and $\alpha=1$. The resulting numerical convergence rates obtained with ADER-DG schemes using polynomial approximation degrees 
from $N=M=2$ to $N=M=5$ are listed in Table \ref{tab.conv1}, together with the chosen values for the effective viscosity $\mu$ and the effective heat conductivity 
coefficient $\kappa$. For the higher order schemes, it was necessary to use smaller values of $\kappa$ and $\mu$, since the scheme otherwise converges to the solution of the viscous problem,
while the reference solution is given by the exact solution of the inviscid problem (compressible Euler equations). From Tab. \ref{tab.conv1} one can observe that high order of 
convergence of the numerical method is achieved also in the stiff limit of the governing PDE system.

\begin{table}  
\caption{Numerical convergence results for ADER-DG schemes applied to the HPR model ($c_s=0.5$, $\alpha=1$) in the low viscosity relaxation limit ($\mu \ll 1, \kappa \ll 1$). 
Results are shown for the density $\rho$ at a final time of $t=1$. The reference solution is given by the exact solution of the inviscid compressible Euler equations.} 
\begin{center} 
\begin{small}
\renewcommand{\arraystretch}{1.0}
\begin{tabular}{ccccccc} 
\hline
  $N_x$ & $\epsilon({L_1})$ & $\epsilon({L_2})$ & $\epsilon({L_\infty})$ & $\mathcal{O}(L_1)$ & $\mathcal{O}(L_2)$ & $\mathcal{O}(L_\infty)$ \\ 
\hline
  \multicolumn{7}{c}{ADER-DG $P_2P_2$ ($\mu=\kappa=10^{-6}$)}  \\
\hline
20	& 9.4367E-03	& 2.2020E-03	& 2.1633E-03	&      &      &		    \\
40	& 1.9524E-03	& 4.4971E-04	& 4.2688E-04	& 2.27 & 2.29	& 2.34  \\
60	& 7.5180E-04	& 1.7366E-04	& 1.4796E-04	& 2.35 & 2.35	& 2.61  \\
80	& 3.7171E-04	& 8.6643E-05	& 7.3988E-05	& 2.45 & 2.42	& 2.41  \\
\hline
  \multicolumn{7}{c}{ADER-DG $P_3P_3$ ($\mu=\kappa=10^{-6}$)}  \\
\hline
10	& 1.7126E-02	& 4.0215E-03	& 3.6125E-03	& 		  &      &        \\
20	& 6.0405E-04	& 1.7468E-04	& 2.1212E-04	& 4.83	& 4.52 & 	4.09  \\
30	& 8.3413E-05	& 2.5019E-05	& 2.7576E-05	& 4.88	& 4.79 & 	5.03  \\
40	& 2.1079E-05	& 6.0168E-06	& 7.6291E-06	& 4.78	& 4.95 & 	4.47  \\
\hline
  \multicolumn{7}{c}{ADER-DG $P_4P_42$ ($\mu = \kappa = 10^{-7}$)}   \\
\hline
10	& 1.5539E-03	& 4.5965E-04	& 5.1665E-04	& 		 &        &       \\ 
20	& 4.3993E-05	& 1.0872E-05	& 1.0222E-05	& 5.14 & 	5.40	& 5.66  \\
25	& 1.8146E-05	& 4.4276E-06	& 4.1469E-06	& 3.97 & 	4.03	& 4.04  \\
30	& 8.6060E-06	& 2.1233E-06	& 1.9387E-06	& 4.09 &	4.03	& 4.17  \\
\hline
  \multicolumn{7}{c}{ADER-DG $P_5P_5$ ($\mu = \kappa = 10^{-7}$)}   \\
\hline
 5	& 1.1638E-02	& 1.1638E-02	& 1.8898E-03	& 		 &       &       \\
10	& 3.9653E-04	& 9.3717E-05	& 6.5319E-05	& 4.88 & 	6.96 & 	4.85 \\
15	& 4.4638E-05	& 1.2572E-05	& 1.9056E-05	& 5.39 & 	4.95 & 	3.04 \\
20	& 9.6136E-06	& 3.0120E-06	& 3.9881E-06	& 5.34 & 	4.97 & 	5.44 \\
\hline 
\end{tabular}
\end{small}
\end{center}
\label{tab.conv1}
\end{table}

\subsection{The first problem of Stokes}  
\label{sec:stokes}

\begin{figure}[!htbp]
\begin{center}
\begin{tabular}{ccc} 
\includegraphics[width=0.3\textwidth]{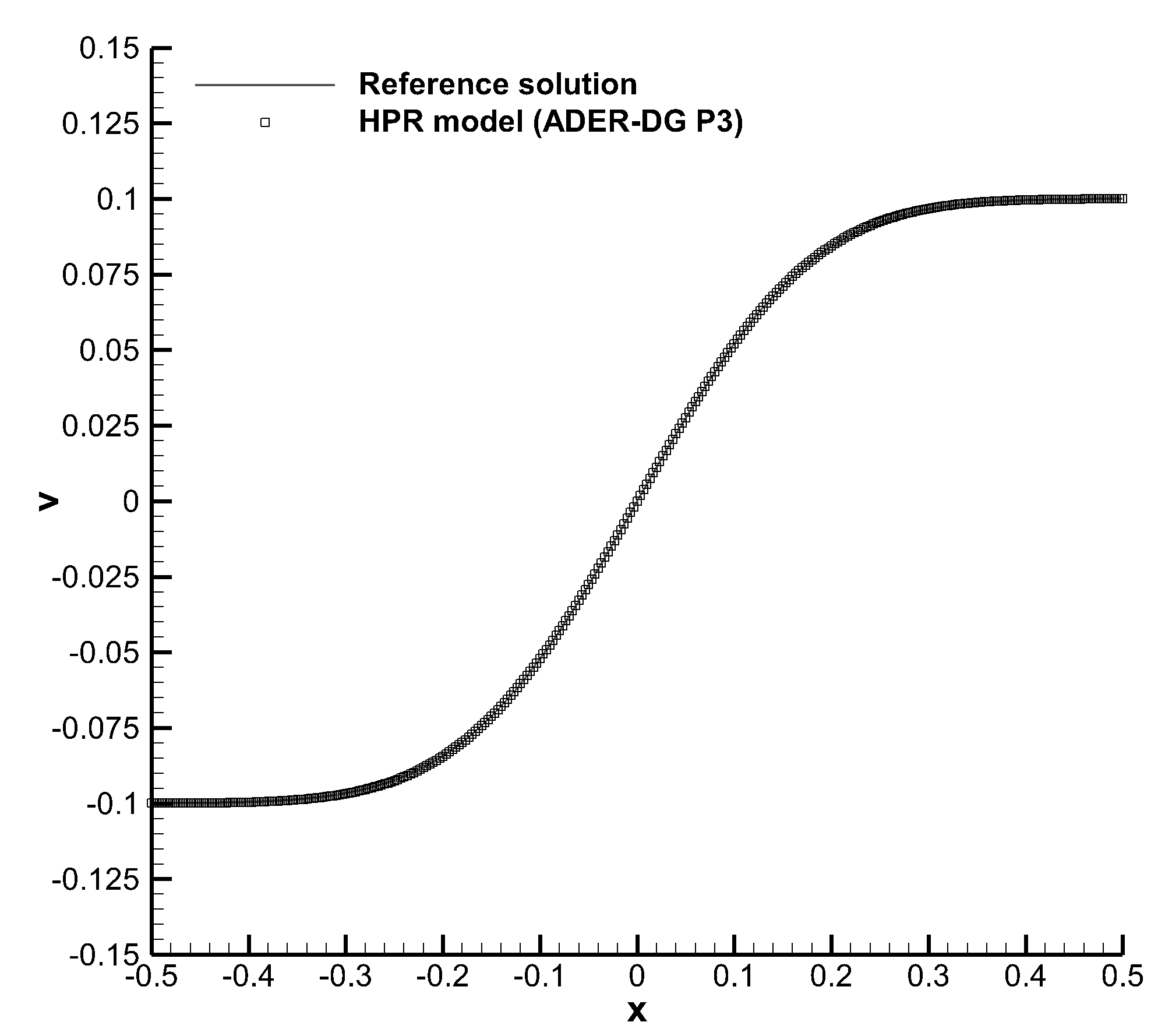} & 
\includegraphics[width=0.3\textwidth]{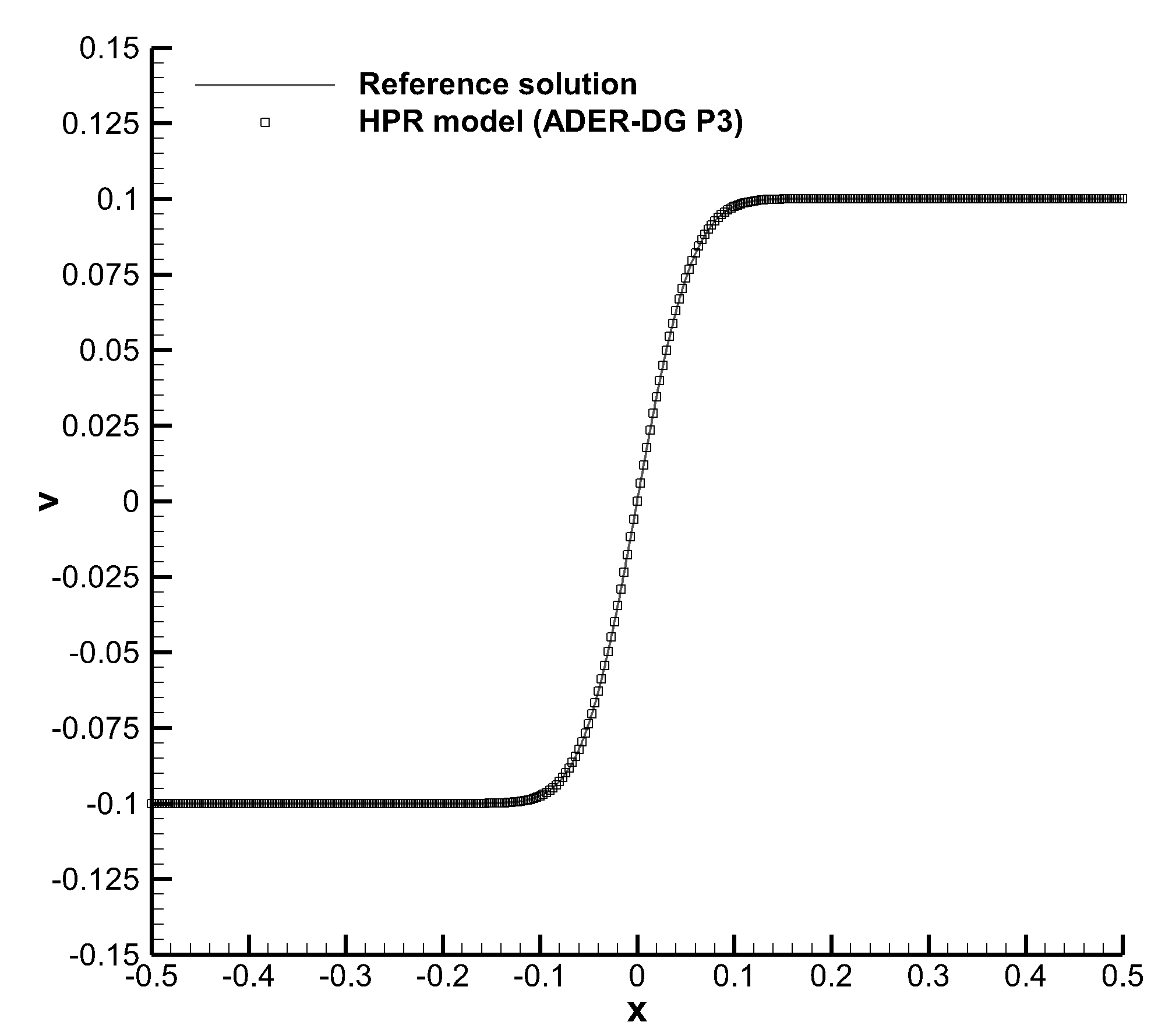} & 
\includegraphics[width=0.3\textwidth]{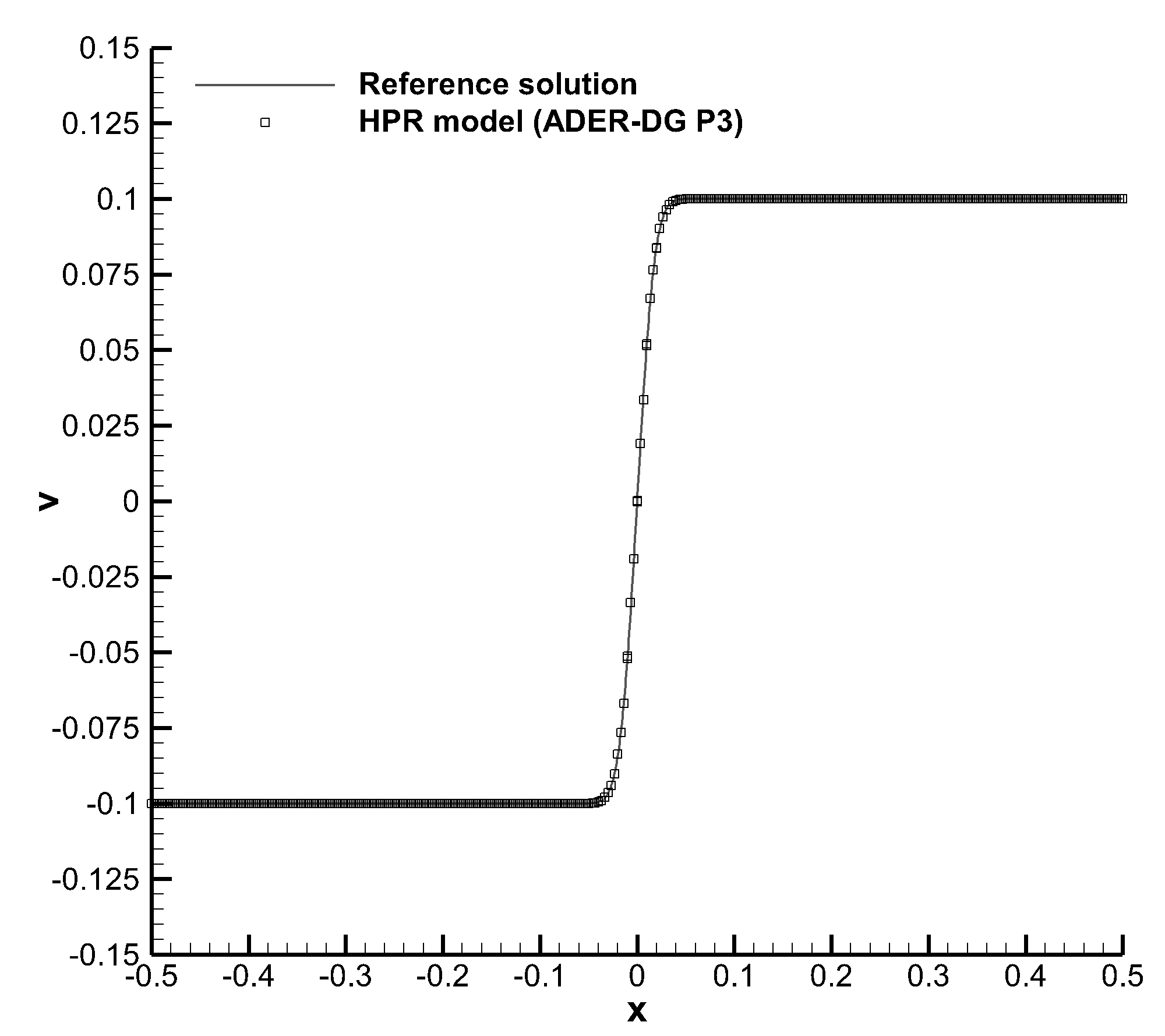}  
\end{tabular} 
\caption{Exact solution of the first problem of Stokes for the Navier-Stokes equations and numerical solution for the hyperbolic model of Peshkov and Romenski (HPR) obtained 
with an ADER-DG $P_3P_3$ scheme at a final time of $t=1.0$ with different viscosities: $\mu=10^{-2}$ (left), $\mu=10^{-3}$ (middle), $\mu=10^{-4}$ (right). } 
\label{fig.stokes}
\end{center}
\end{figure}

There are very few test problems for which an exact analytical solution of the unsteady Navier-Stokes equations is known. 
One of those is the first problem of Stokes \cite{BLTheory}, which consists of the time-evolution of an infinite incompressible 
shear layer. To get an almost incompressible behavior, we run the simulation at a low Mach number of $M=0.1$. 
The computational domain is $\Omega = [-0.5;+0.5] \times [-0.05;+0.05]$, with periodic boundary conditions in $y$ direction. 
At the boundaries in $x$ direction the initial condition is imposed. The initial condition of the problem is given by 
$\rho=1$, $u=0$, $p=1/\gamma$, $\mathbf{A}=\mathbf{I}$, $\mathbf{J}=0$, while the velocity component $v$ is $v=-v_0$ for $x<0$
and $v=+v_0$ for $x\geq 0$. The physical parameters of this test problem are set to $v_0=0.1$, $\gamma=1.4$, $c_v=1$, $\rho_0=1$, $c_s=1$ 
and $\alpha=\kappa=0$. Simulations are performed with an ADER-DG $P_3P_3$ scheme ($N=M=3$) on a grid composed of $100 \times 10$ elements 
up to a final time of $t=1$. 
The exact solution of the incompressible Navier-Stokes equations for the velocity component $v$ is given by 
\begin{equation}
\label{eqn.stokes} 
  v(x,t) = v_0 \textnormal{erf}\left( \halb \frac{x}{\sqrt{\mu t}} \right),  
\end{equation} 
and serves as a reference solution for the HPR model. 
The comparison between the Navier-Stokes reference solution \eqref{eqn.stokes} and the numerical results obtained for the HPR 
model are presented in Fig. \ref{fig.stokes}, where one can observe an excellent agreement between the two for various viscosities
$\mu$.

\subsection{Laminar boundary layer over a flat plate} 
\label{sec:blasius}

The laminar flow over a flat plate has been studied by Prandtl in his famous paper \cite{Prandtl1904}, where the concept of boundary layers was introduced in 
fluid mechanics for the first time. The boundary layer equations proposed by Prandtl were then solved for the first time in the special case of a 
laminar flow over a flat plate by Blasius in \cite{Blasius1908}. For an overview of boundary layer theory, the reader is referred to the well-known textbook 
by Schlichting and Gersten \cite{BLTheory}. In the case of incompressible flow, the boundary layer equations take the following simple form: 
\begin{equation}
\label{eqn.comp.bl}
   f_{\eta \eta \eta} + f f_{\eta \eta} = 0, \qquad \textnormal{with} \qquad f(0)=0, \qquad f_\eta(0)=0, \qquad \lim \limits_{\eta \to \infty} f_\eta(\eta) = 1.  
\end{equation}
Here, $f=f(\eta)$ is the dimensionless stream function in the similarity variable $\eta = y \sqrt{\frac{U_\infty}{2 \nu x}}$, while the axial
flow velocity is given by $u = U_\infty f_\eta$. Nowadays, the boundary layer equation \eqref{eqn.comp.bl} can be solved by any standard ODE solver 
in combination with a classical shooting technique. In this paper, however, we use a shooting method based on the ODE solver proposed in \cite{ADERNSE}, 
which is a special case of the space-time Galerkin predictor of the ADER approach, but applied to the simple case of a pure ODE. 

The setup of the proposed numerical test case is as follows: the computational domain is $\Omega=[0;1.5] 
\times [0;0.4]$ and is discretized with $75 \times 100$ rectangular elements. At $y=0$ we impose a no-slip 
wall boundary condition and the chosen Reynolds number of the flow is $\Re=10^3$. The initial condition and the
physical parameters for the computational setup are $\gamma=1.4$, $\rho_0=1$, $c_v=1$, $c_s=8$, $\rho=1$, $u=U_\infty=1$, 
$v=V_\infty=0$, $p=100/\gamma$, $\mathbf{A}=\mathbf{I}$, $\mathbf{J}=0$, $\alpha=\kappa=0$ and $\mu=10^{-3}$.  
The Mach number of this setup is therefore $M_\infty=0.1$. At $x=0$ the inflow boundary condition is given
by the free stream data, i.e. by the initial condition. Simulations are run up to $t=10$ using a third order
ADER-WENO finite volume scheme ($N=0$, $M=2$). 

\begin{figure}[!htbp]
\begin{center}
\begin{tabular}{cc} 
\includegraphics[width=0.6\textwidth]{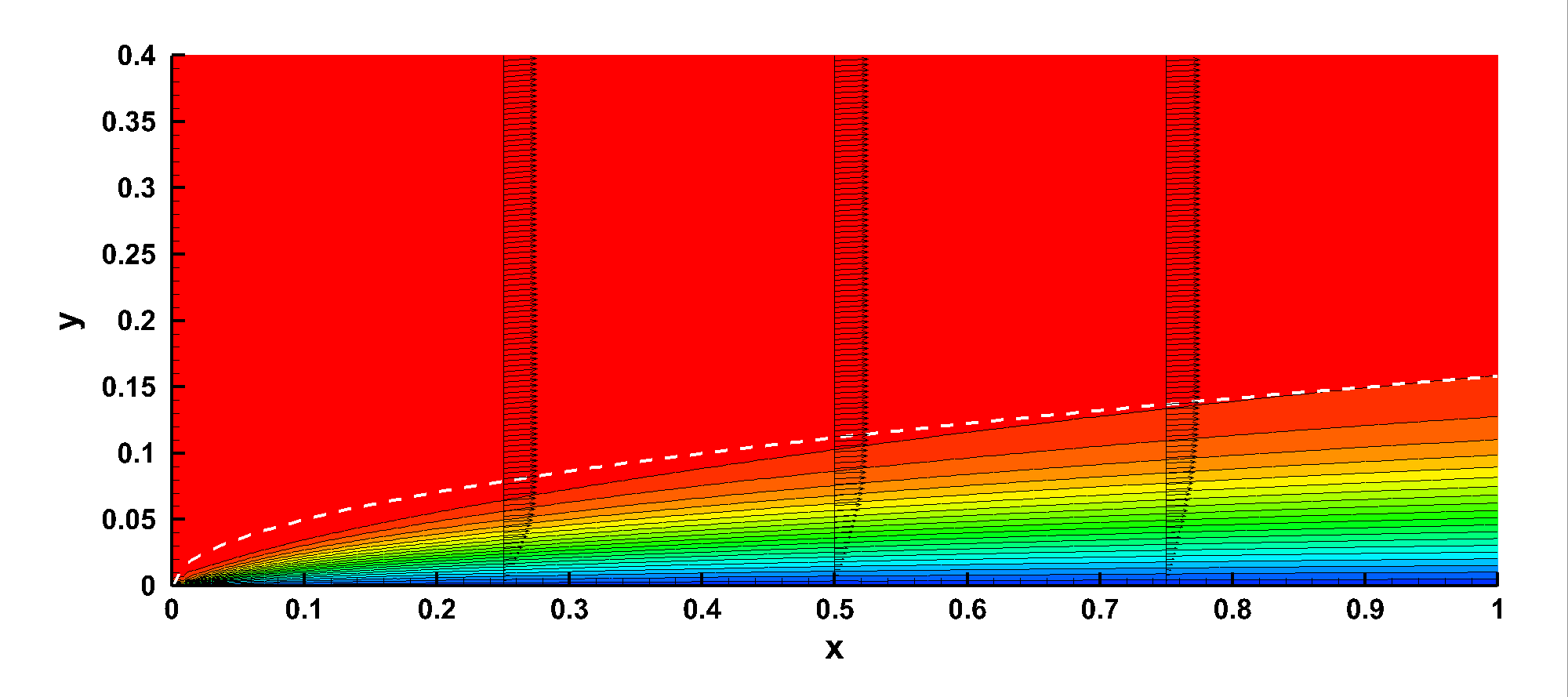} & 
\includegraphics[width=0.3\textwidth]{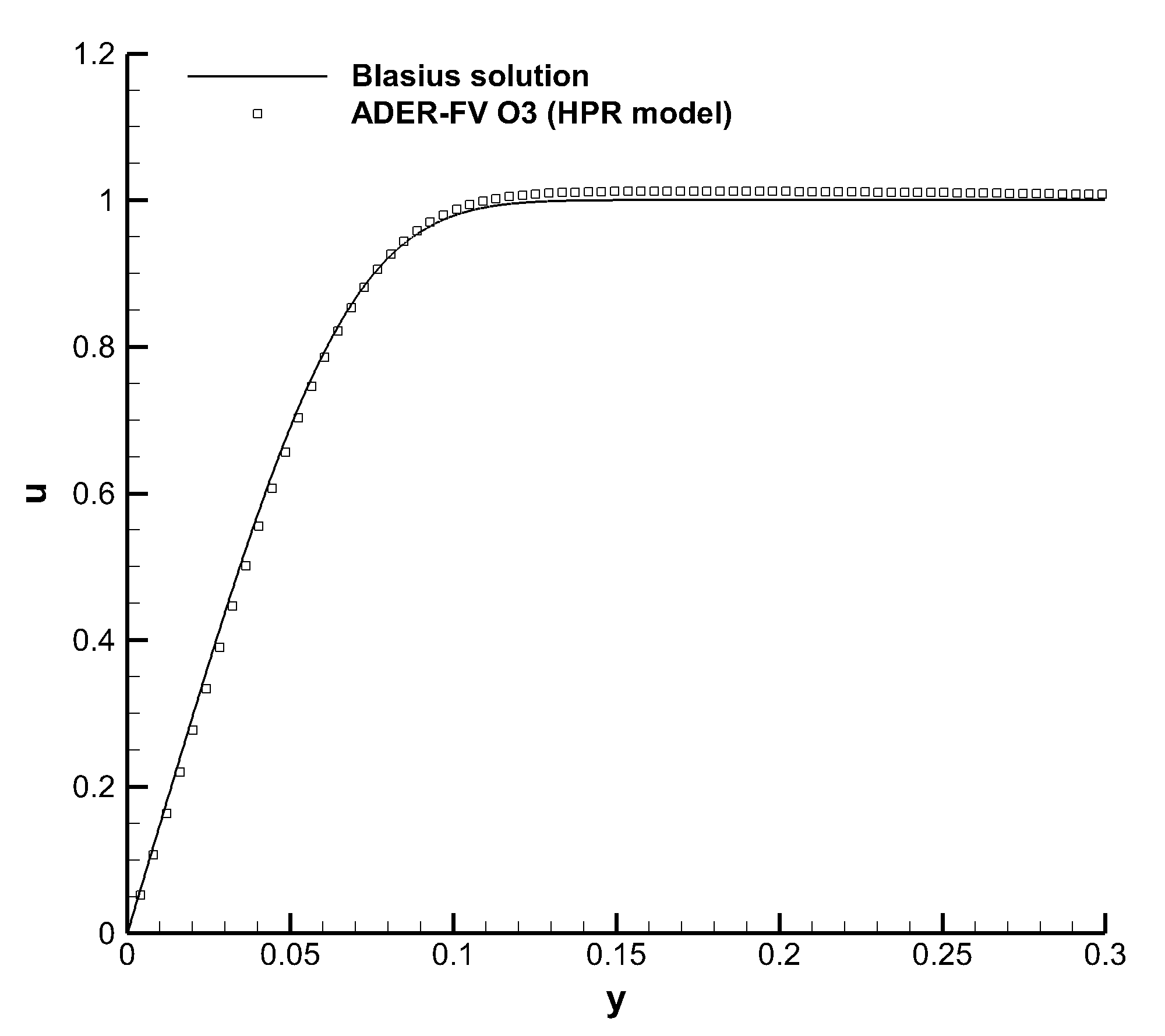}  \\
\end{tabular} 
\caption{Laminar boundary layer over a flat plate. Numerical solution for the hyperbolic model of Peshkov and Romenski (HPR) obtained 
with a third order ADER-WENO finite volume scheme at a final time of $t=10.0$. Left: boundary layer thickness $\delta_{0.99}$ (dashed white line), 
velocity contours and some velocity profiles. Right: vertical cut through the velocity profile along the line $x=0.5$ and comparison with the Blasius 
solution.  } 
\label{fig.blasius1}
\end{center}
\end{figure}

\begin{figure}[!htbp]
\begin{center}
\begin{tabular}{cc} 
\includegraphics[width=0.49\textwidth]{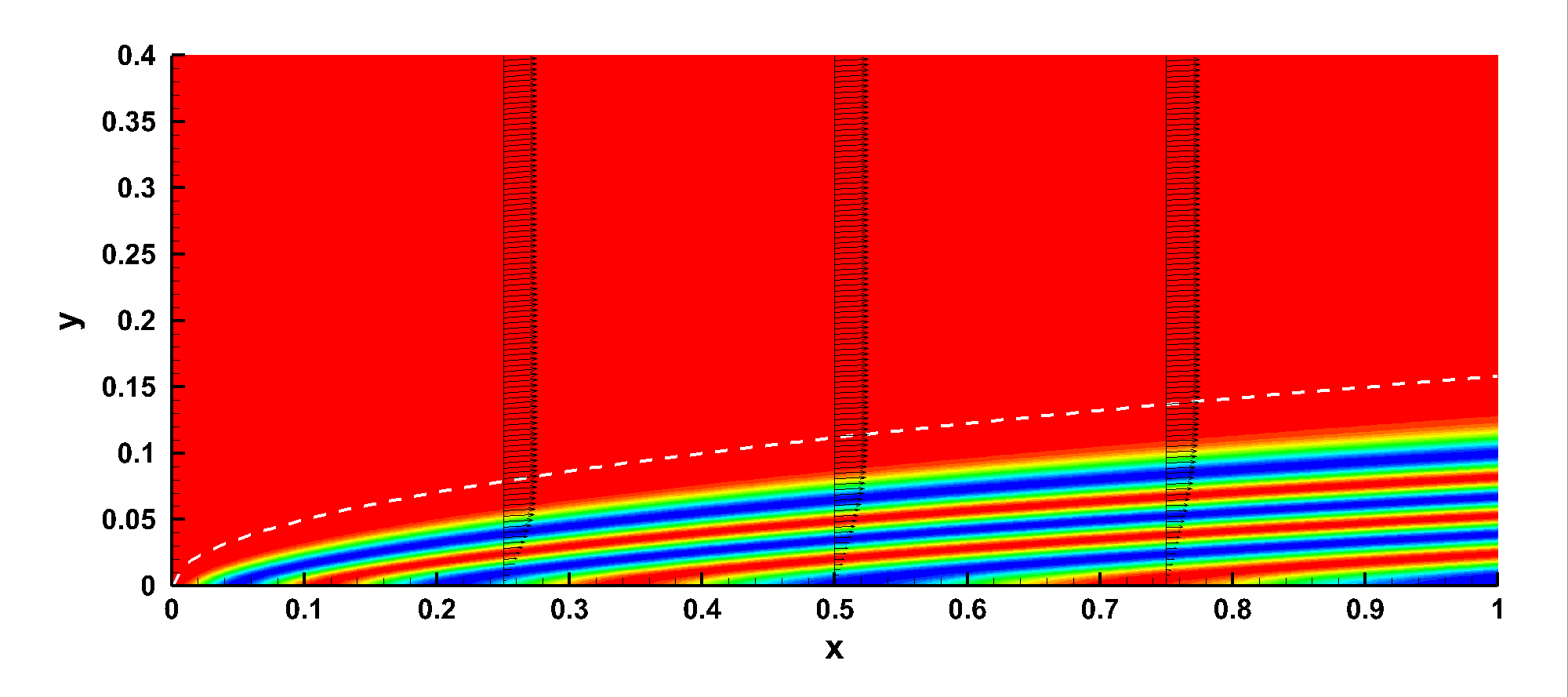} & 
\includegraphics[width=0.49\textwidth]{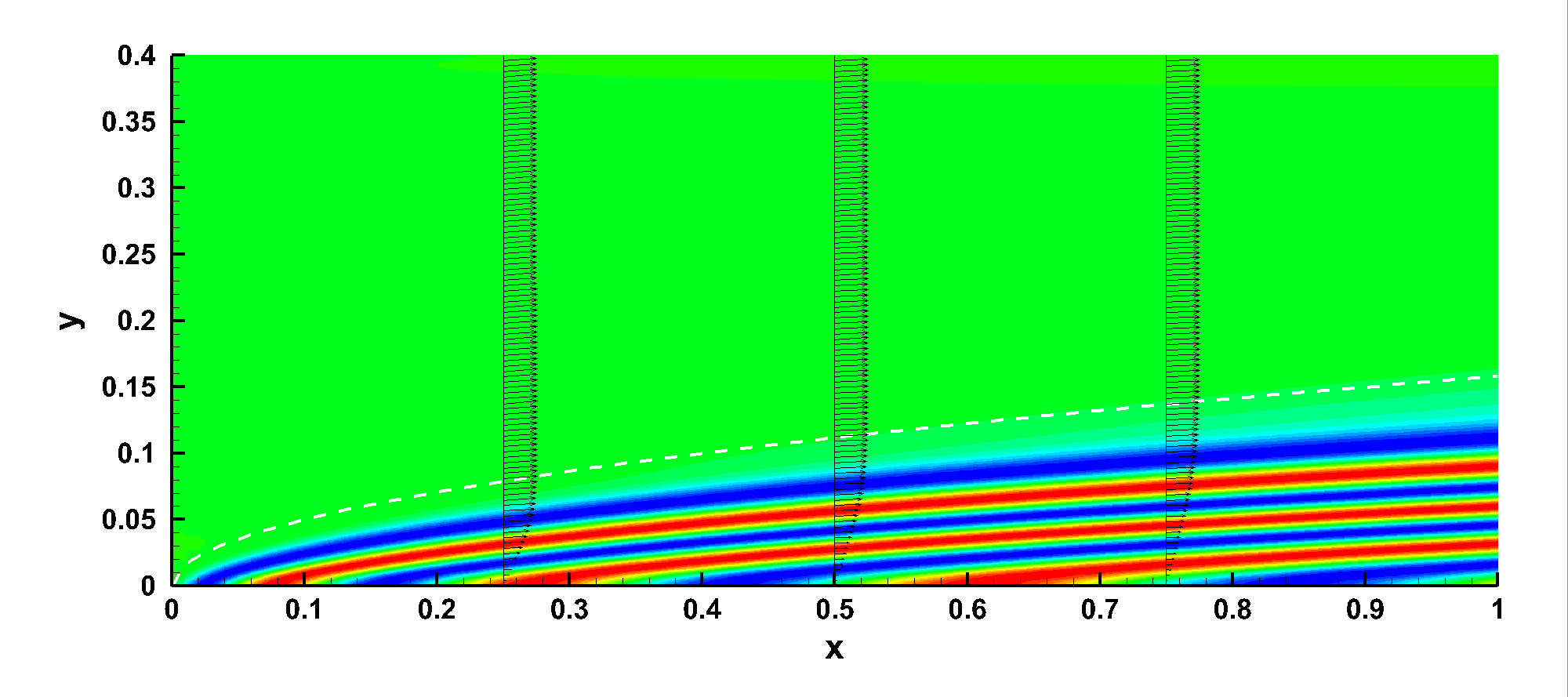}   
\end{tabular} 
\caption{Laminar boundary layer over a flat plate. Numerical solution for the hyperbolic model of Peshkov and Romenski (HPR) obtained 
with a third order ADER-WENO finite volume scheme at a final time of $t=10.0$. 41 equidistant color contours in the interval [-1,1] for the 
distortion tensor component $A_{11}$ (left) and $A_{12}$ (right). } 
\label{fig.blasius2}
\end{center}
\end{figure}

In Fig. \ref{fig.blasius1} the computational results obtained for the HPR model are shown, together with a 1D cut through 
the numerical solution at $x=0.5$ and a comparison with the Blasius reference solution is made. 
A good agreement between the numerical solution of the HPR model and the Blasius solution can be noted, despite the 
fact that two completely different mathematical models have been used to obtain them. 
This confirms the validity of the HPR model in the stiff relaxation limit when $\tau_1 \to 0$, where it is able to 
accurately reproduce the known results from Navier-Stokes theory. For the sake of completeness, in Fig. \ref{fig.blasius2} 
we show two of the components of the distortion tensor $\mathbf{A}$. 

\subsection{Hagen-Poiseuille flow in a duct}   
\label{sec:duct}
\begin{figure}[!htbp]
\begin{center}
\begin{tabular}{c} 
\includegraphics[width=0.95\textwidth]{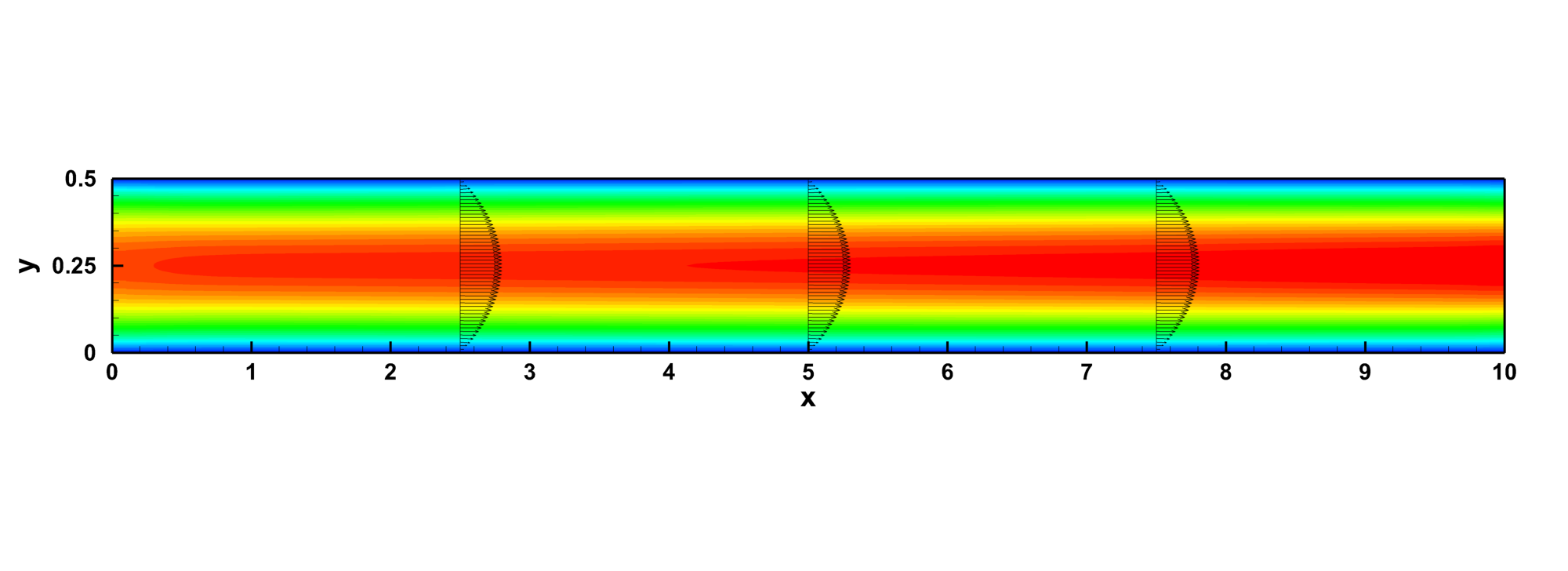}  \\ 
\includegraphics[width=0.45\textwidth]{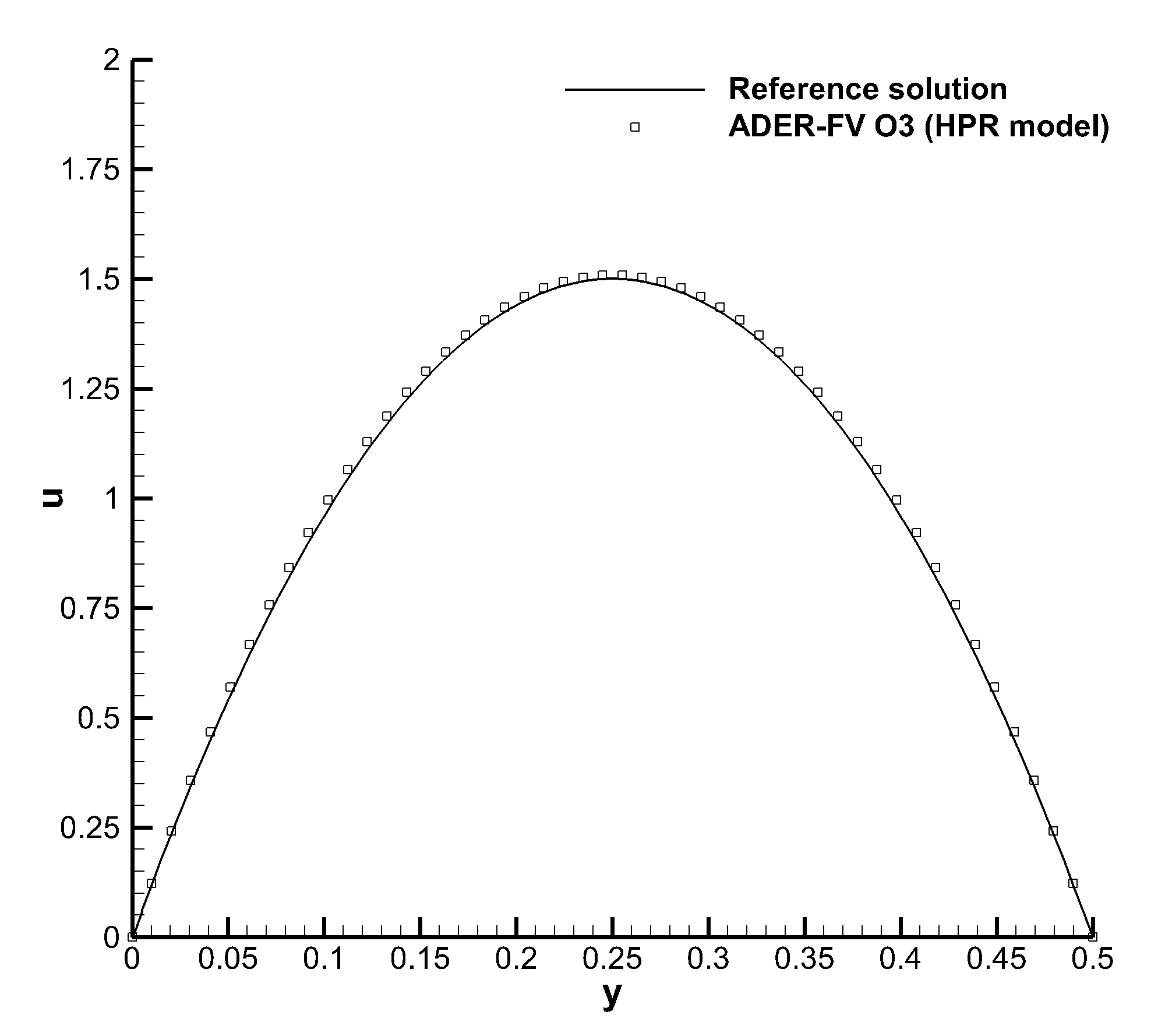}  
\end{tabular} 
\caption{Steady laminar Hagen-Poiseuille flow in a duct at $Re=50$. Exact solution of the Navier-Stokes equations and numerical solution for the hyperbolic model of Peshkov and Romenski (HPR) at a 
final time of $t=10.0$ obtained with a third order ADER-WENO finite volume scheme. Axial velocity contours with some velocity profiles (top) and 1D cut along the $y$ axis at $x=5$ (bottom).  } 
\label{fig.poiseuille}
\end{center}
\end{figure}

Here, we consider the steady flow of a viscous Newtonian fluid in a rectangular duct of length $L$ and height $h$ in the presence of a constant pressure 
gradient $\Delta p < 0$, and choosing the $x-$ axis as the direction of motion. This test, referred to as the \emph{Hagen-Poiseuille flow}, has a well 
known solution of the Navier-Stokes equations \cite{Landau-Lifshitz6} with a parabolic velocity profile given by 
\begin{equation} 
\label{eqn.hp}
v = \halb \frac{\Delta p}{L} \frac{\rho}{\mu} y(y-h). 
\end{equation} 
Although Eq.~(\ref{eqn.hp}) has been derived for an incompressible fluid, we still expect to obtain a good numerical agreement with it, as long as our simulations are performed in the low Mach number regime. For that purpose, we use the following physical parameters and initial condition for our simulation: $\gamma=1.4$, $\rho_0=1$,  
$c_v=1$, $c_s=8$, $p=100/\gamma$, $u=v=0$, $\mathbf{A}=\mathbf{I}$, $\mathbf{J}=0$, $\alpha=\kappa=0$ and $\mu=10^{-2}$. The pressure gradient is imposed between the left
inflow and the right outlet as $\Delta p = -4.8$, leading to a mean flow velocity of $\bar u = 1$ and a maximum flow velocity of $u_{\max}=1.5$.  
We have solved the problem in the computational domain $\Omega=[0,10]\times[0,0.5]$ covered by $100\times50$ cells, applying a third order ADER-WENO finite volume scheme 
 ($N=0$, $M=2$) to the HPR model. The computational results are shown at time $t=10$ in Fig.~\ref{fig.poiseuille}, referring to a Reynolds number of $\Re=\frac{\bar u h \rho}{\mu}=50$. 
In the top panel the laminar flow is very well reproduced, with only a moderate increase of the flow velocity from $x=0$ to $x=10$. In the bottom panel we perform a 
direct comparison to the Navier-Stokes reference solution, by plotting the velocity across the flow as measured at $x=5$. We conclude that the HPR model can successfully 
solve this classical test of laminar, steady viscous flow in a duct. 

\subsection{The lid-driven cavity}   
\label{sec:cavity}

\begin{figure}[!htbp]
\begin{center}
\begin{tabular}{cc} 
\includegraphics[width=0.45\textwidth]{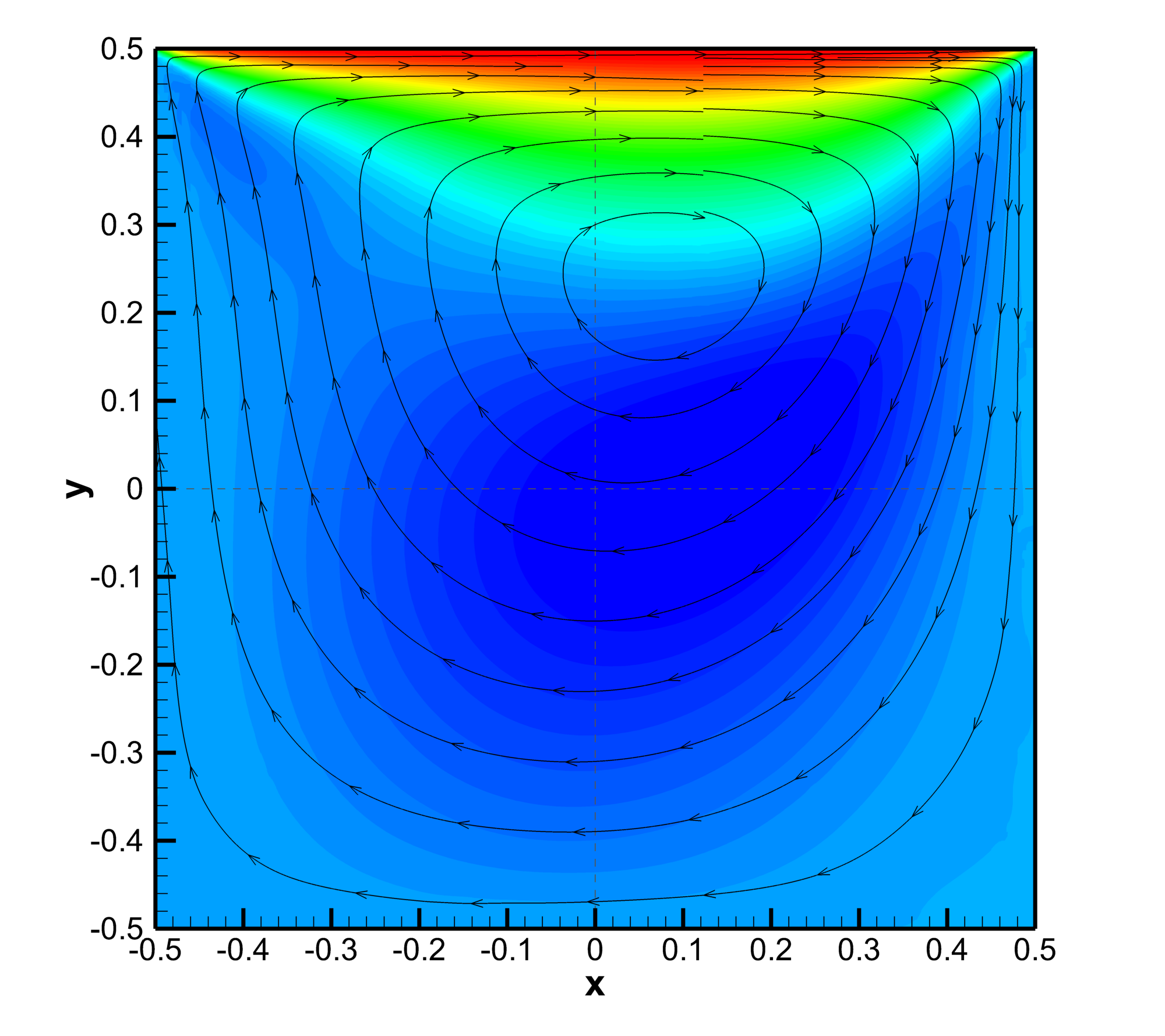} & 
\includegraphics[width=0.45\textwidth]{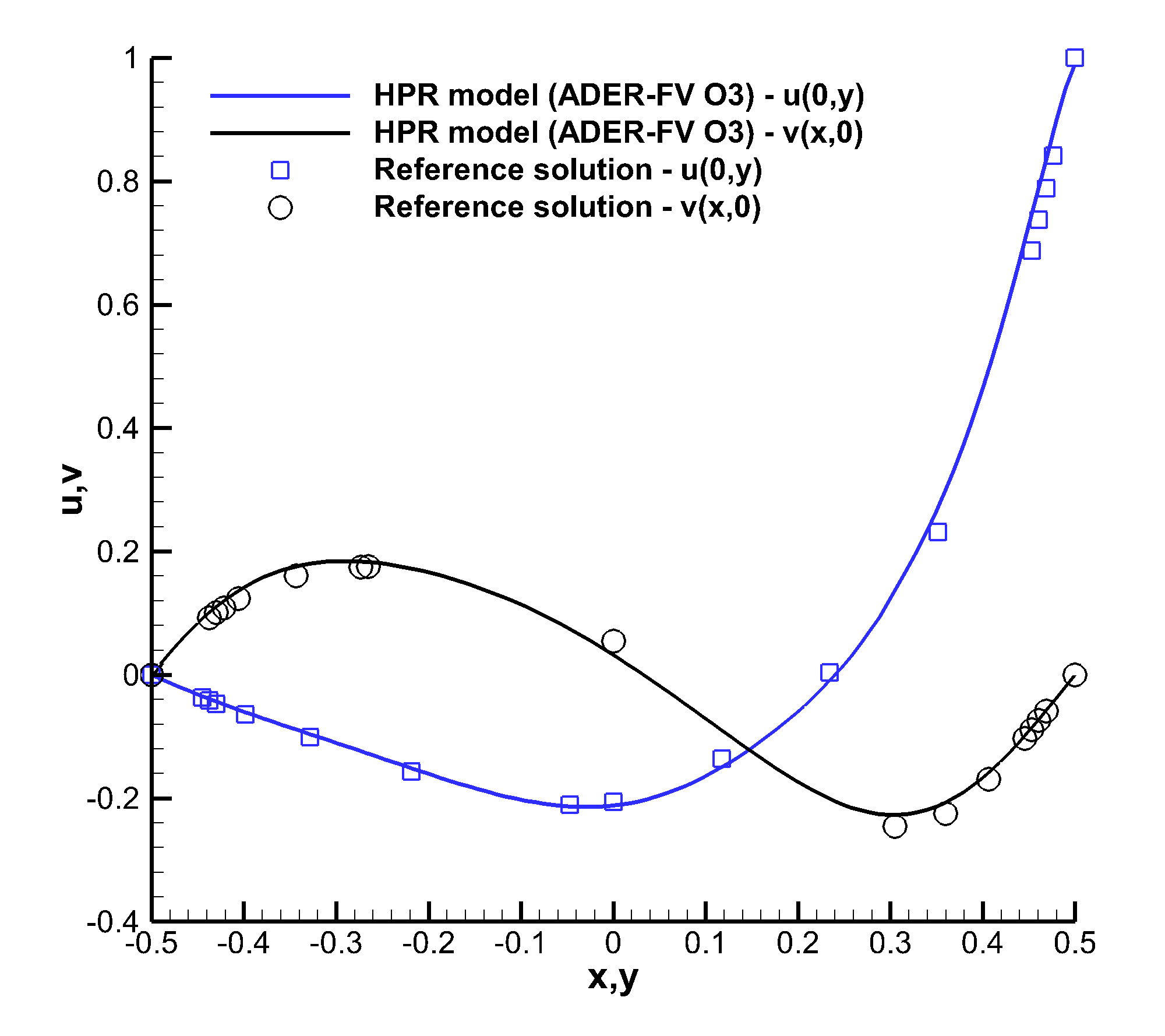}  \\ 
\includegraphics[width=0.45\textwidth]{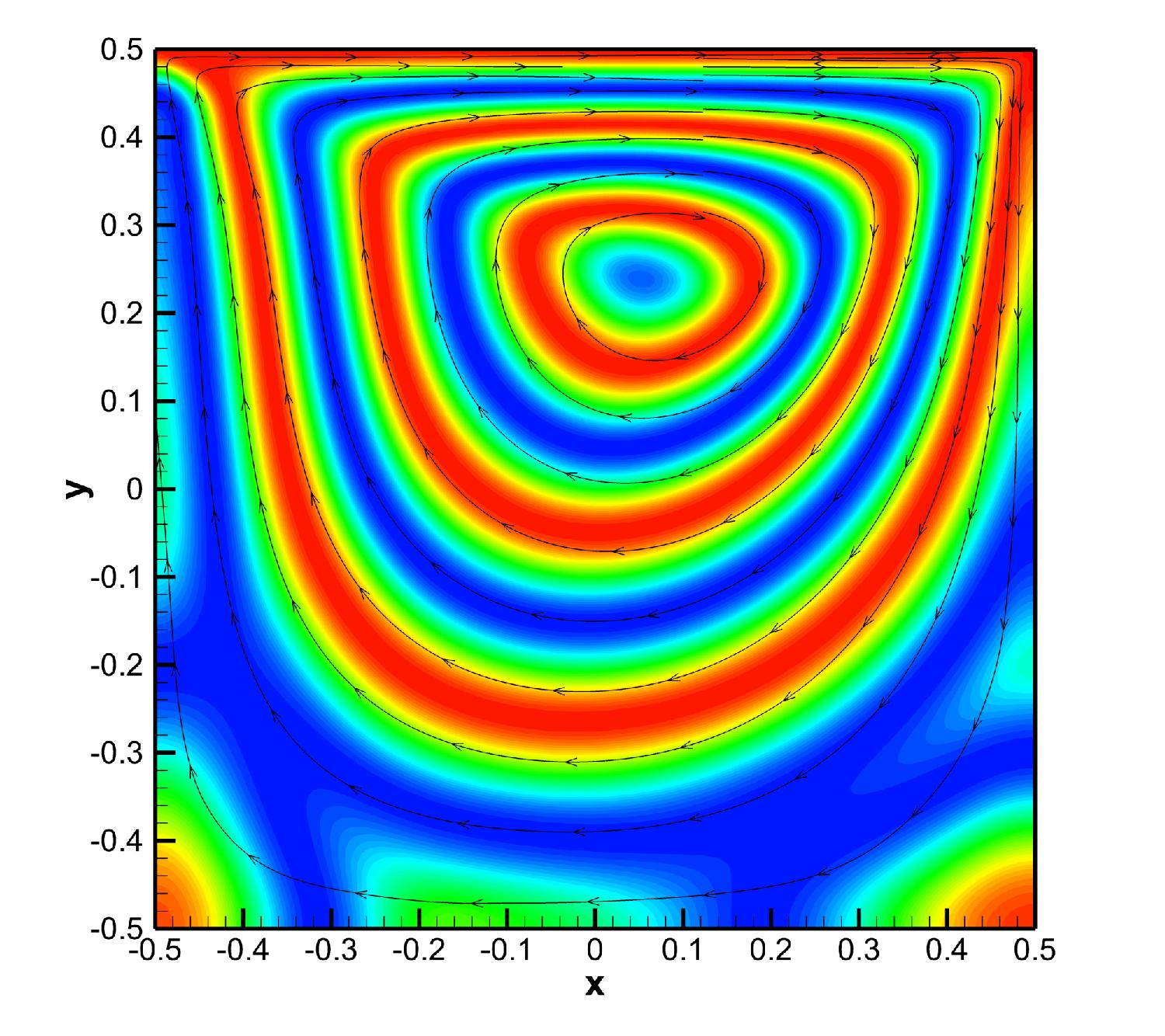} & 
\includegraphics[width=0.45\textwidth]{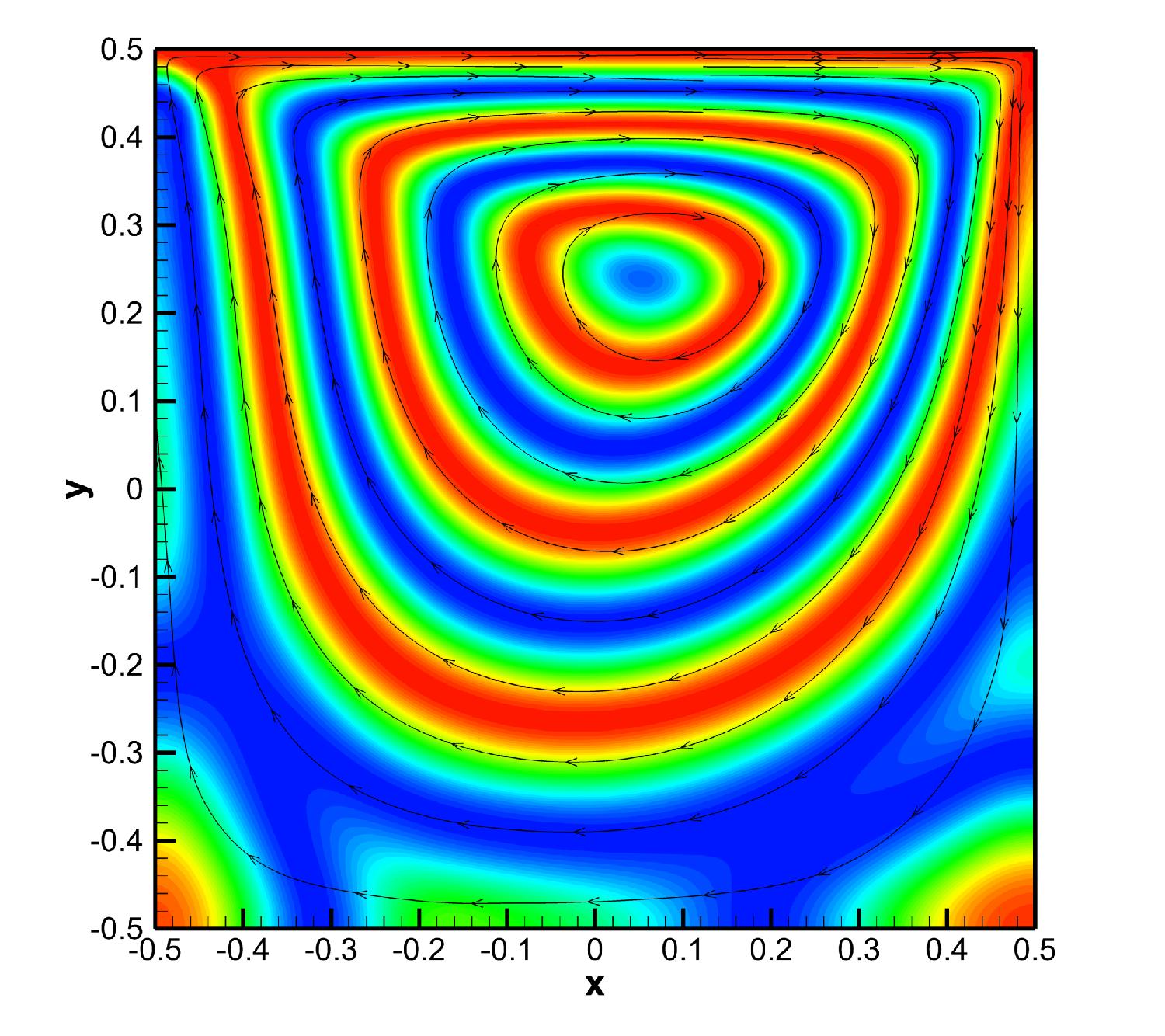}  
\end{tabular} 
\caption{Lid driven cavity at Re=100. Reference solution by Ghia et al. \cite{Ghia1982} and numerical solution for the hyperbolic model of Peshkov and Romenski (HPR) at a final 
time of $t=10.0$ obtained with a third order ADER-WENO finite volume scheme. Streamlines and $u$ velocity contours (top left) and 1D cuts along the $x$ and the $y$ axis (top right). 
The components of the distortion tensor are shown in the bottom row, with 41 equidistant contour colors in the interval [-1,1]: $A_{11}$ (bottom left) and $A_{12}$ (bottom right). } 
\label{fig.cavity}
\end{center}
\end{figure}

The lid-driven cavity is a classical benchmark problem for numerical methods applied to the incompressible Navier-Stokes equations, see \cite{Ghia1982,Tavelli2014}. 
However, it can also be used for compressible flow solvers in the low Mach number regime, see \cite{DumbserCasulli2015}. 
The computational domain is the box $\Omega = [-0.5,0.5] \times [-0.5,0.5]$, which is initialized with a density of $\rho=1$, a velocity of $u=v=0$ and a pressure of 
$p=100/\gamma$. The rest of the parameters is set to $\gamma=1.4$, $c_v = 1$, $c_s = 8$, $\rho_0=1$, $\mathbf{A}=\mathbf{I}$ and $\mathbf{J}=0$. 
The dynamic viscosity is chosen as $\mu=10^{-2}$, while heat conduction is neglected, i.e. $\alpha=\kappa=0$. 
The flow is driven by the upper boundary, whose velocity is set to $\mathbf{v}=(1,0)$. On the other three boundaries, a no-slip wall boundary condition $\mathbf{v}=0$ is imposed. 
We run a third order ADER-WENO finite volume scheme  ($N=0$, $M=2$) on a grid composed of $100 \times 100$ elements until a final time of $t=10$. The reference Mach number of this test case with 
respect to the speed of the lid is $M = 0.1$. The computational results are presented in Fig. \ref{fig.cavity}, where also a comparison with the Navier-Stokes reference solution 
of Ghia et al. \cite{Ghia1982} is shown. We note a very good agreement between the numerical solution of the HPR model and the solution of the incompressible Navier-Stokes 
equations. In the bottom panels of Fig. \ref{fig.cavity}, we plot two components of the distortion tensor $\mathbf{A}$, which is very useful to visualize the main structures 
of the flow. 

\subsection{Double shear layer} 
\label{sec:dsl}

\begin{figure}[!htbp]
\centering 
\begin{tabular}{cc} 
\includegraphics[width=0.45\textwidth]{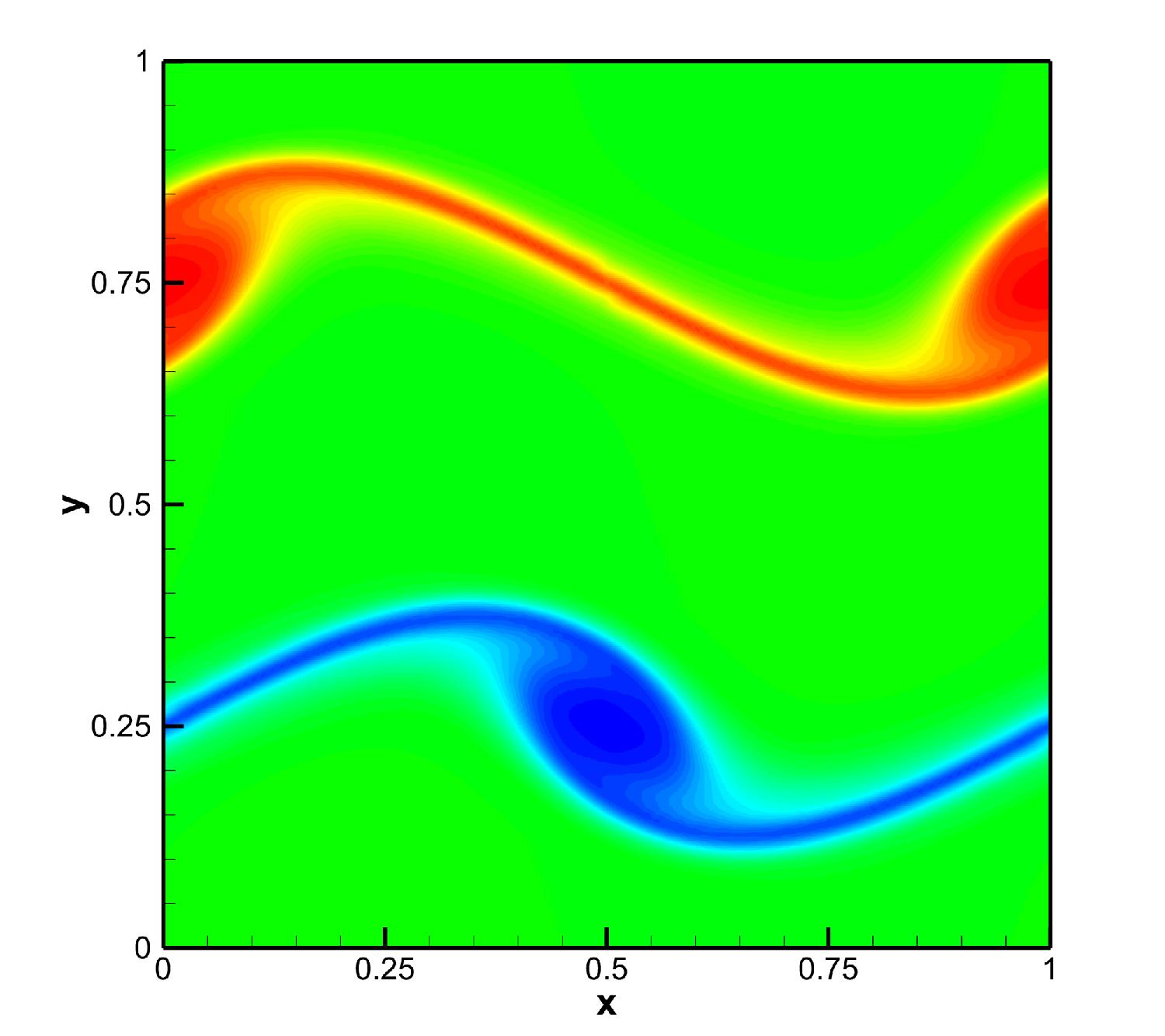}   & 
\includegraphics[width=0.45\textwidth]{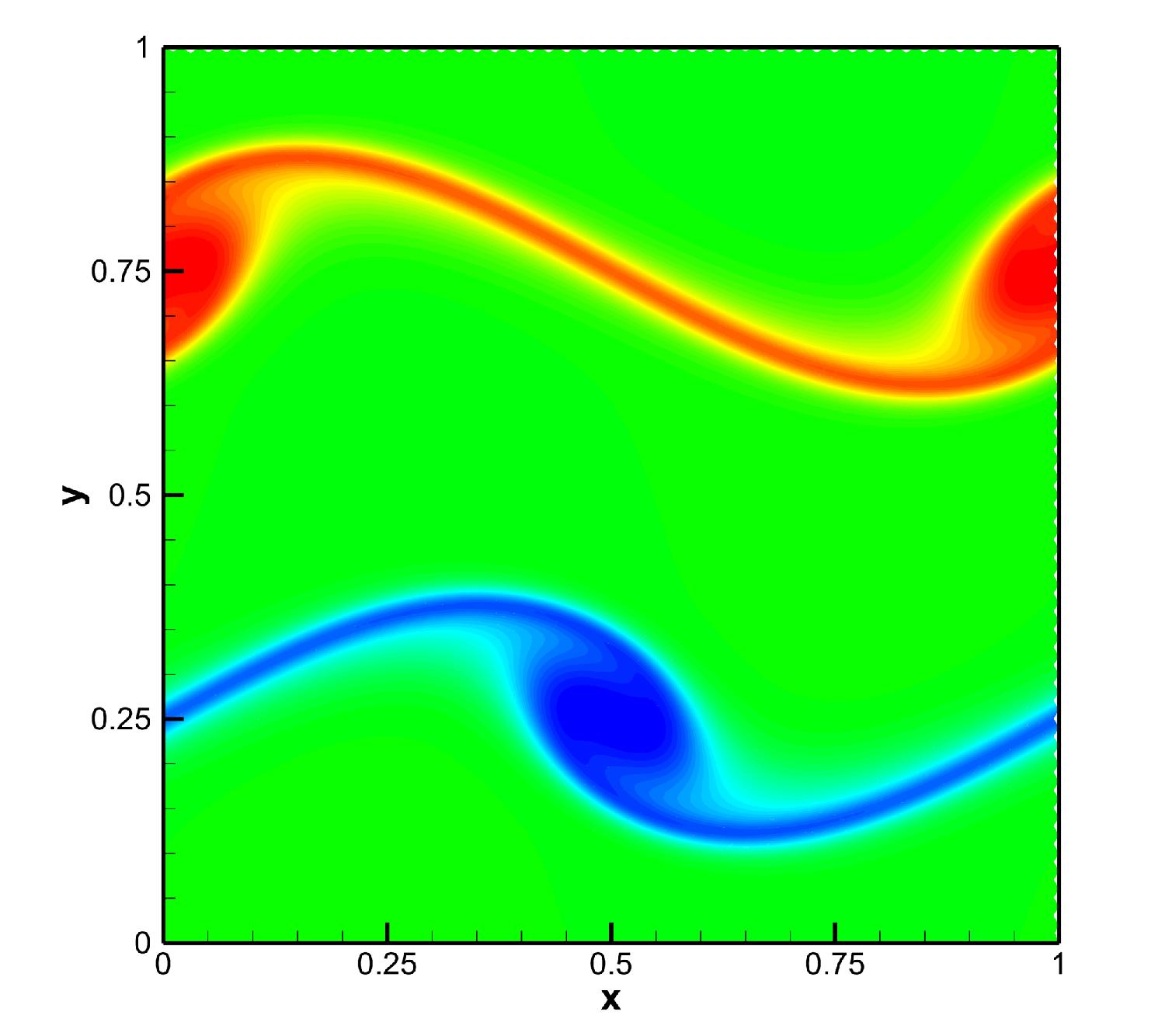}          \\  
\includegraphics[width=0.45\textwidth]{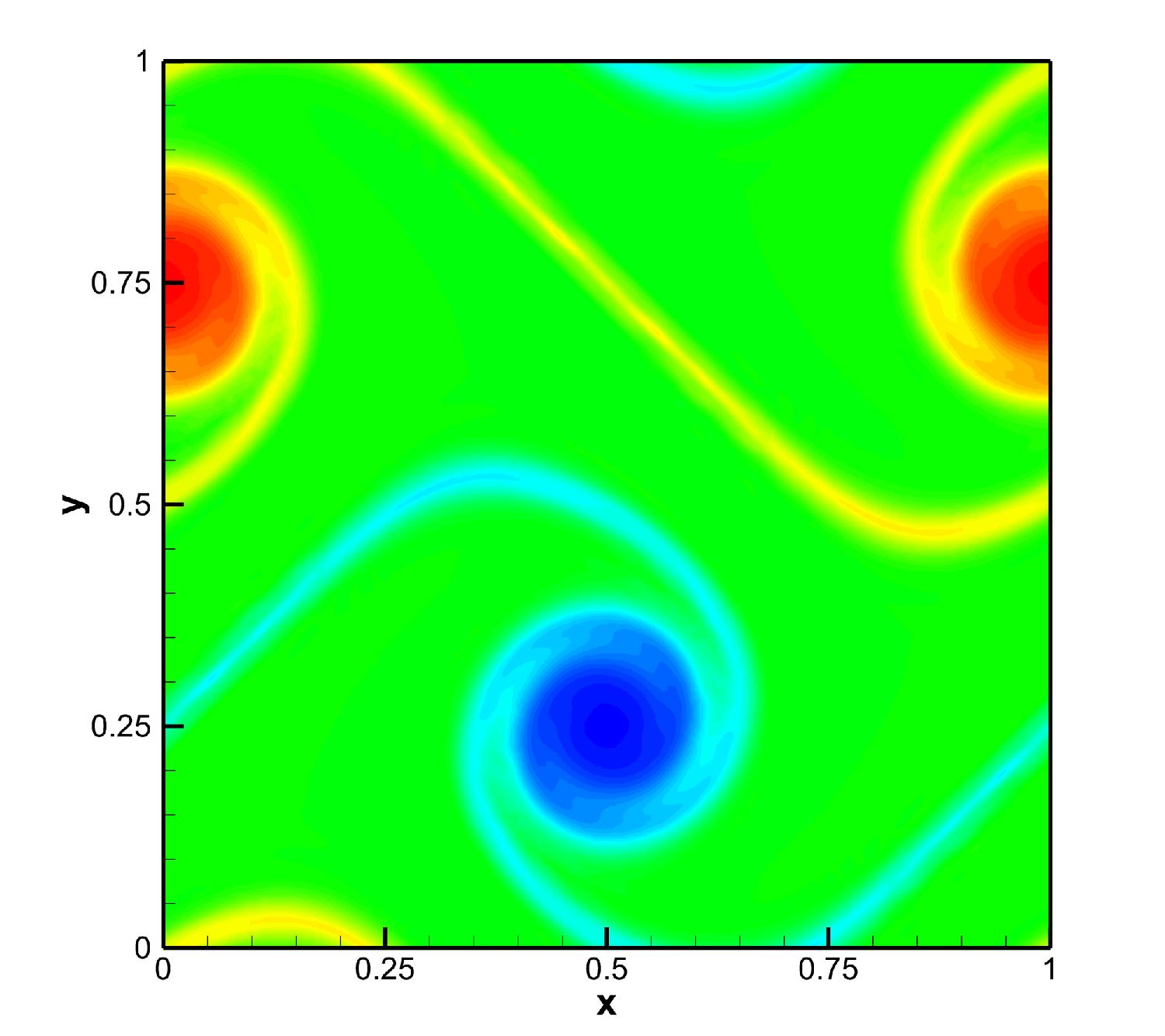}   & 
\includegraphics[width=0.45\textwidth]{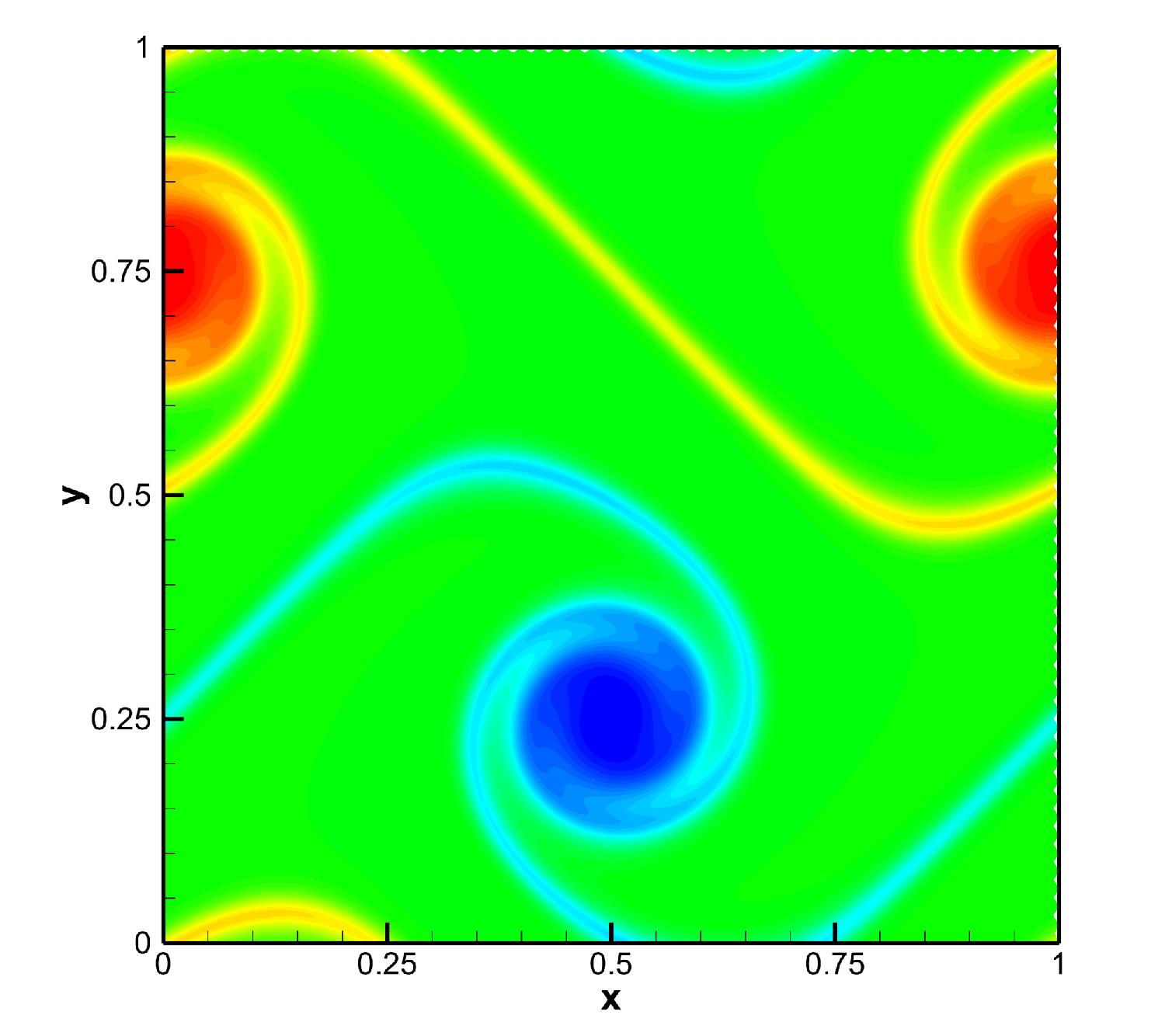}          \\  
\includegraphics[width=0.45\textwidth]{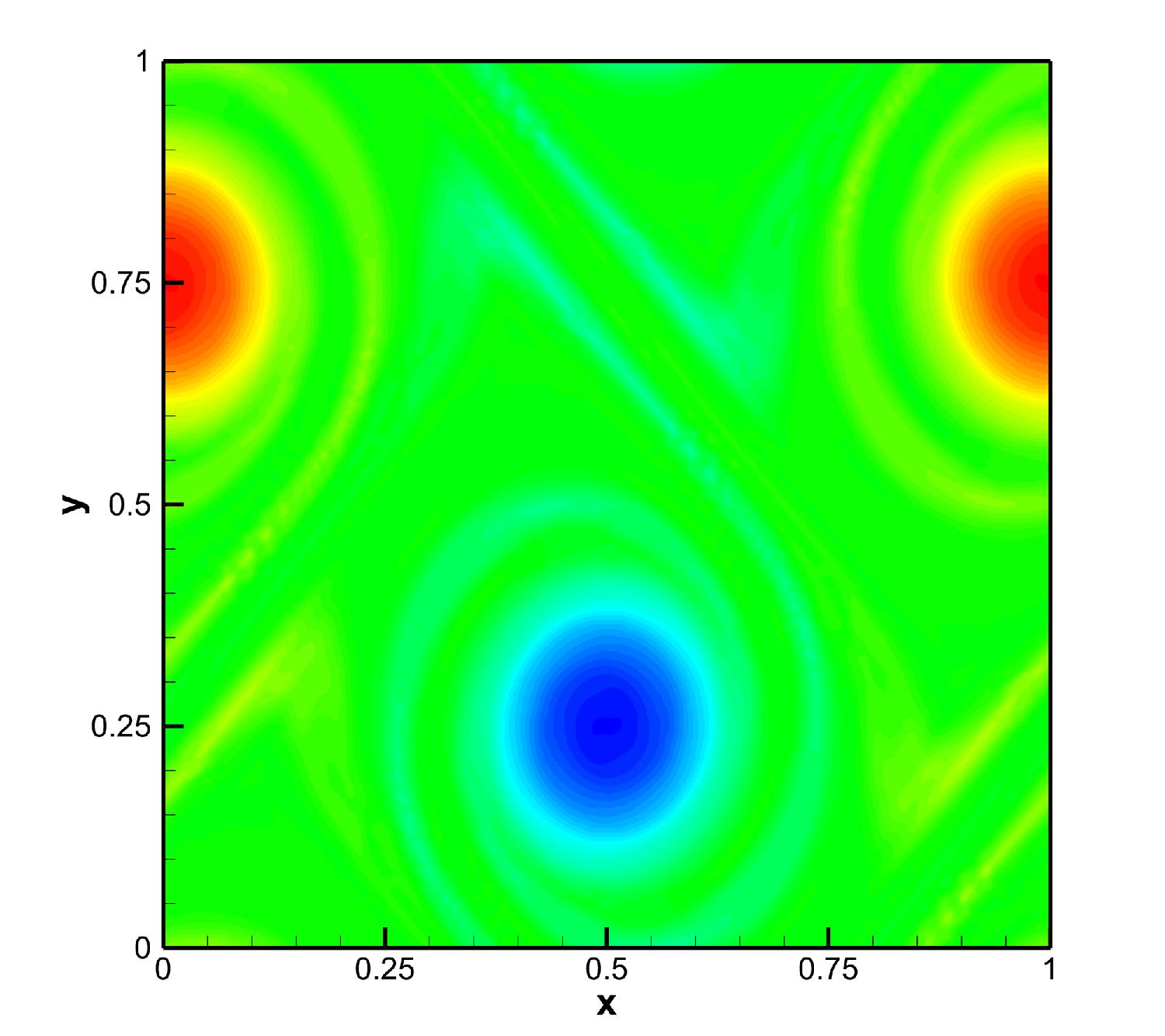}   & 
\includegraphics[width=0.45\textwidth]{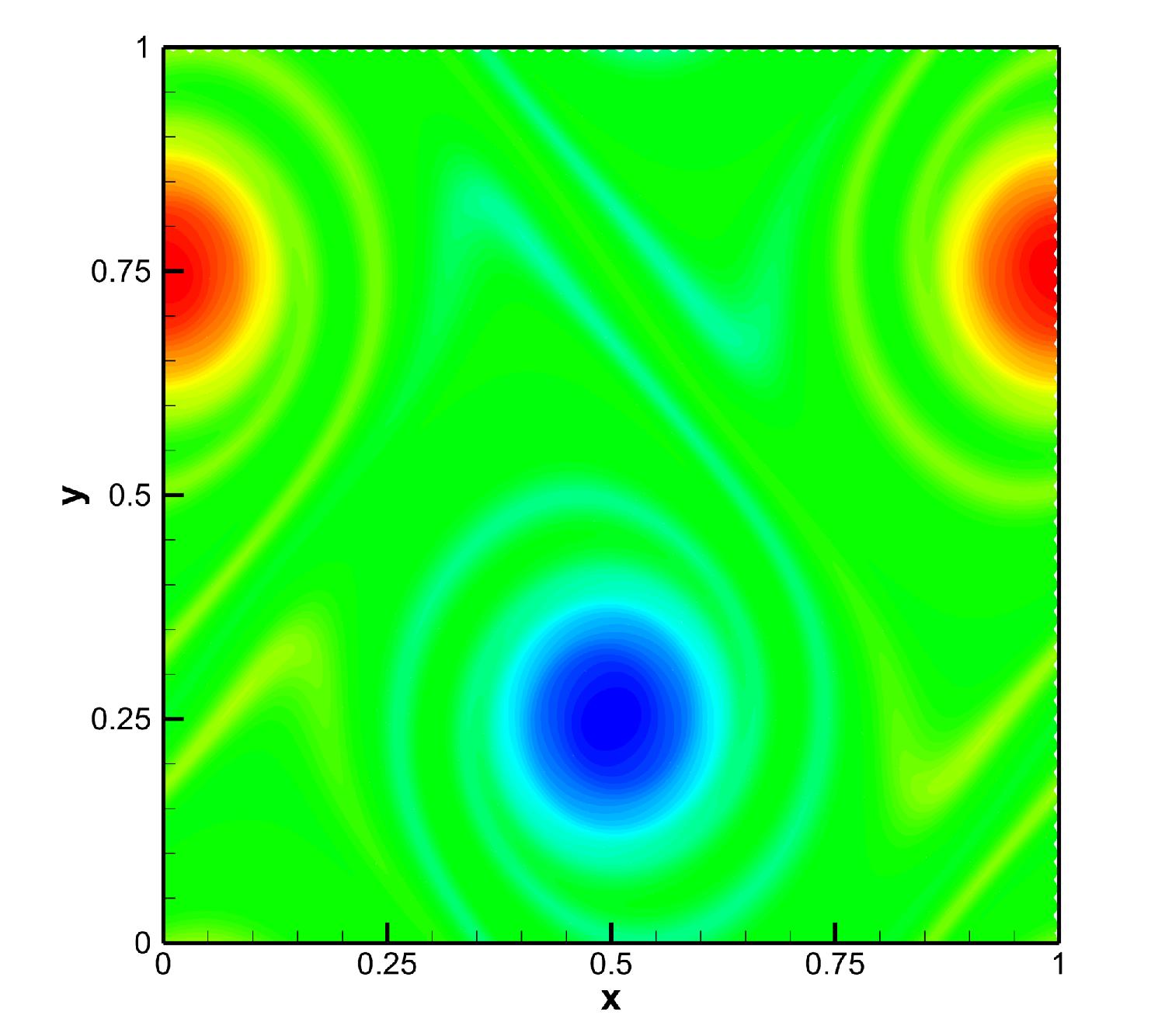}             
\end{tabular} 
\caption{Vorticity contours for the double shear layer with a viscosity of $\nu=2 \cdot 10^{-4}$. Right: numerical solution of the incompressible Navier-Stokes equations obtained with the 
staggered semi-implicit space-time DG scheme of Tavelli and Dumbser \cite{Tavelli2014,Tavelli2015}. Left: numerical solution for the hyperbolic model of  Peshkov  and Romenski (HPR) at times $t=0.8$, $t=1.2$, 
$t=1.8$ from top to bottom obtained with a fourth order ADER-WENO finite volume scheme. }  
\label{fig.dsl}
\end{figure}

\begin{figure}[!htbp]
\begin{center}
\begin{tabular}{cc} 
\includegraphics[width=0.45\textwidth]{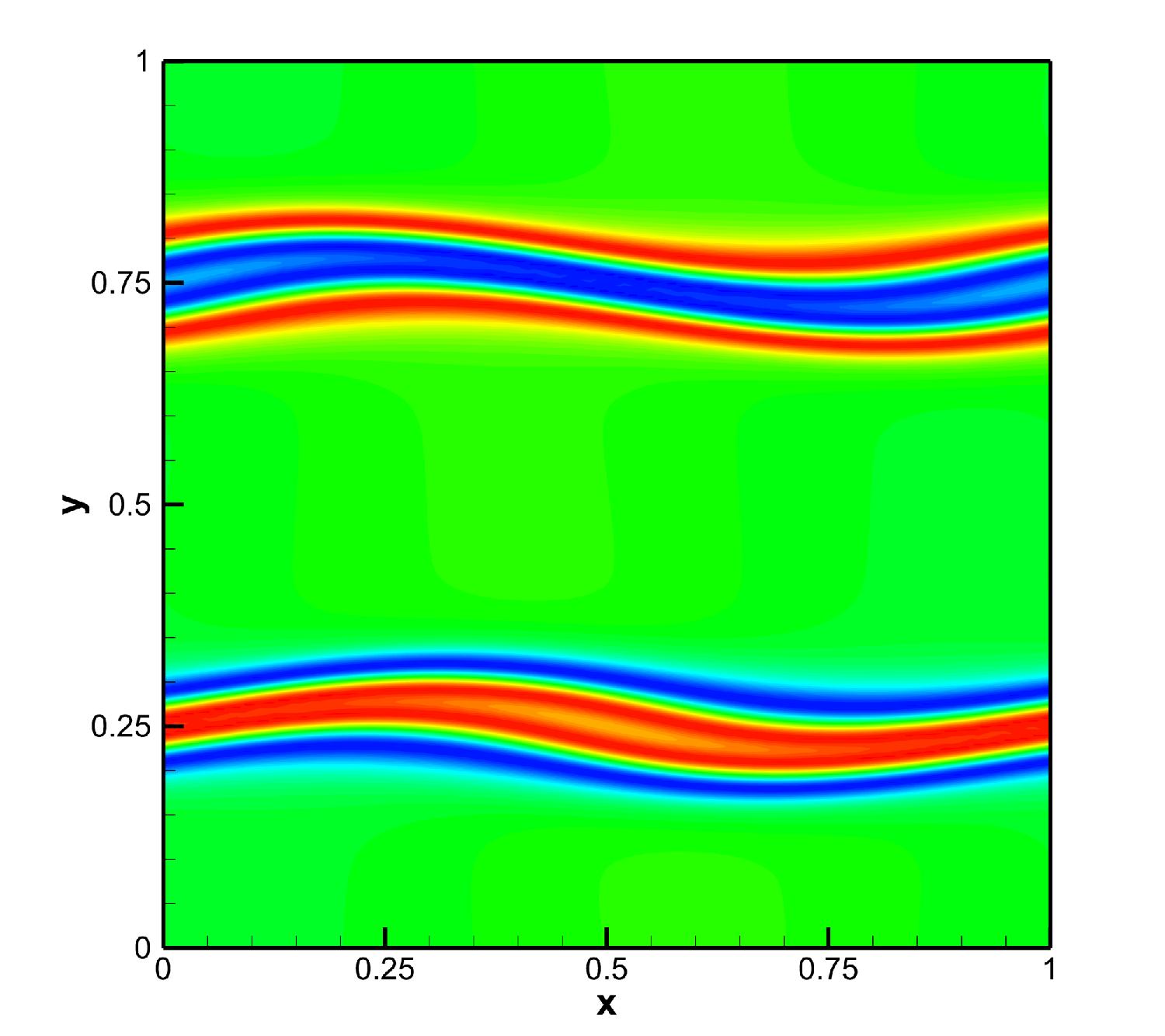}   & 
\includegraphics[width=0.45\textwidth]{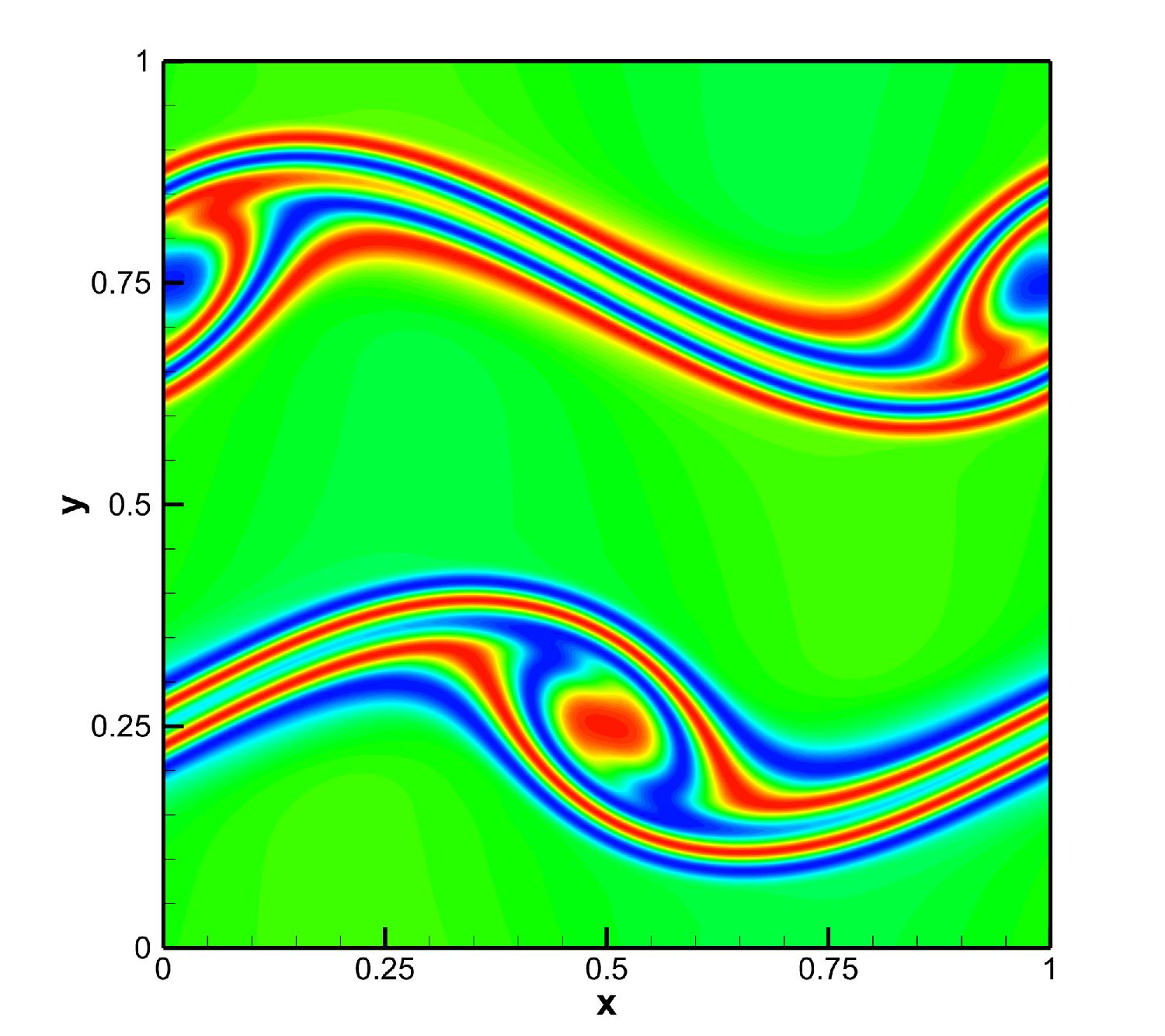}   \\  
\includegraphics[width=0.45\textwidth]{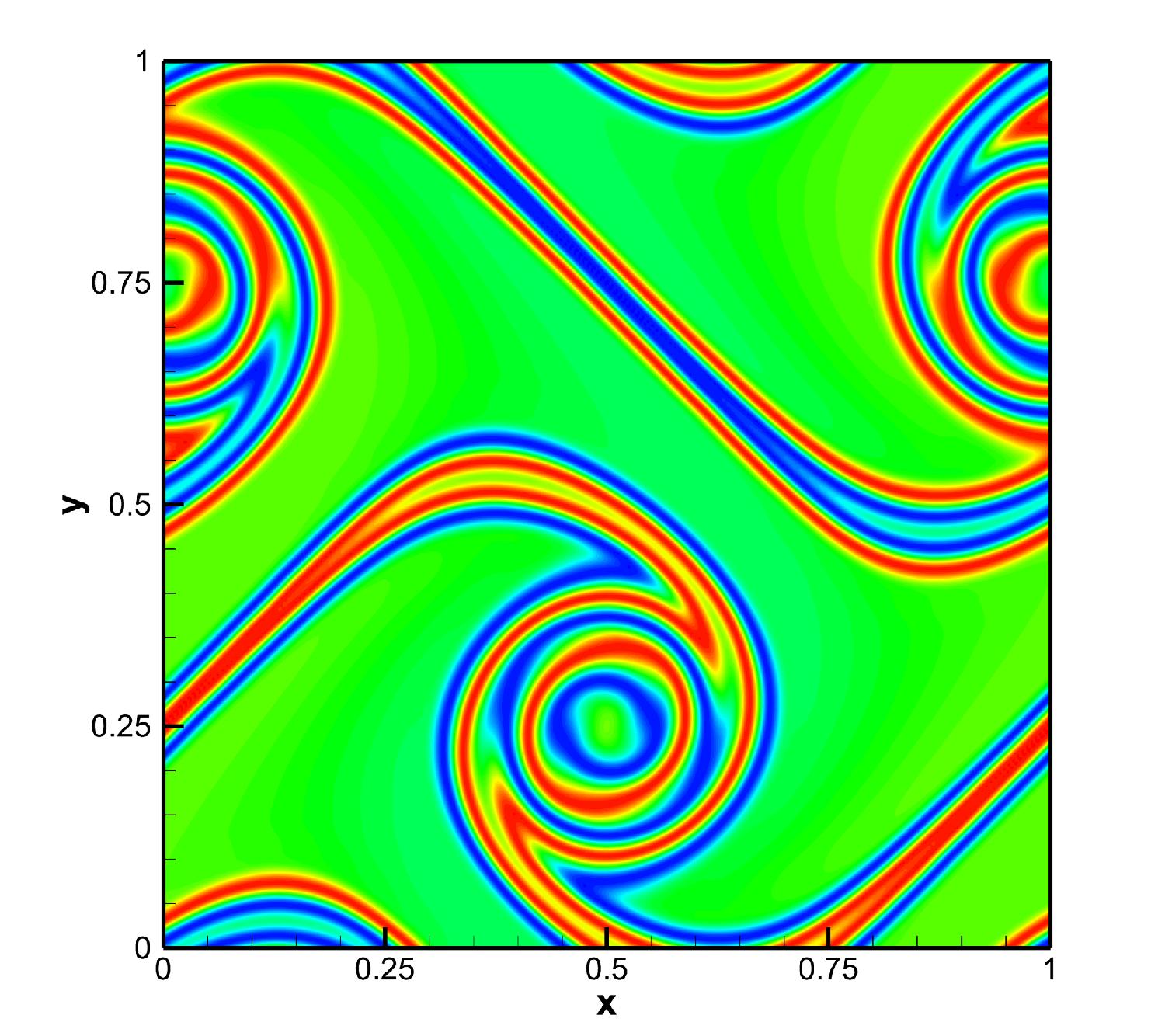}   & 
\includegraphics[width=0.45\textwidth]{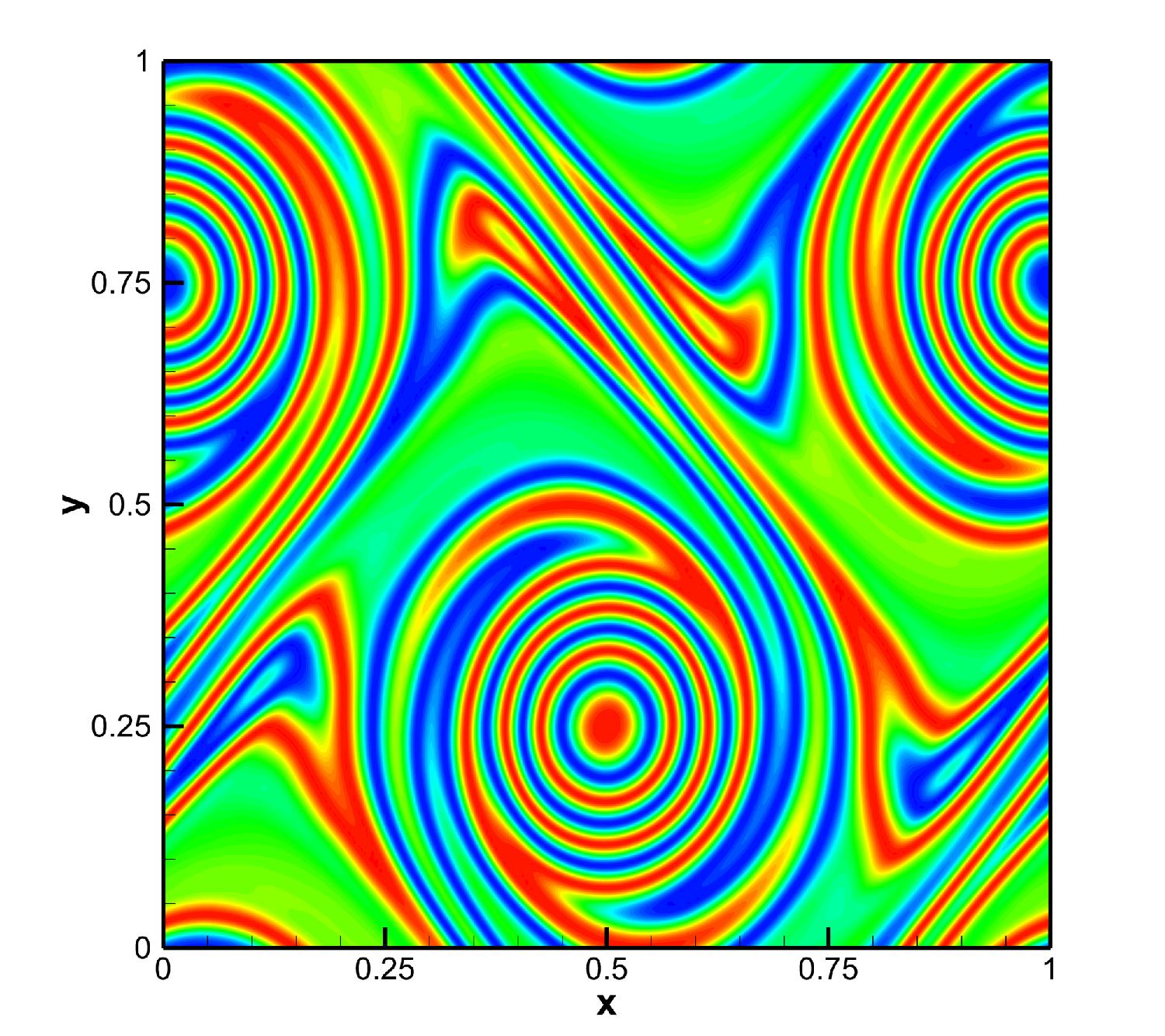}   \\  
\end{tabular} 
\caption{Double shear layer: 41 equidistant contour colors in the interval [-1,1] for component $A_{12}$ of the distortion tensor of the HPR model at 
times $t=0.4$, $t=0.8$, $t=1.2$ and $t=1.8$. } 
\label{fig.dslA12}
\end{center}
\end{figure}

The numerical scheme is now applied to a double shear layer benchmark problem, see \cite{Bell1989,Tavelli2015}, which contains a high initial velocity gradient.  
The computational domain is defined as $\Omega=[0,1]^2$ and periodic boundary conditions are imposed everywhere. As initial condition we consider the 
following perturbed double shear layer profile: 
\begin{equation}
    u=\left\{
    \begin{array}{l}
    	\tanh\left( \tilde{\rho} (y-0.25) \right), \qquad \textnormal{ if } y \leq 0.5, \\
    	\tanh\left( \tilde{\rho} (0.75-y) \right), \qquad \textnormal{ if } y > 0.5,
    \end{array}
    \right.
    \label{eq:DSL_0} 
\end{equation} 
\begin{equation} 		
    v= \delta \sin(2\pi x), \qquad \rho  = 1, \qquad p = \frac{100}{\gamma},
\label{eq:DSL_2}
\end{equation}
where $\tilde{\rho}$ is a parameter that determines the slope of the shear layer; and $\delta$ is the amplitude of the initial perturbation. For the present test we set $\delta=0.05$; $\tilde{\rho}=30$; $\nu= \mu / \rho_0 = 2 \cdot 10^{-4}$. The other parameters are $\gamma = 1.4$, $\rho_0=1$, $c_v=1$, $c_s=8$, $\mathbf{A}=\mathbf{I}$, $\mathbf{J}=0$, $\alpha=\kappa=0$. Simulations are carried out up to a final time of $t=1.8$ using a fourth order ADER-WENO finite volume scheme  ($N=0$, $M=3$) on a grid composed of $200 \times 200$ cells. In Fig. \ref{fig.dsl} the computational results 
obtained with the HPR model are compared with a numerical reference solution based on the solution of the incompressible Navier-Stokes equations. We can note an excellent agreement between the two. 
The reference solution has been obtained with the staggered space-time discontinuous Galerkin finite element method recently proposed in \cite{Tavelli2015}. 
In both models, the two thin shear layers evolve into several vortices, as observed in \cite{Bell1989}, and overall the small flow structures seem to be relatively well resolved also at 
the final time $t=1.8$. In Fig. \ref{fig.dslA12} we plot the time evolution of the distortion tensor component $A_{12}$. One can again observe that the components of the tensor 
$\mathbf{A}$ are excellent candidates for \textit{flow visualization}, since they reveal even more details of the flow structures than the vorticity plotted in the previous Figure \ref{fig.dsl}. 
The main advantage here is that in the framework of the HPR model, the tensor $\mathbf{A}$ is one of the main variables already contained in the state vector $\mathbf{Q}$ of the governing PDE
system, while vorticity needs to be computed from the velocity field via some post-processing technique.

\subsection{Von Karman vortex street} 
\label{sec:karman}

\begin{figure}[!htbp]
\begin{center}
\begin{tabular}{c} 
\includegraphics[width=0.95\textwidth]{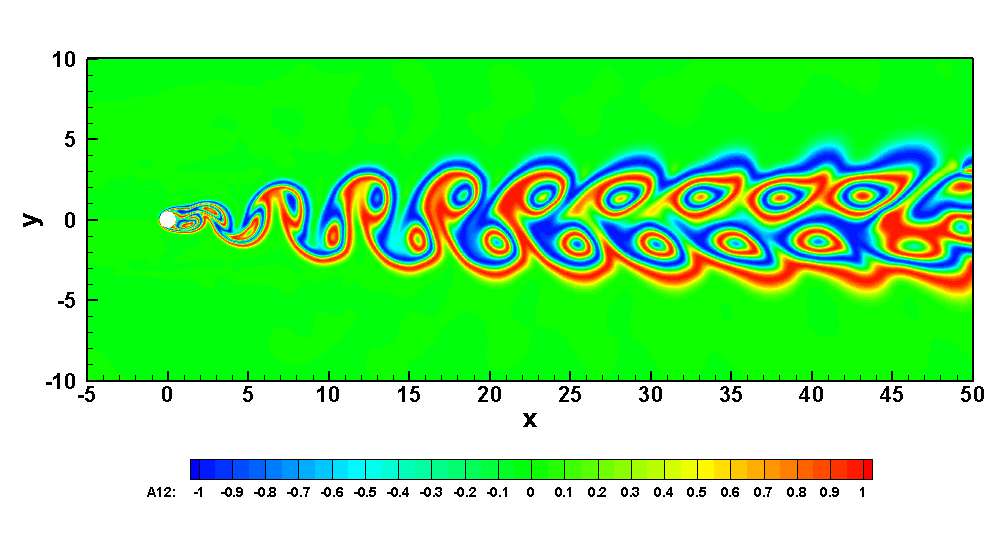} 
\end{tabular} 
\caption{Flow around a circular cylinder ($\Re=150$) at time $t=500$. The quantity $A_{12}$ is shown via 41 equidistant contour levels in the interval [-1,1]. 
The typical von Karman vortex street is clearly visible. }  
\label{fig.karmanstreet}
\end{center}
\end{figure}

\begin{figure}[!htbp]
\begin{center}
\begin{tabular}{cc} 
\includegraphics[width=0.45\textwidth]{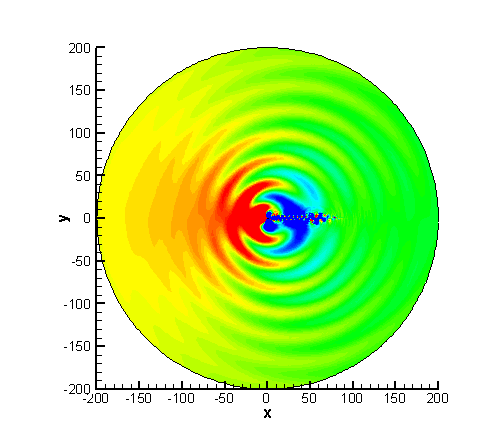} & 
\includegraphics[width=0.45\textwidth]{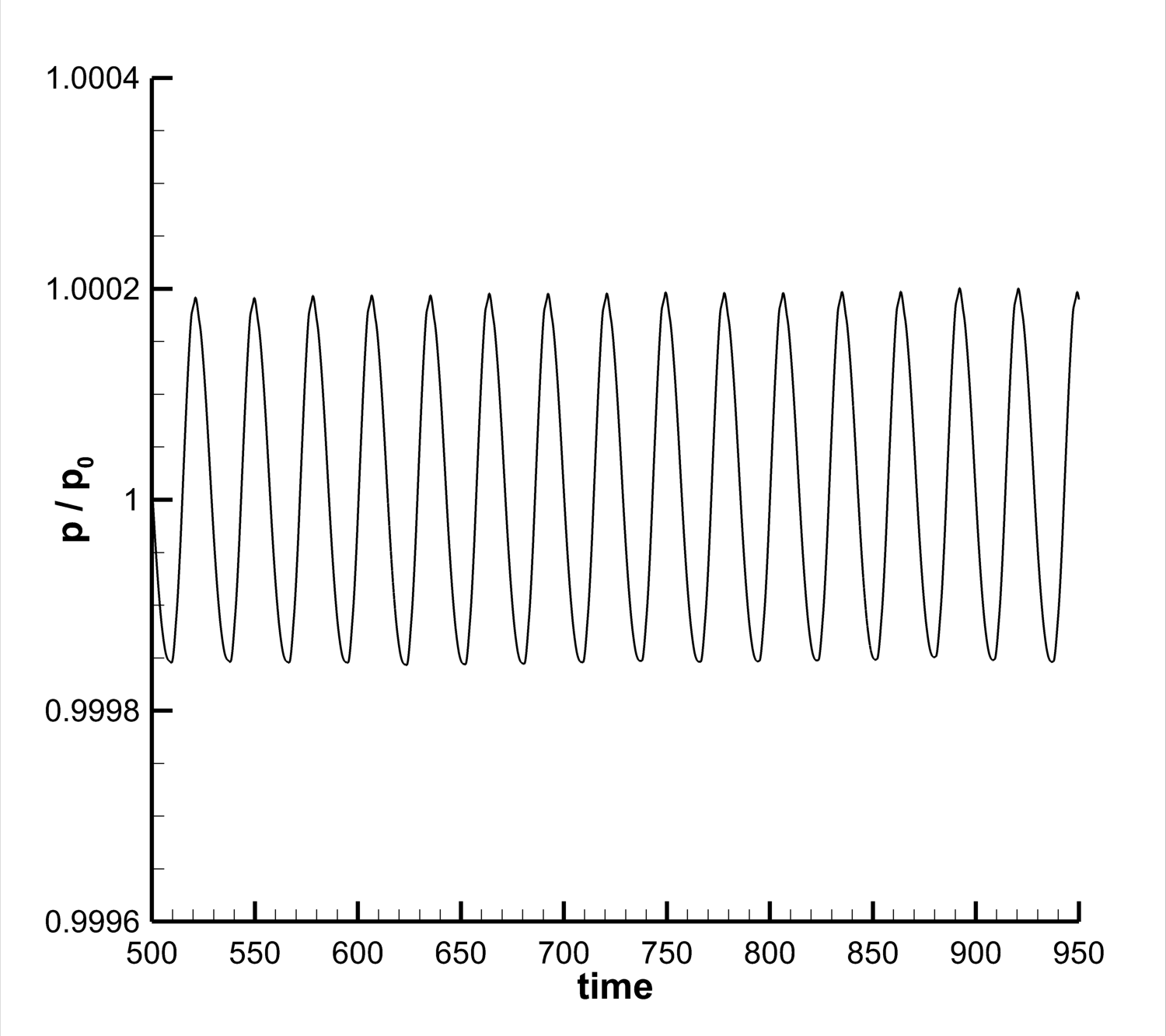}   
\end{tabular} 
\caption{Flow around a circular cylinder at $\Re=150$. Pressure color contours at $t=500$, indicating the sound field generated by the unsteady vortex street (left) 
and sound pressure level (SPL) $p / p_0$ registered at the point $x=(0,50)$ (right). The corresponding Strouhal number of the signal is $\textnormal{St}=0.175$. } 
\label{fig.karmanspl}
\end{center}
\end{figure}

In this section, we solve a test problem used in \cite{Mueller2008} and \cite{ADERNSE} in the context of sound generation 
by a von Karman vortex street that is shed behind a circular cylinder. The circular obstacle has a 
diameter of $d=1$, the Reynolds number based on the diameter is $\Re=\frac{\rho_0 u_0 d}{ \mu } = 150$ and the 
Mach number of the flow is $M_0 = u_0/c_0 = 0.2$. We use $\gamma=1.4$, $c_v=1$, $\rho_0 = 1$, $c_s=0.8$ and $\alpha = \kappa = 0$. 
The initial condition for the HPR model is $\rho=\rho_0$, $u=u_0=0.2$, $v=0$, $p=p_0=1/\gamma$, $\mathbf{A}=\mathbf{I}$ and $\mathbf{J}=0$. 
Computations are performed with a third order ADER WENO finite volume scheme  ($N=0$, $M=2$) on an unstructured triangular mesh 
\cite{DumbserKaeser06b,DumbserKaeser07,DumbserZanotti} that consists of 145,646 triangles. 
The computational domain is a circle with diameter $D=200$ and the simulation is run until a final time of $t=950$. 
A plot of the distortion tensor component $A_{12}$ in the vicinity of the cylinder is depicted in Fig. \ref{fig.karmanstreet}. 
As already stated in the previous examples, the tensor $\mathbf{A}$ turns out to be very useful for \textit{flow visualization}, since all
details of the flow structure can be easily recognized in Fig. \ref{fig.karmanstreet}. An instantaneous plot at $t=500$ of the sound pressure 
field generated by this unsteady flow is illustrated in Fig. \ref{fig.karmanspl} and agrees reasonably well with those computed in 
\cite{Mueller2008,ADERNSE}. The time history of the sound pressure level (SPL) at the point $x=(0,50)$ is also presented in Fig. \ref{fig.karmanspl}. 
From the analysis of the SPL signal of our numerical simulations we obtain a Strouhal number of $\textnormal{St}=\frac{f d}{u_0} = 0.175$, which is in 
reasonable agreement with the value of $\textnormal{St}=0.183$ found by M\"uller \cite{Mueller2008} and with the value of $\textnormal{St}=0.182$ reported in 
\cite{ADERNSE}.

\subsection{Compressible mixing layer} 
\label{sec:mixinglayer}

\begin{figure}[!htbp]
\begin{center}
\begin{tabular}{l} 
\includegraphics[width=0.7\textwidth]{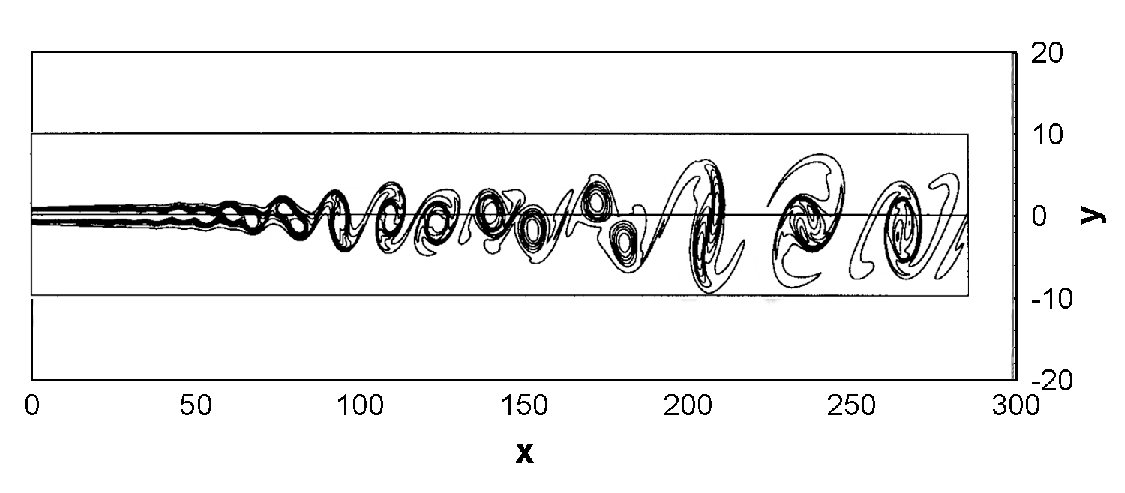}       \\ 
\includegraphics[width=0.7\textwidth]{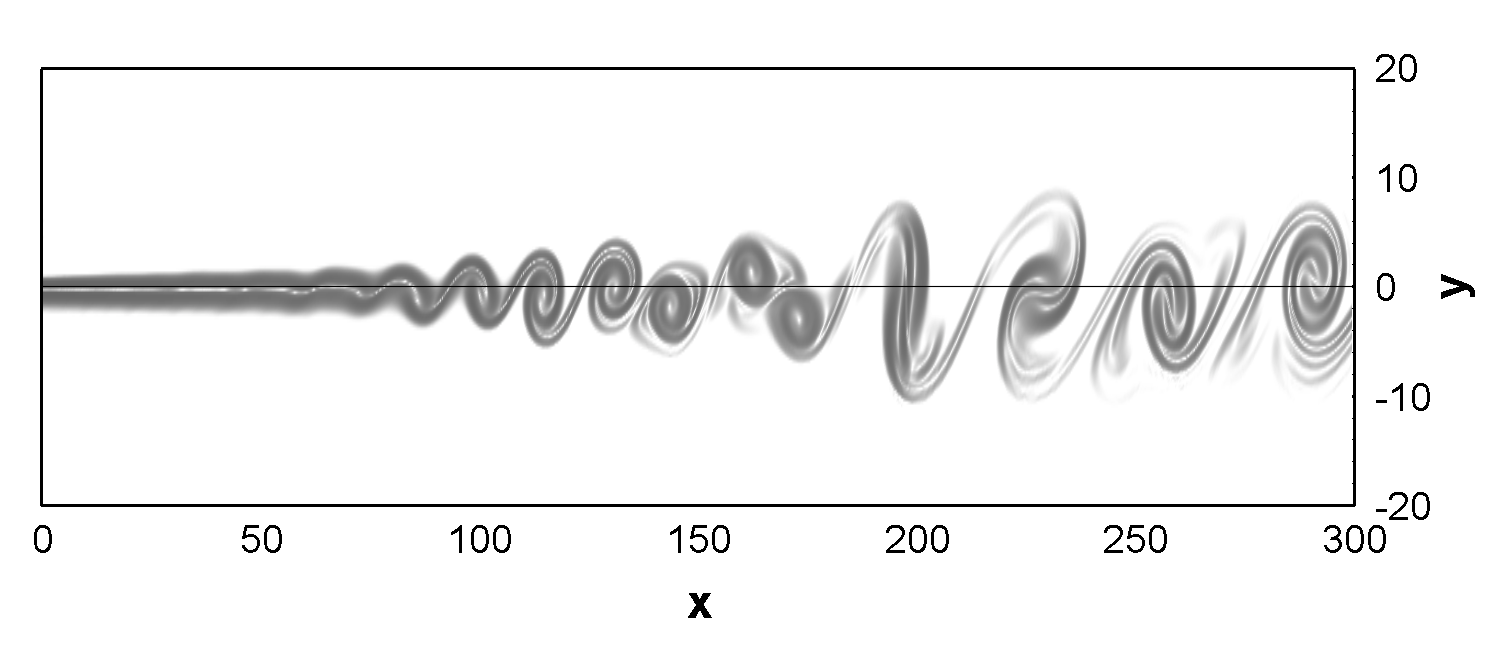} \\  
\includegraphics[width=0.68\textwidth]{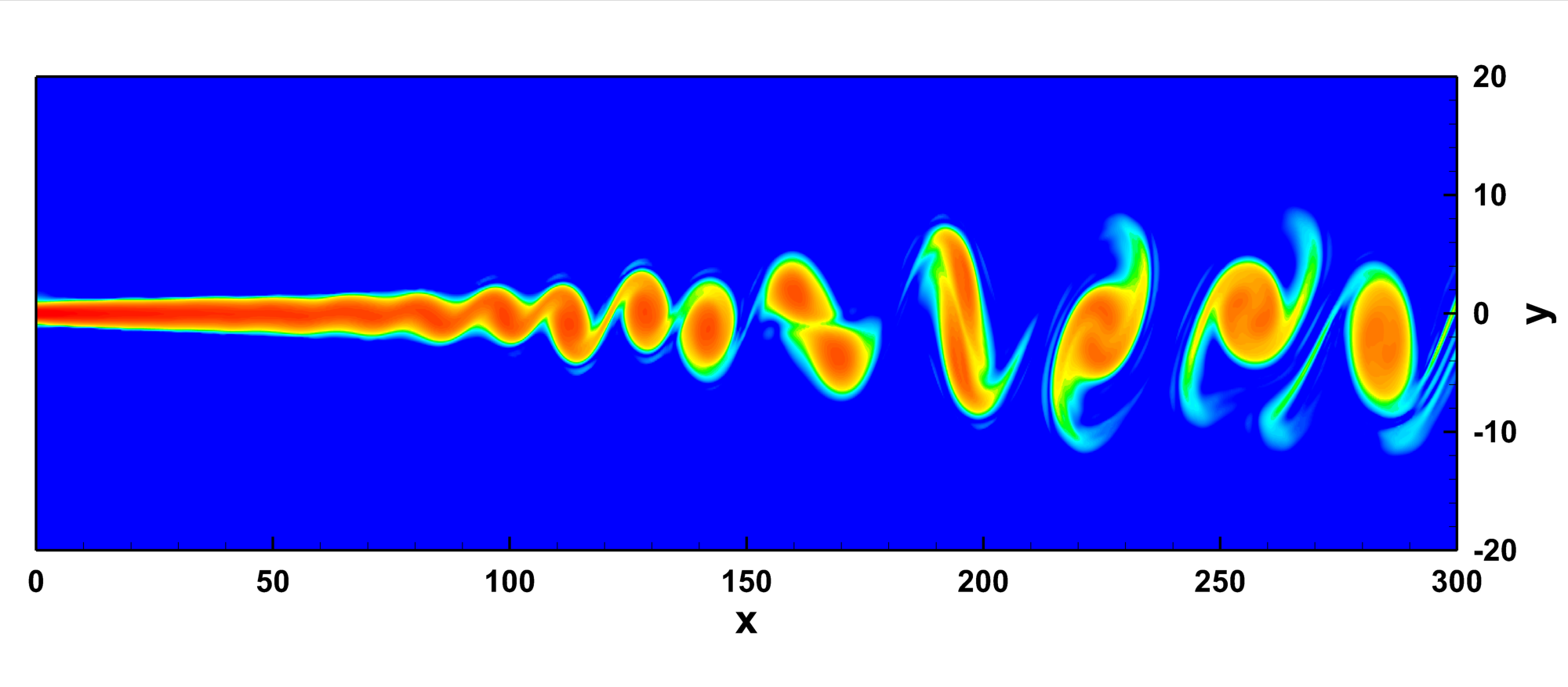}            \\ 
\includegraphics[width=0.68\textwidth]{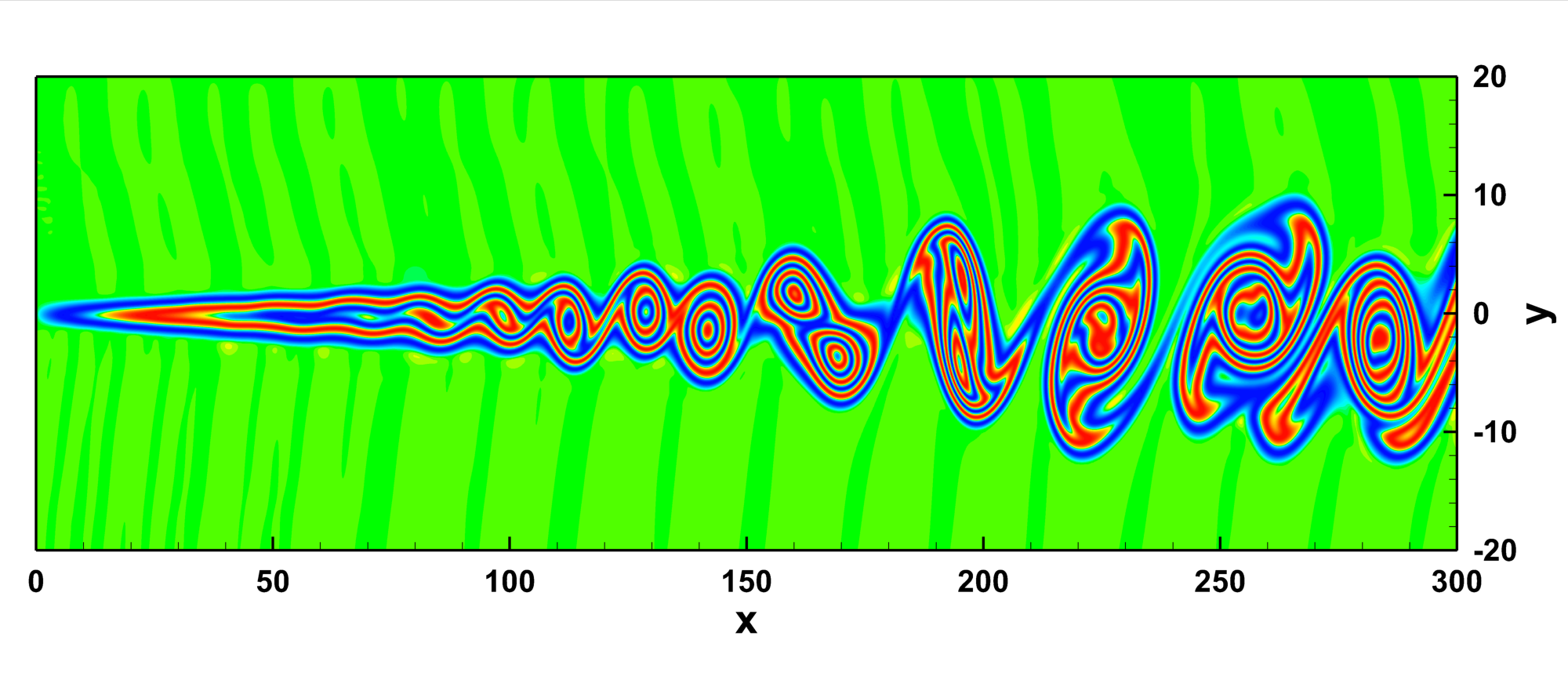}    
\end{tabular} 
\caption{Compressible mixing layer at $Re=500$ at a final time of $t=1596.8$. The reference solutions have been obtained by solving the compressible Navier-Stokes equations with a
compact finite difference scheme, see Colonius and Moin \cite{colonius} (top row) and with a sixth order $P_NP_M$ scheme ($P_3P_5$), see Dumbser \cite{ADERNSE} (second row). 
The numerical solution obtained with a third order ADER-WENO finite volume scheme ($N=0$, $M=2$) for the hyperbolic model of Peshkov and Romenski (HPR) is shown at the bottom. 
In rows 1-3 the vorticity magnitude is plotted, while in the fourth row (for the HPR solution) the quantity $A_{12}$ is shown via 41 equidistant contour levels in the interval [-1,1].  } 
\label{fig.mix}
\end{center}
\end{figure}

Next, we simulate the behavior of an unsteady compressible mixing layer, following partially the setup presented in \cite{colonius,babucke,stedg2,ADERNSE}. 
For $y \to \infty$ the limit of the axial flow velocity is $u_\infty = 0.5$, while for $y \to -\infty$ the axial velocity component tends to $u_{-\infty} = 0.25$. 
Here, we use the simplified initial condition for the flow velocity 
\begin{equation}
  u_0 = \frac{1}{8} \tanh \left( 2 y \right) + \frac{3}{8}, \qquad v_0 = 0, 
\end{equation} 
while the initial condition for the other flow variables is simply given by $\rho = 1$, $p=p_0=1/\gamma$, $\mathbf{A}=\mathbf{I}$ and $\mathbf{J}=0$. 
The vorticity thickness at the inflow, with respect to which all lengths are made dimensionless, is 
\begin{equation}
  \theta = \frac{u_\infty-u_{-\infty}}{ \max{ \left( \left. \frac{\partial u}{\partial y} \right|_{x=0} \right)} } := 1,
\end{equation}
and the Reynolds number based on this vorticity thickness is 
\begin{equation}
   \textnormal{Re}_{\theta} = \frac{\rho_0  u_\infty \theta}{\mu} = 500. 
\end{equation}
Thus, the parameters of the HPR model are set to $\gamma=1.4$, $c_v=2.5$, $\rho_0=1$, $c_s=0.8$, $\mu = 10^{-3}$, $\alpha = \kappa = 0$, i.e. heat conduction is again neglected in this example.  
At the inflow ($x=0$), the flow quantities $\rho$, $u$, $v$ and $p$ are perturbed as follows: $\rho(0,y,t) = \rho_0 + 0.05 \, \delta$, 
$u(0,y,t) = u_0 + \delta$, $v(0,y,t) = v_0 + 0.6 \, \delta$ and $p(0,y,t) = p_0 + 0.2 \, \delta$, with 
\begin{equation}
  \delta = -10^{-3} \exp( -0.25 y^2)   
                          \left(    \cos \left( \omega t \right) 
													        + \cos \left( \halb \omega t - 0.028 \right) 
																	+ \cos \left( \frac{1}{4} \omega t + 0.141 \right) 
																	+ \cos \left( \frac{1}{8} \omega t + 0.391 \right) \right) 
\end{equation} 
and the fundamental frequency of the mixing layer $\omega = 0.3147876$. The computational domain is defined by $\Omega = [0,400] \times [-50, 50]$ and is covered by a 
Cartesian grid of $1600 \times 800$ elements. The final simulation time is set to $t=1596.8$. The computational results obtained with the HPR model are
depicted in Fig. \ref{fig.mix}, where also two numerical reference solutions based on the solution of the compressible Navier-Stokes equations, see \cite{colonius,ADERNSE}, 
are shown. Overall, we observe a reasonable qualitative agreement between the HPR model and the Navier-Stokes reference solution. In the last row of Fig. 
\ref{fig.mix} we show again a component of the tensor $\mathbf{A}$, which seems to reveal the vortex flow structures of the mixing layer even better and in a clearer 
way than the usual vorticity. This underlines again the potential of using $\mathbf{A}$ for the purpose of \textit{flow visualization}. 

\subsection{2D Taylor-Green vortex}  
\label{sec:tgv}

Another typical test problem used for the verification of numerical methods for the incompressible Navier-Stokes equations is the Taylor-Green vortex problem, which is another one 
of the rare examples where an exact analytical solution of the unsteady Navier-Stokes equations is known. In two space dimensions, the solution reads 
\begin{eqnarray}
    u(x,y,t)&=&\sin(x)\cos(y)e^{-2\nu t}, \label{eq:TG_0} \\
    v(x,y,t)&=&-\cos(x)\sin(y)e^{-2\nu t}, \label{eq:TG_1} \\
    p(x,y,t)&=& C + \frac{1}{4}(\cos(2x)+\cos(2y))e^{-4\nu t}.
\label{eq:TG_2}
\end{eqnarray}
The computational domain is $\Omega=[0,2\pi]^2$ with four periodic boundaries everywhere. We carry out the numerical simulations based on the HPR model up to a final time of $t=10$, 
using a fourth order ADER-DG $P_3P_3$ scheme ($N=M=3$) on a computational grid composed of $50 \times 50$ elements. For the HPR model, the following set of parameters has been chosen: $\gamma=1.4$, 
$\rho_0=1$, $\mu=10^{-2}$, $c_v=1$, $c_s=10$, $\alpha=\kappa=0$. The initial conditions for the velocity and the pressure are given by \eqref{eq:TG_0}-\eqref{eq:TG_2}, where the additive 
constant in the pressure field is set to $C=100/\gamma$. The distortion tensor and the heat flux are initialized as usual with $\mathbf{A}=\mathbf{I}$ and $\mathbf{J}=0$. 

\begin{figure}[!htbp]
\begin{center}
\begin{tabular}{cc} 
\includegraphics[width=0.45\textwidth]{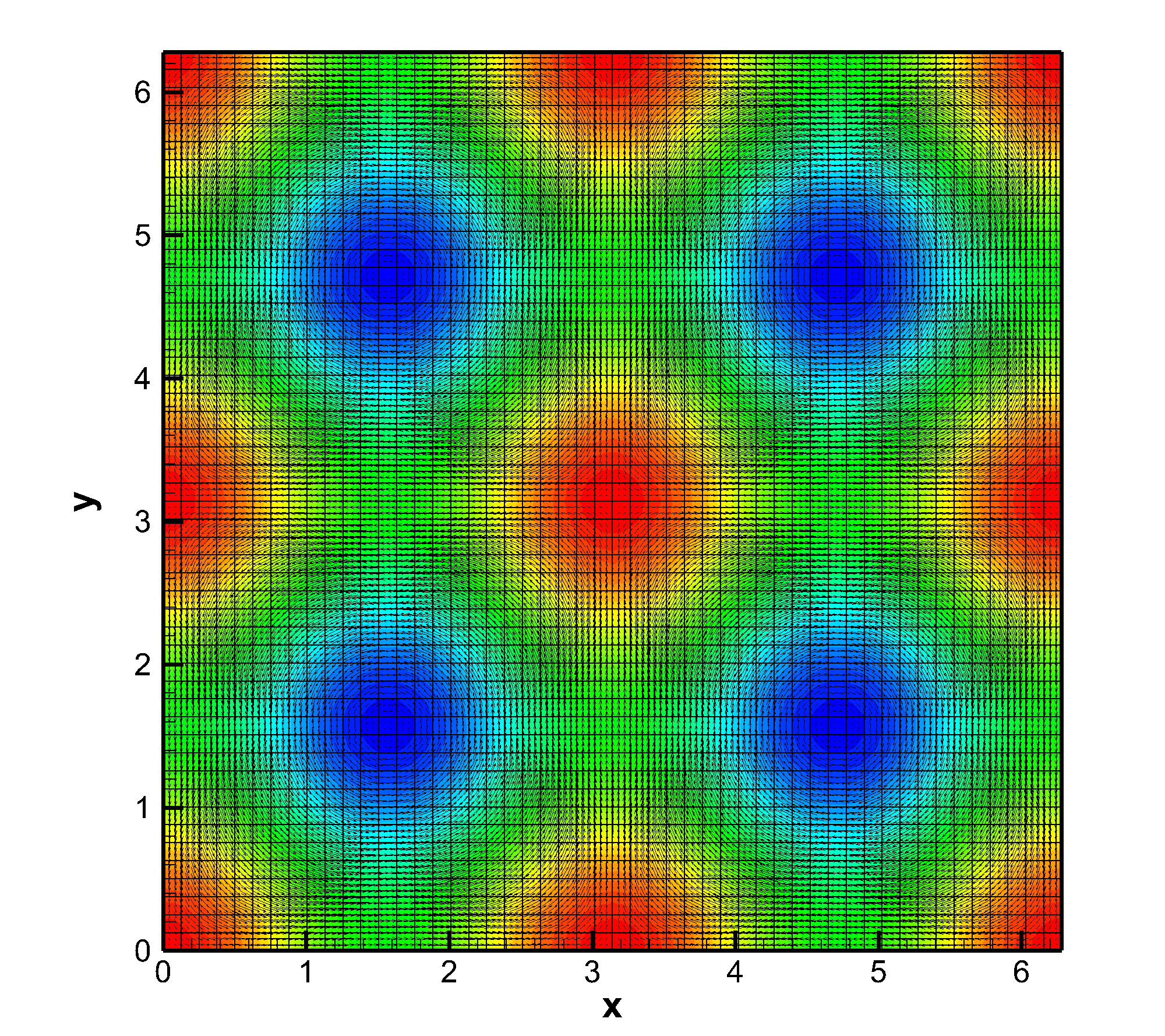}   & 
\includegraphics[width=0.45\textwidth]{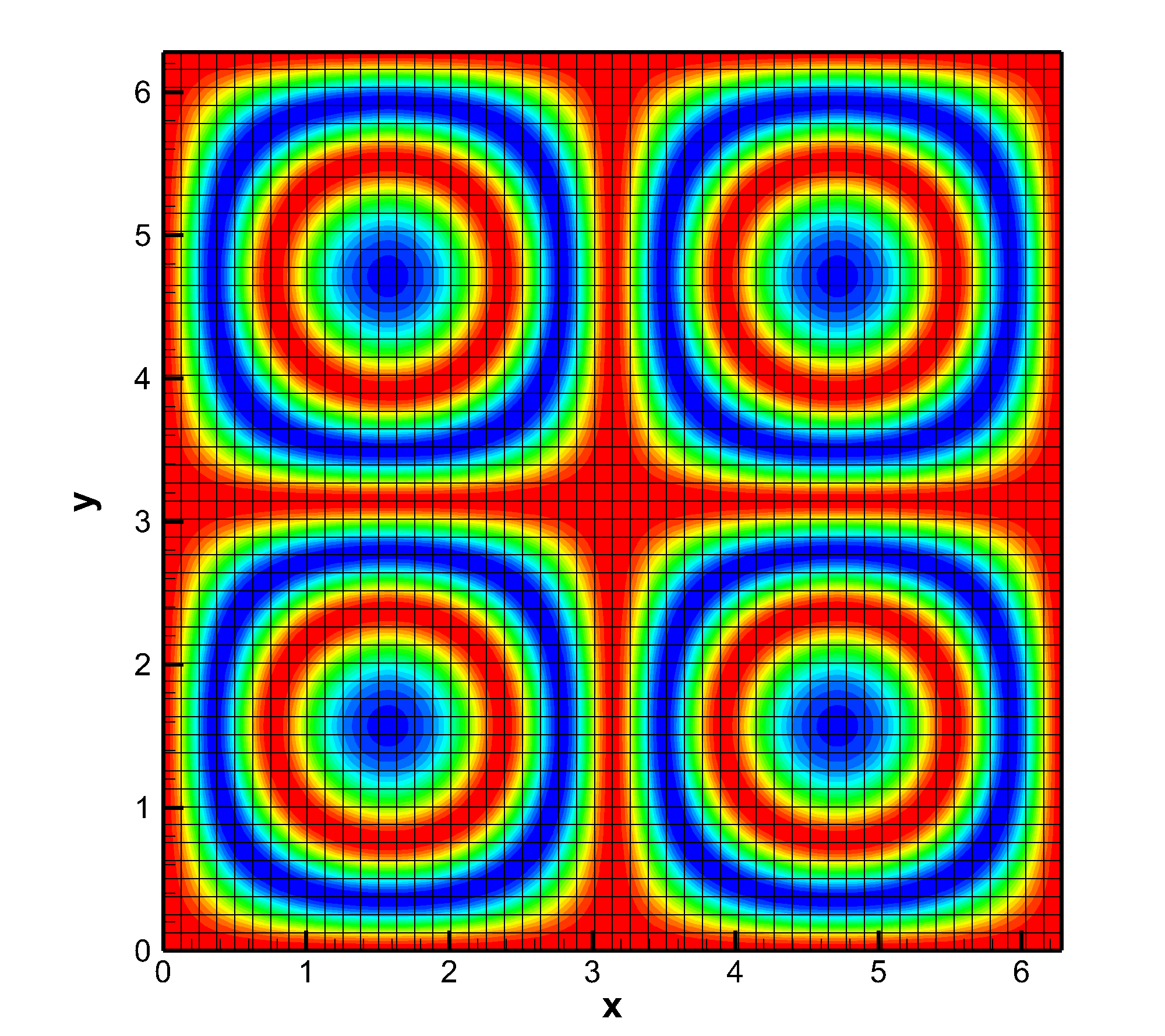} \\  
\includegraphics[width=0.45\textwidth]{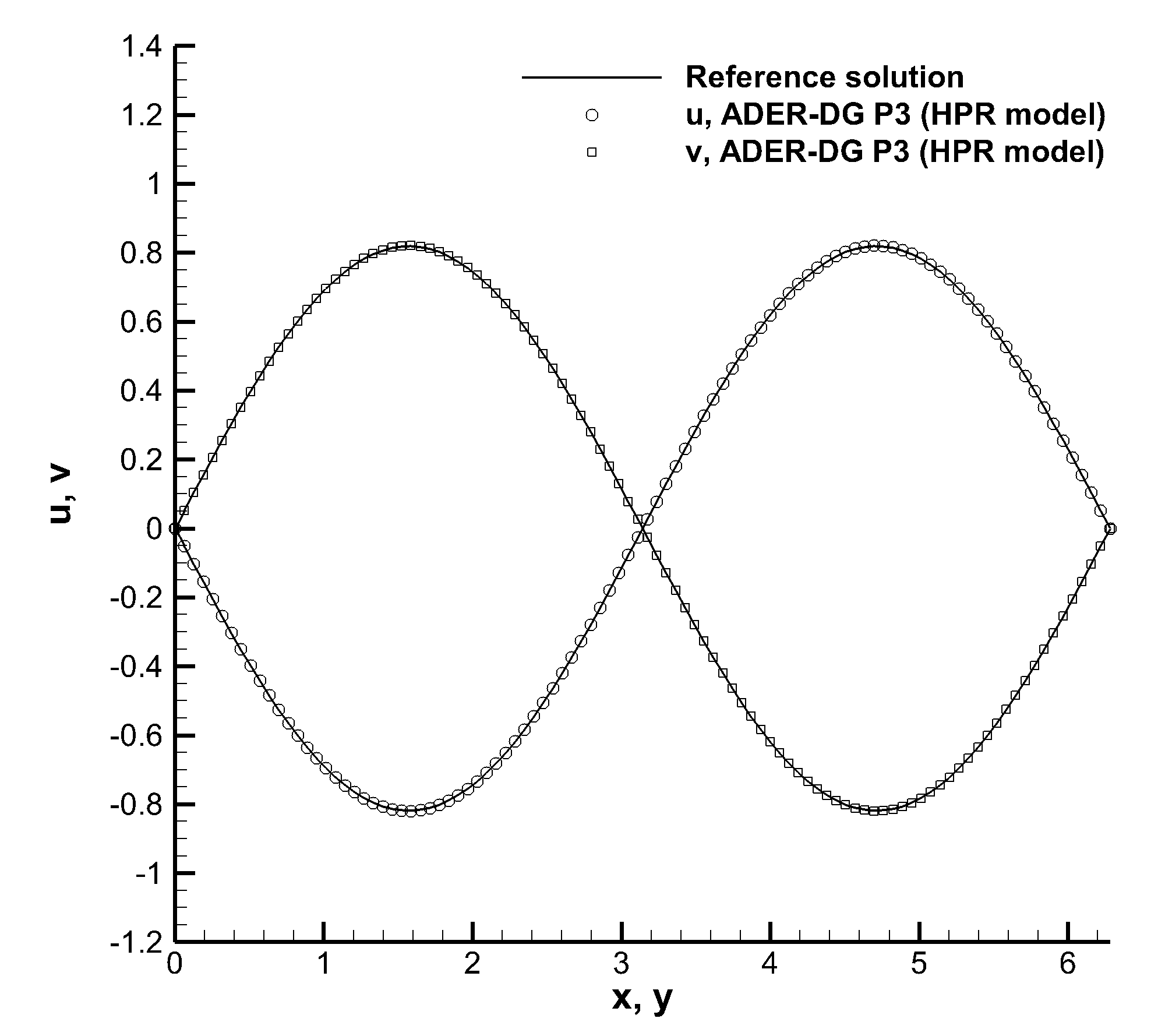}   & 
\includegraphics[width=0.45\textwidth]{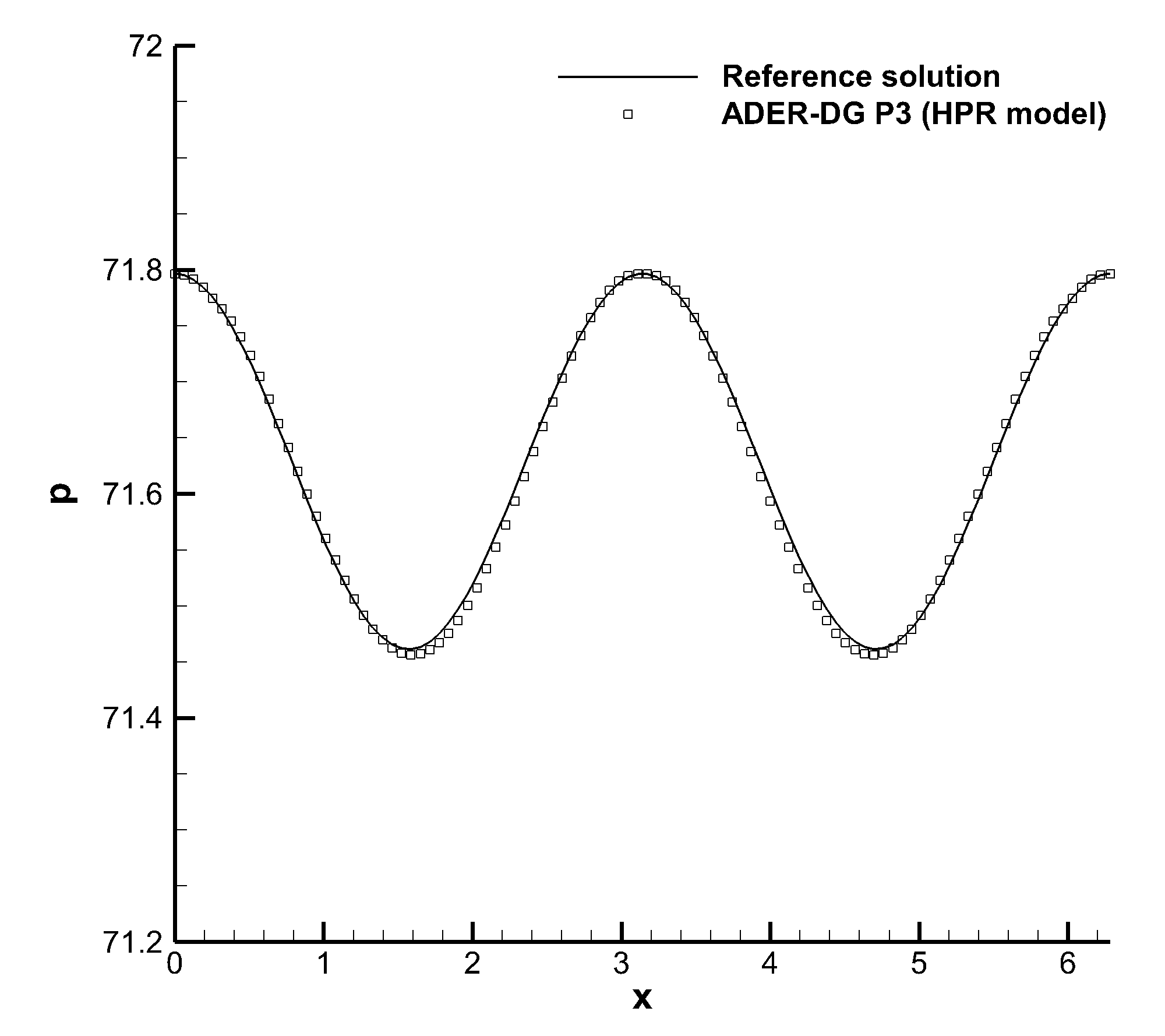}    
\end{tabular} 
\caption{Taylor-Green vortex with a viscosity of $\nu=10^{-2}$. Exact solution of the Navier-Stokes equations and numerical solution for the hyperbolic model of Peshkov and Romenski (HPR) at a final 
time of $t=10.0$ obtained with an ADER-DG $P_3P_3$ scheme ($N=M=3$). Pressure contours and velocity vectors (top left). 41 color contours in the interval [-1,1] of the distortion tensor 
component $A_{11}$ (top right). 1D cuts along the $x$ and the $y$ axis for velocity components $u$ and $v$ (bottom left) and 1D cut along the $x$-axis for the pressure $p$.  } 
\label{fig.tgv}
\end{center}
\end{figure}

The computational results are depicted in Fig. \ref{fig.tgv}, where also a comparison with the exact solution of the incompressible Navier-Stokes equations is shown. Overall, one can note an 
excellent agreement between the HPR model and the reference solution, both for velocity and pressure. The distortion tensor component $A_{11}$ is also drawn in Fig. \ref{fig.tgv} and
reveals the vortex structures of the flow. 

\subsection{3D Taylor-Green vortex}  
\label{sec:tgv3d}

\begin{figure}[!htbp]
\begin{center}
\begin{tabular}{cc} 
\includegraphics[width=0.45\textwidth]{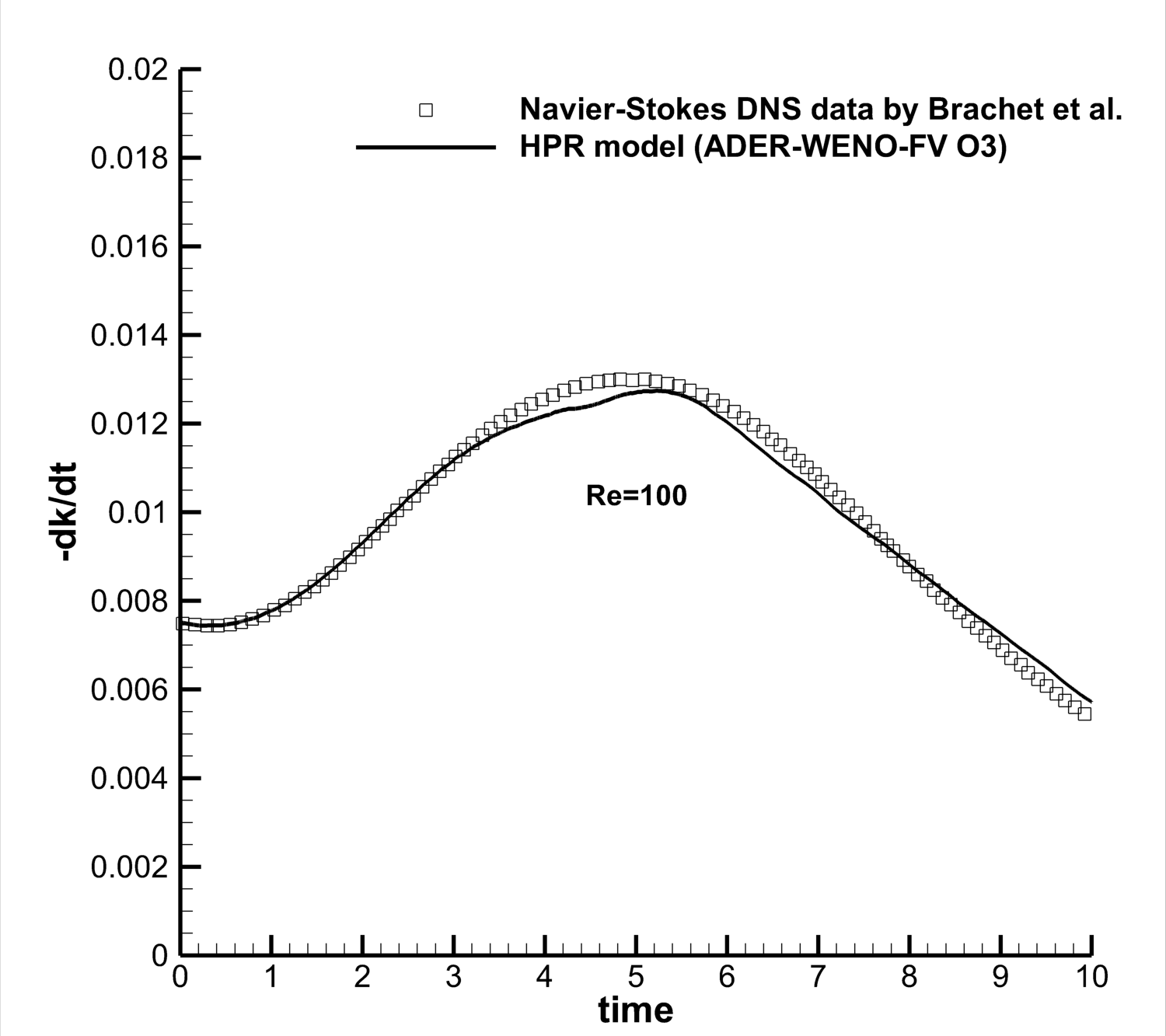}     & 
\includegraphics[width=0.45\textwidth]{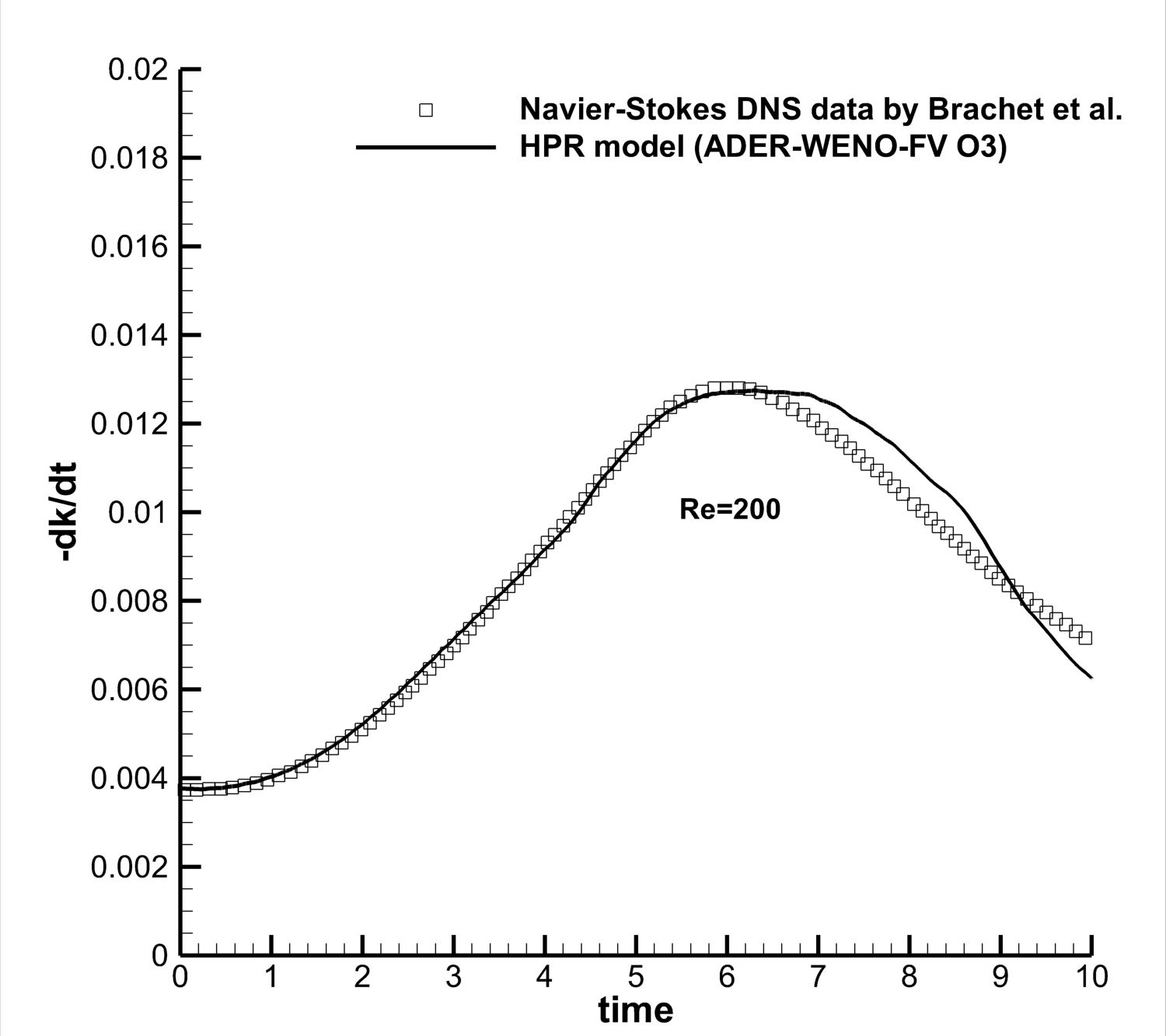}     
\end{tabular}    
\caption{Time evolution of the kinetic energy dissipation rate for the 3D Taylor-Green vortex obtained at a Reynolds number of $\Re=100$ (left) and $\Re=200$ (right). The reference solution is given 
by the DNS data of Brachet et al. \cite{Brachet1983} for the incompressible Navier-Stokes equations, which is compared with the numerical solution for the hyperbolic model of Peshkov and Romenski 
(HPR) until a final time of $t=10$ using a third order $P_0P_2$ ADER-WENO finite volume scheme ($N=0$, $M=2$) on a Cartesian grid of $224^3$ elements. }
\label{fig.dkdt}
\end{center}
\end{figure}

\begin{figure}[!htbp]
\begin{center}
\begin{tabular}{cc} 
\includegraphics[width=0.45\textwidth]{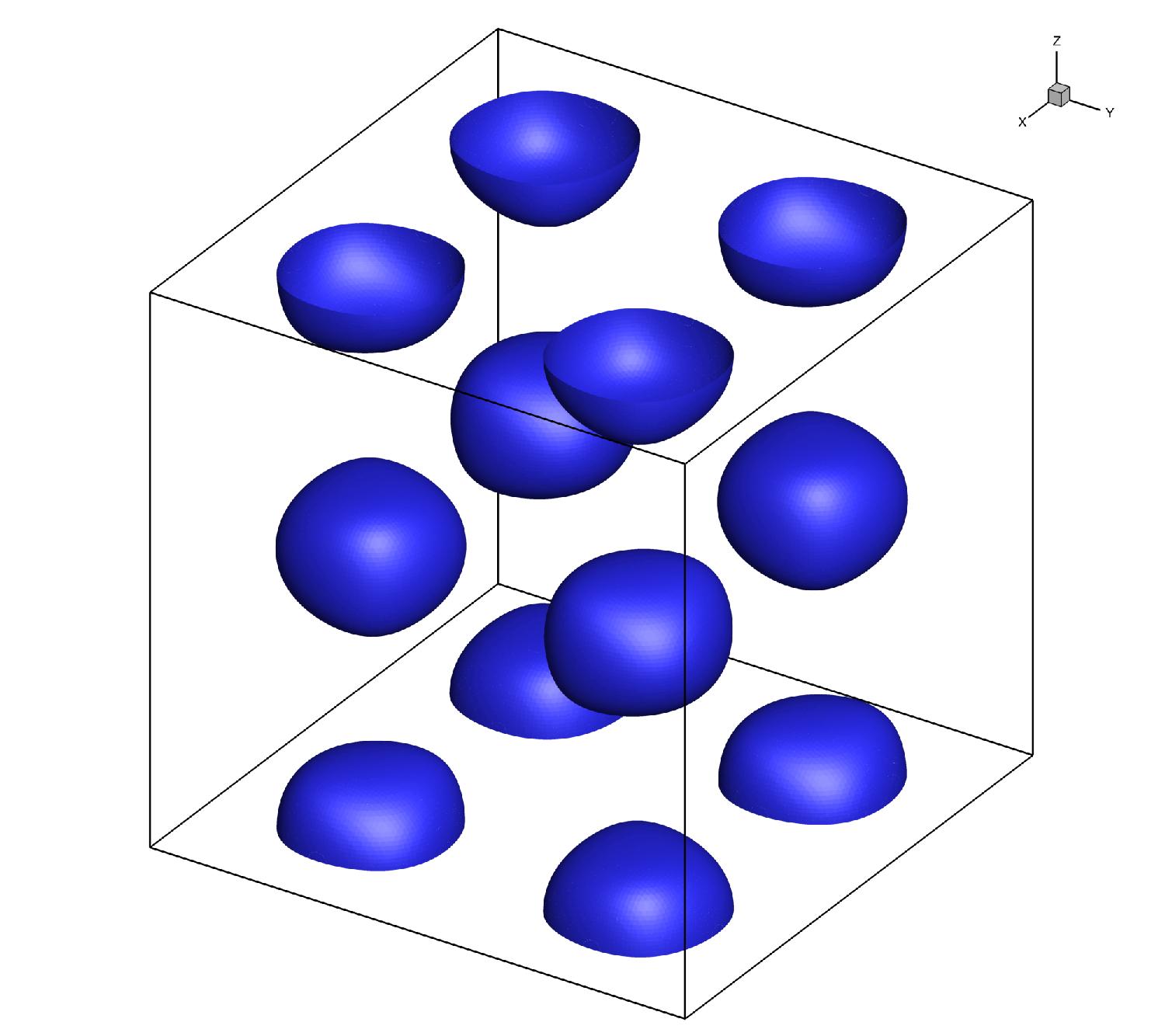}   & 
\includegraphics[width=0.45\textwidth]{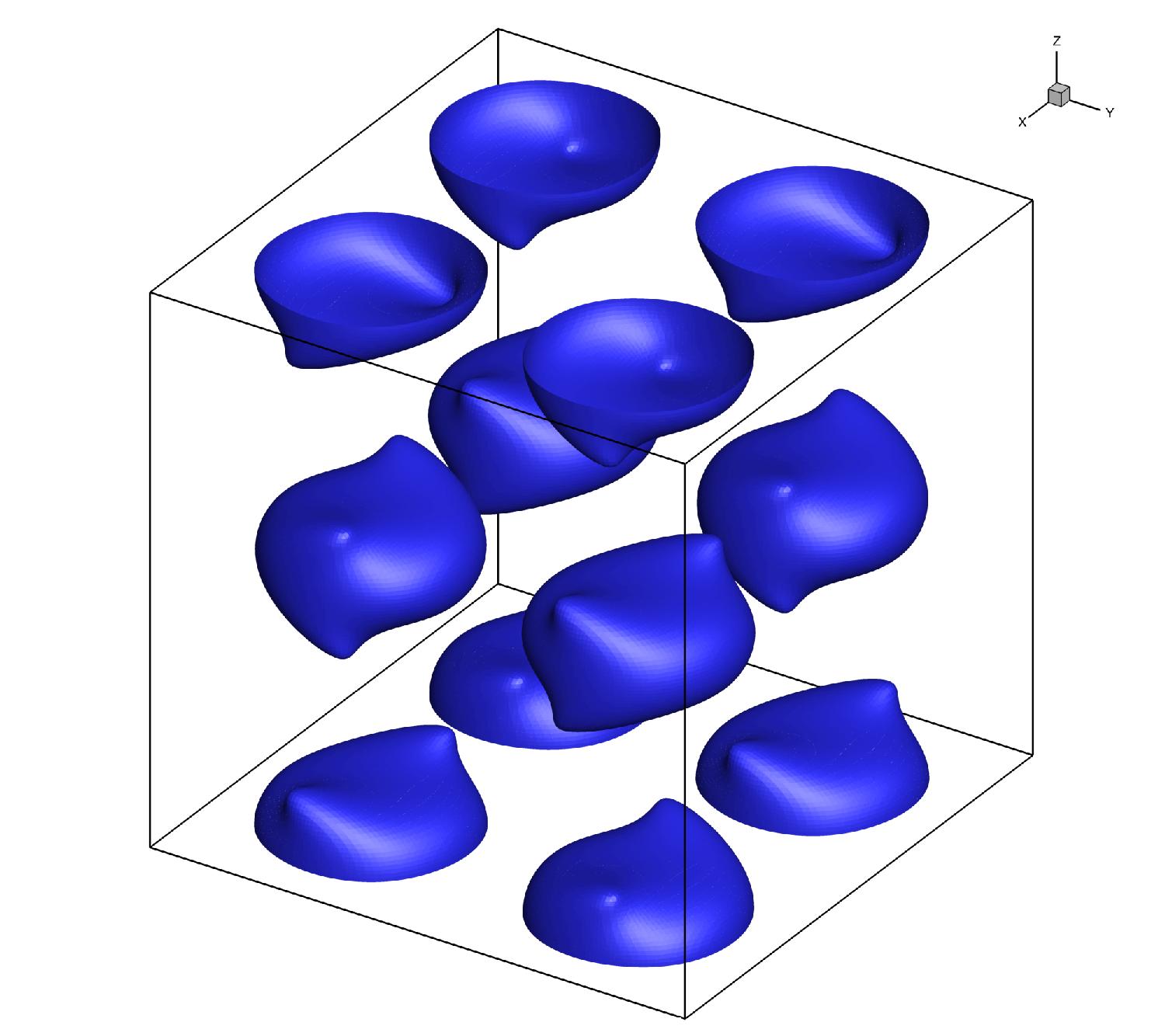}  \\  
\includegraphics[width=0.45\textwidth]{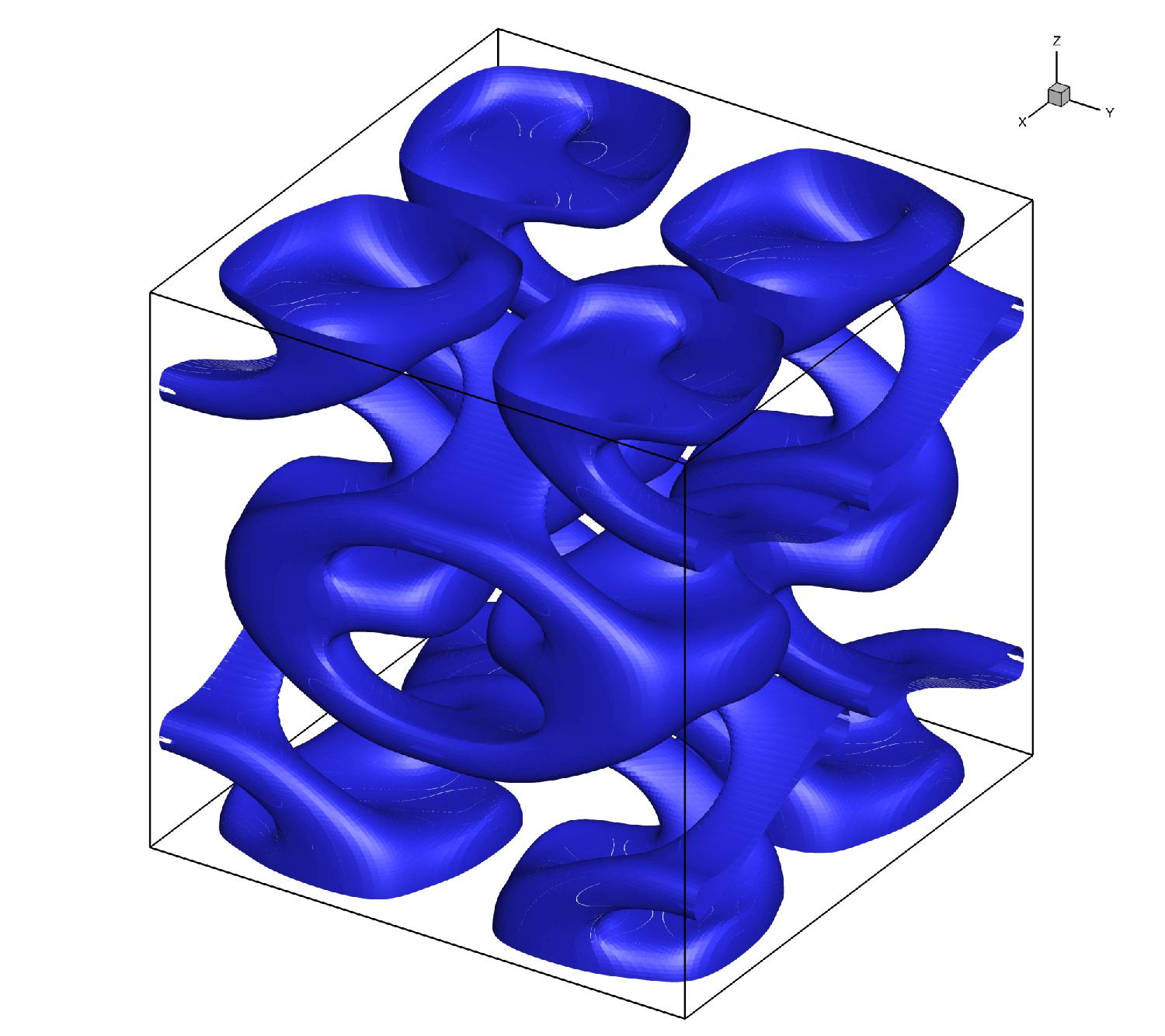}   & 
\includegraphics[width=0.45\textwidth]{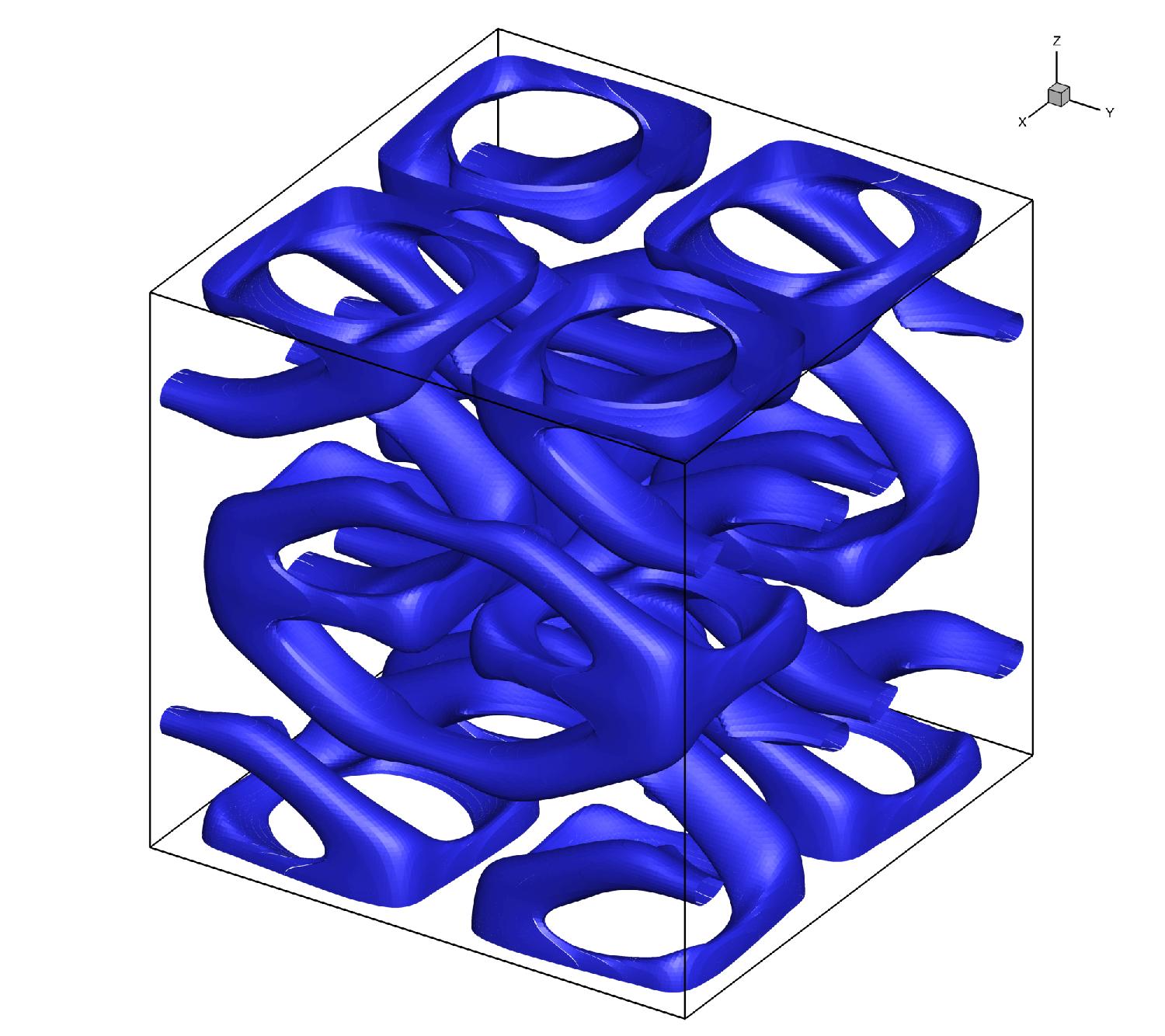}  \\   
\includegraphics[width=0.45\textwidth]{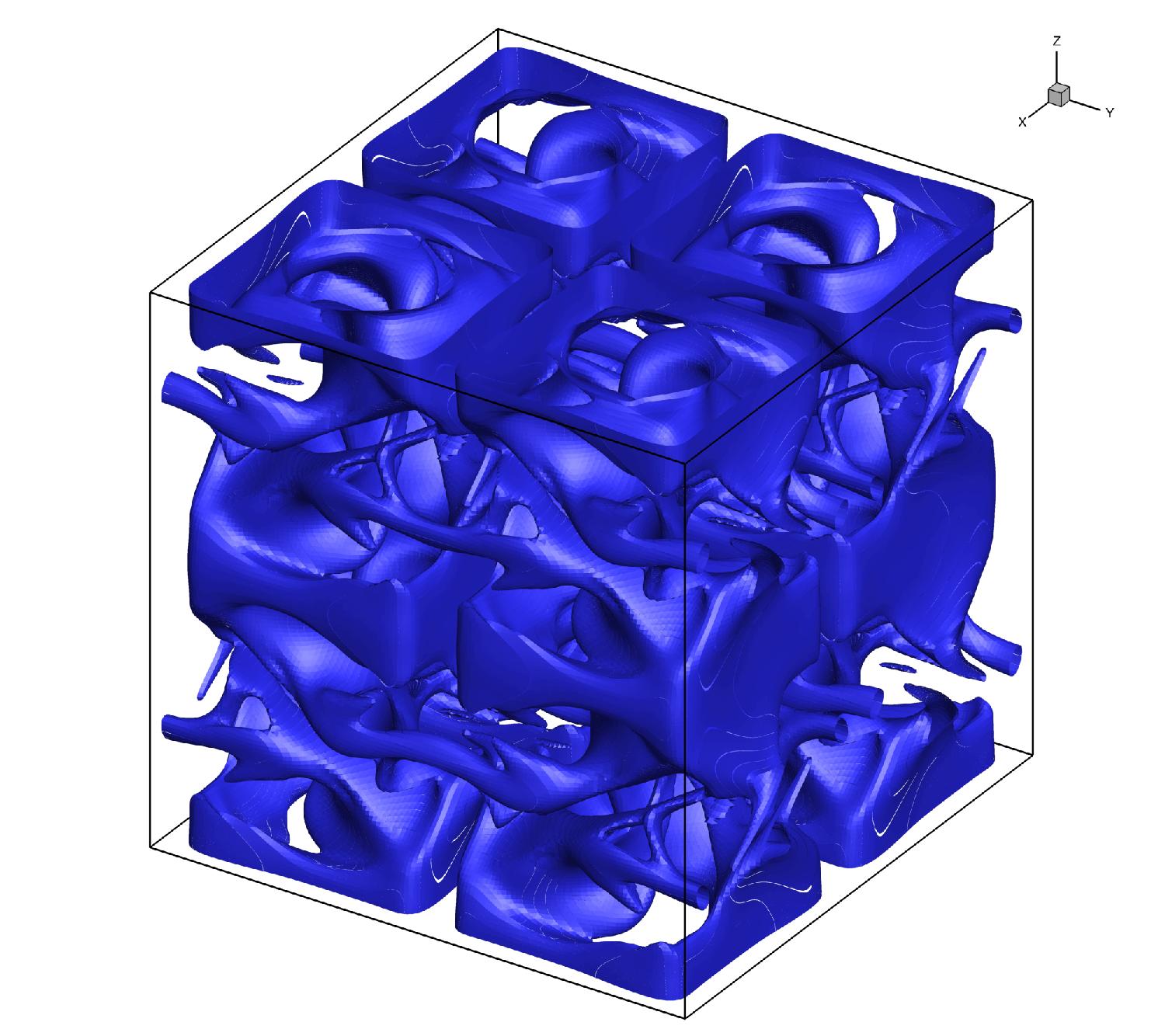}   & 
\includegraphics[width=0.45\textwidth]{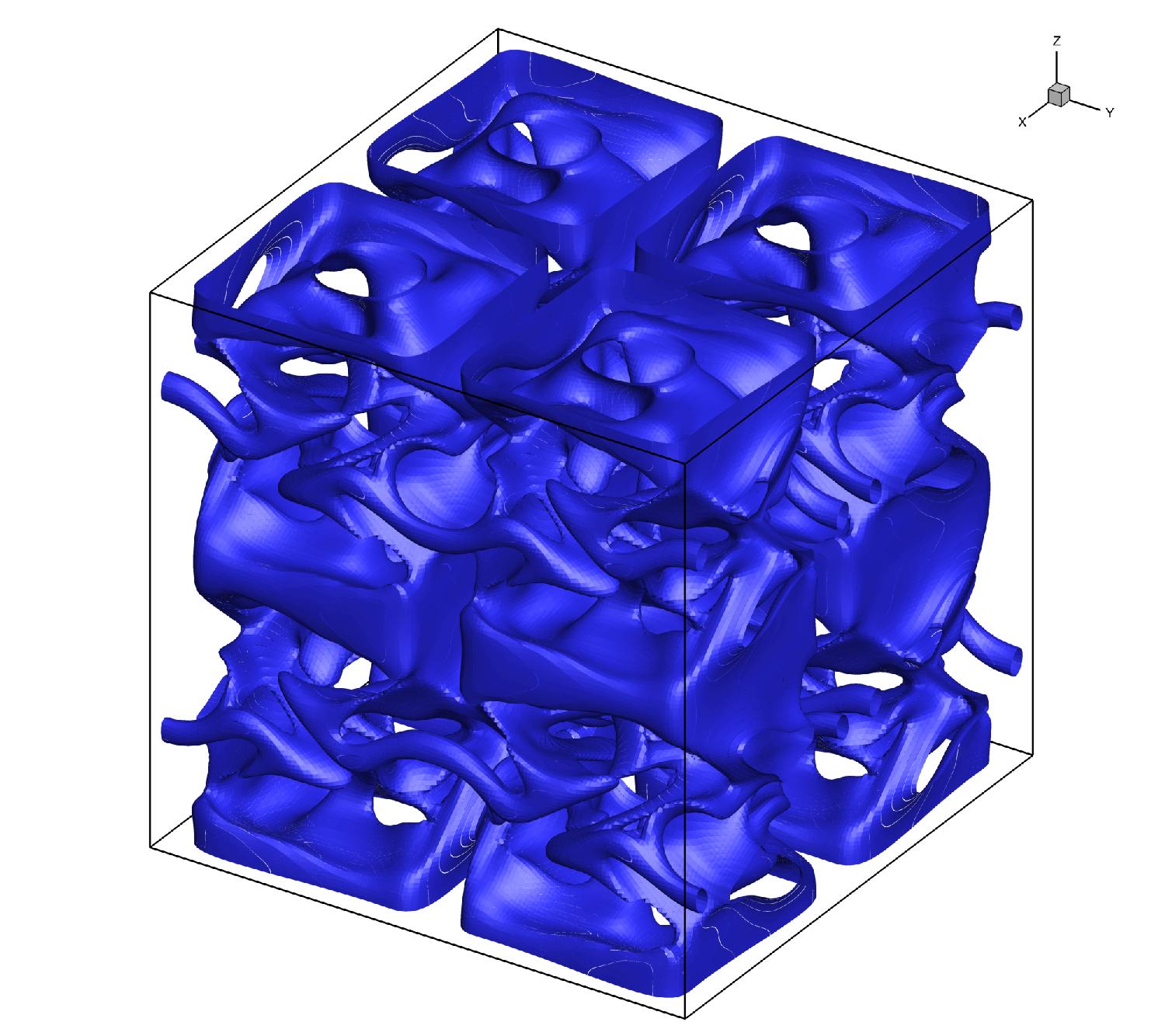}    
\end{tabular} 
\caption{Isocontour surfaces of the tensor component $A_{11}$ obtained with the first order hyperbolic model of Peshkov \& Romenski 
for the 3D Taylor-Green vortex (Re$=200$) for different intermediate times: $t=1$ (top left), $t=2$ (top right), 
$t=3$ (center left) and $t=4$ (center right), 
$t=8$ (bottom left) and $t=10$ (bottom right). } 
\label{fig.tgv3d}
\end{center}
\end{figure}

\begin{figure}[!htbp]
\begin{center}
\begin{tabular}{ccc} 
\includegraphics[width=0.32\textwidth]{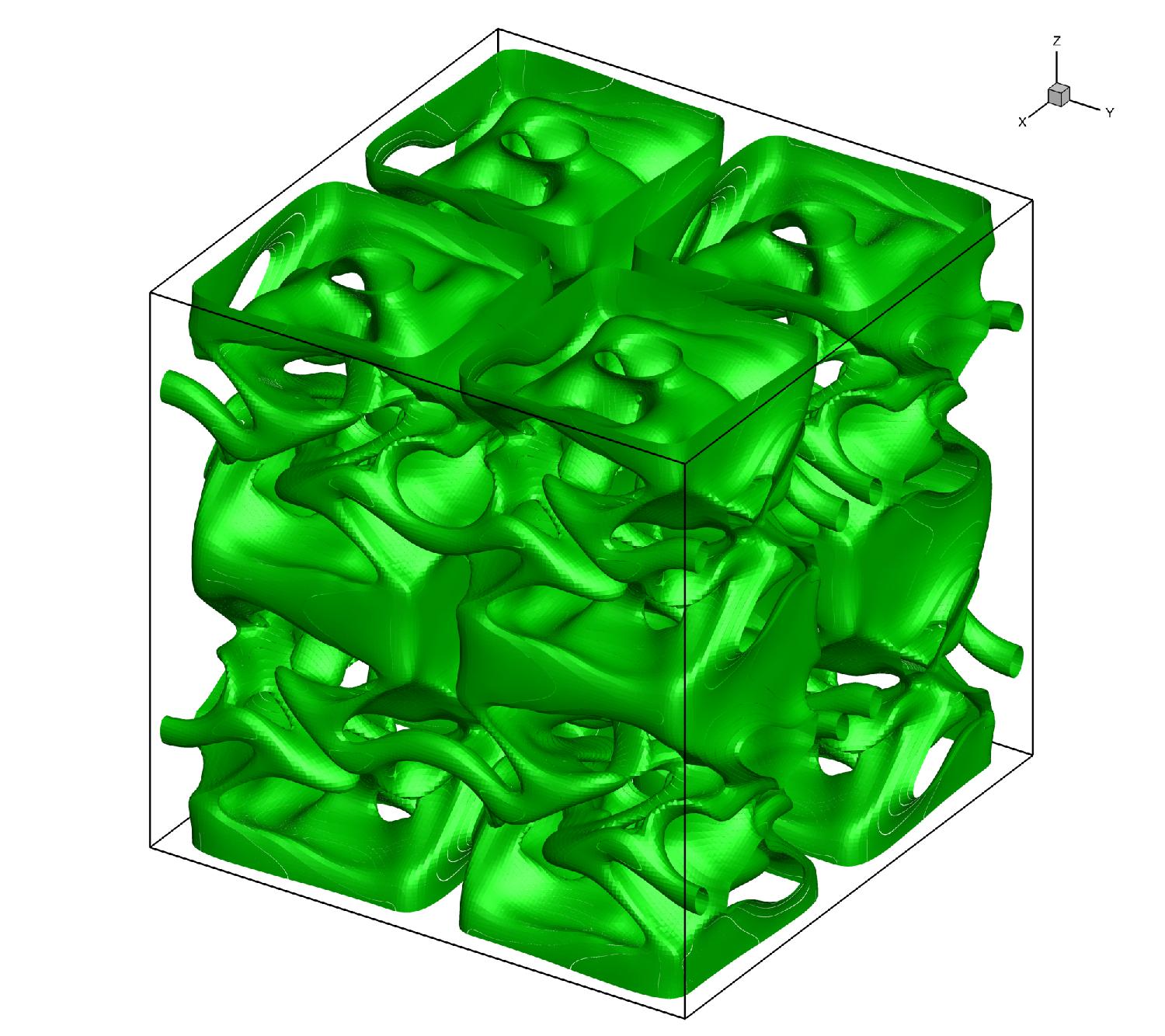}   & 
\includegraphics[width=0.32\textwidth]{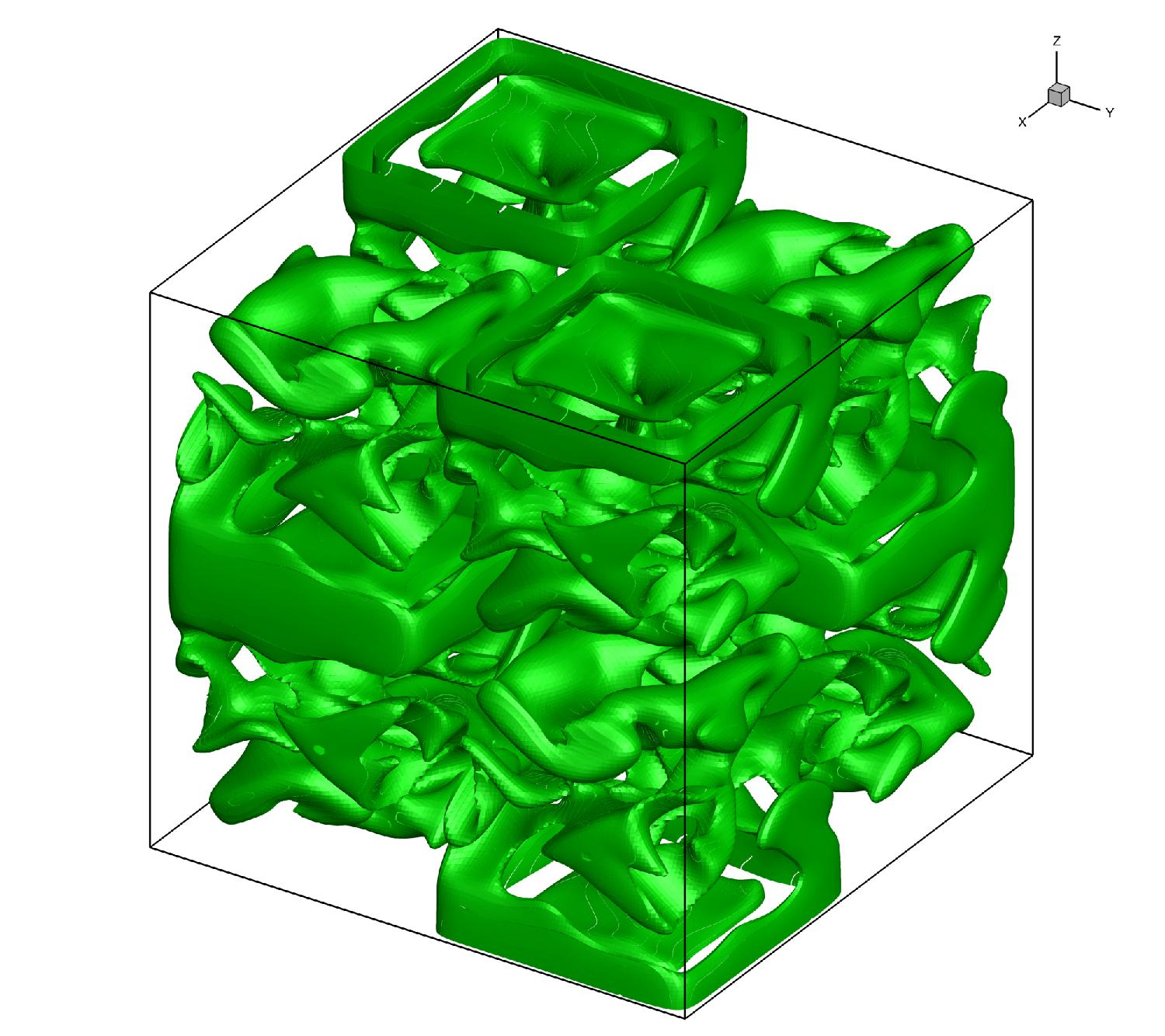}   &   
\includegraphics[width=0.32\textwidth]{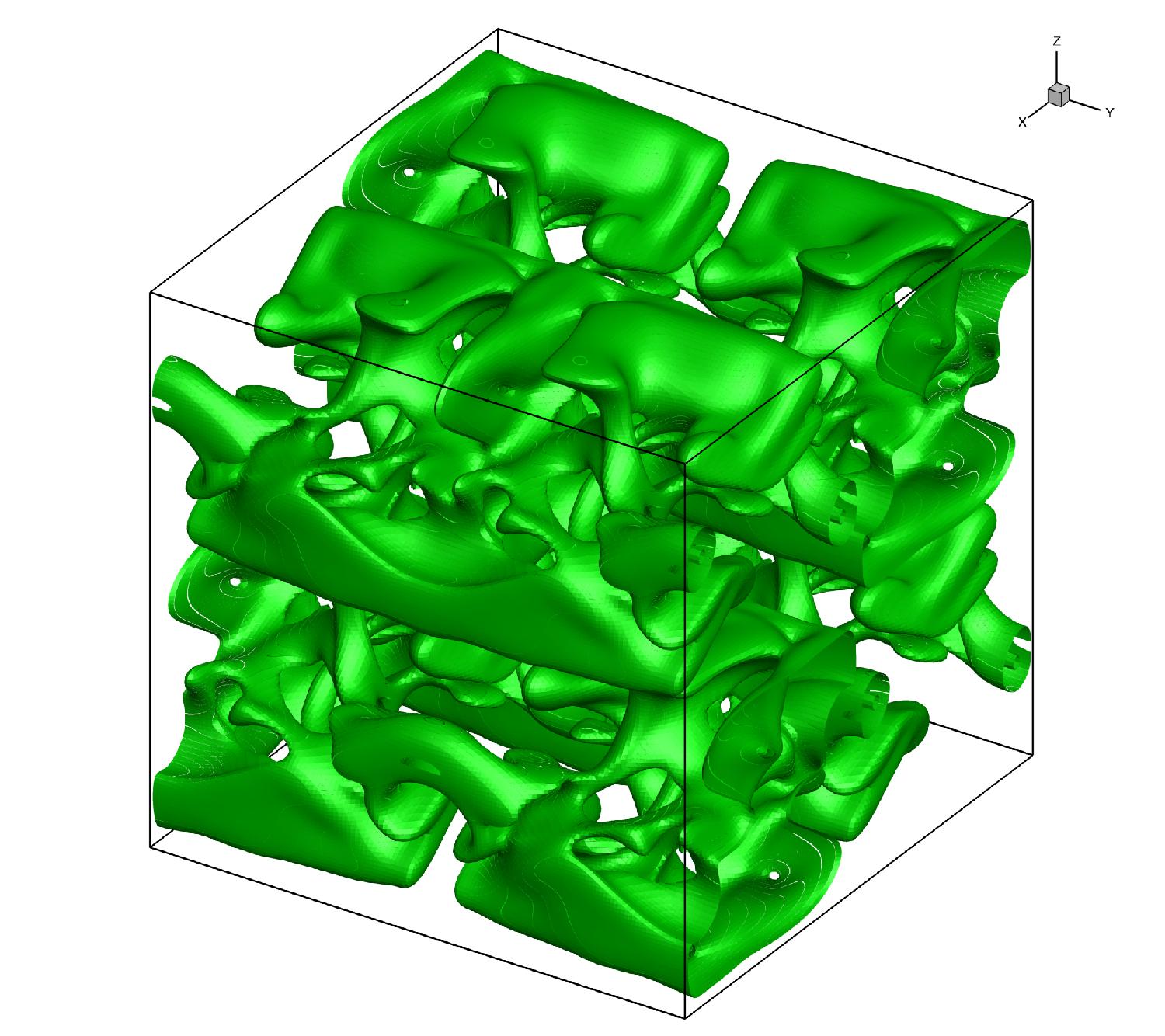}   \\ 
\includegraphics[width=0.32\textwidth]{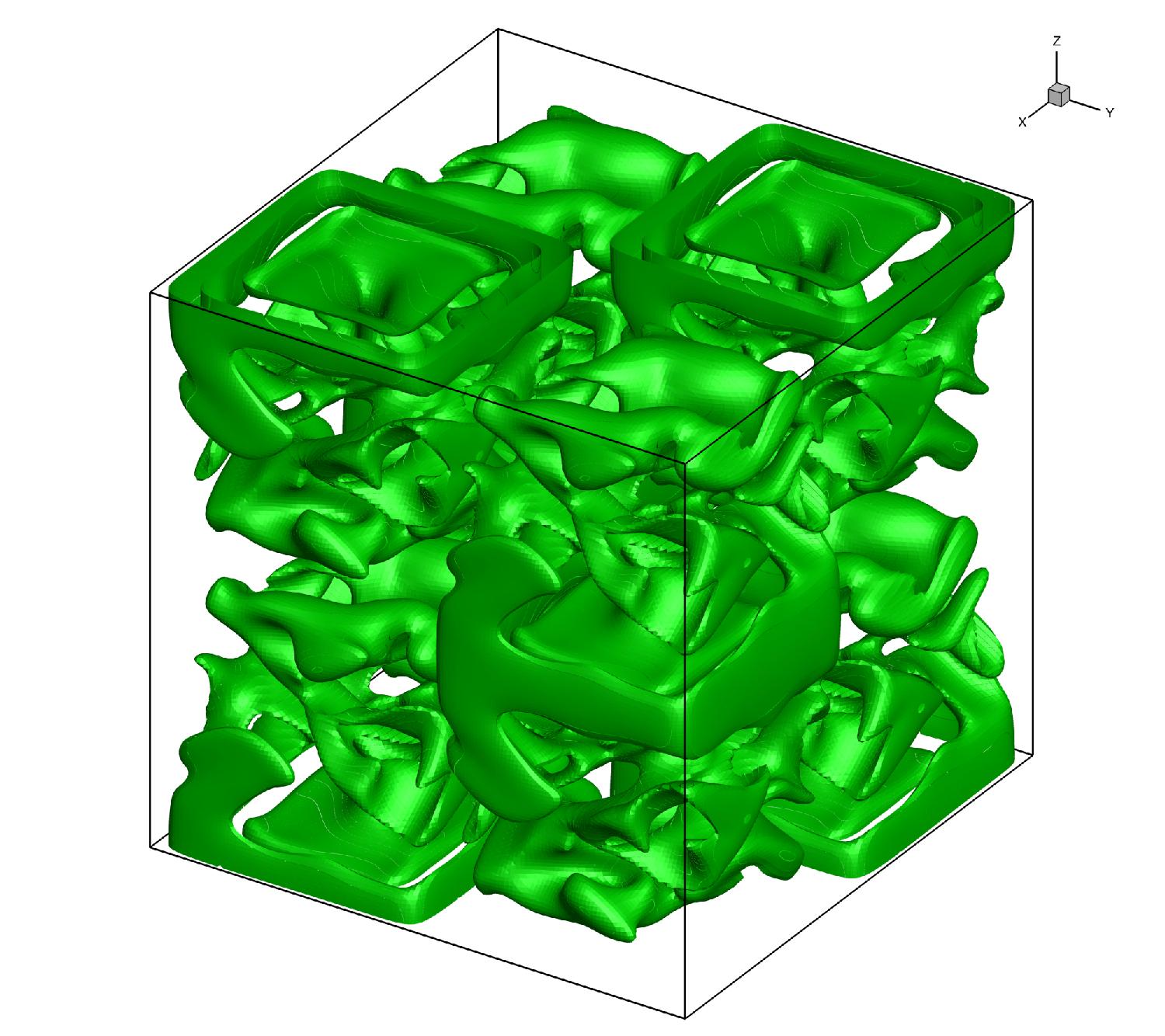}   & 
\includegraphics[width=0.32\textwidth]{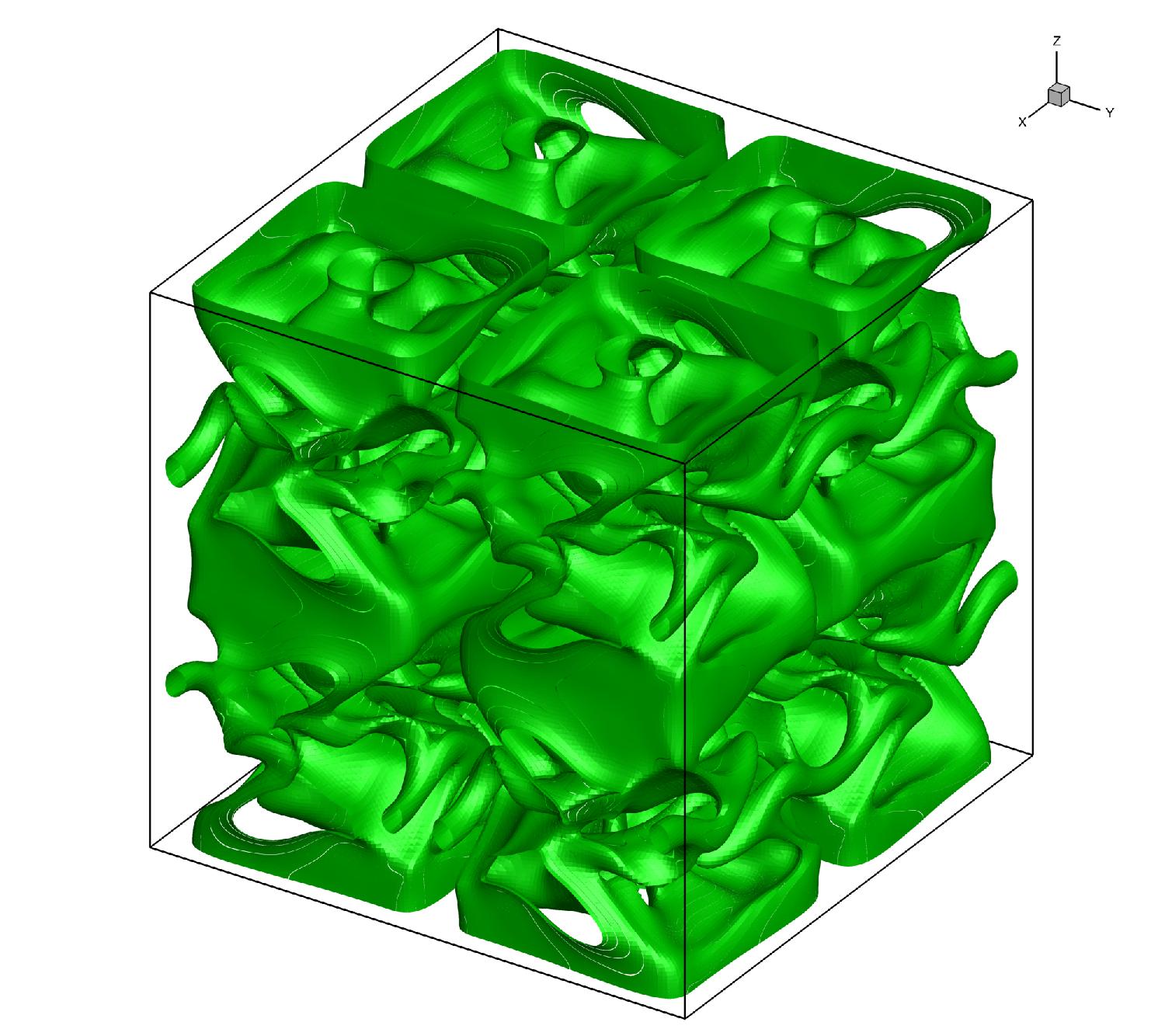}   &   
\includegraphics[width=0.32\textwidth]{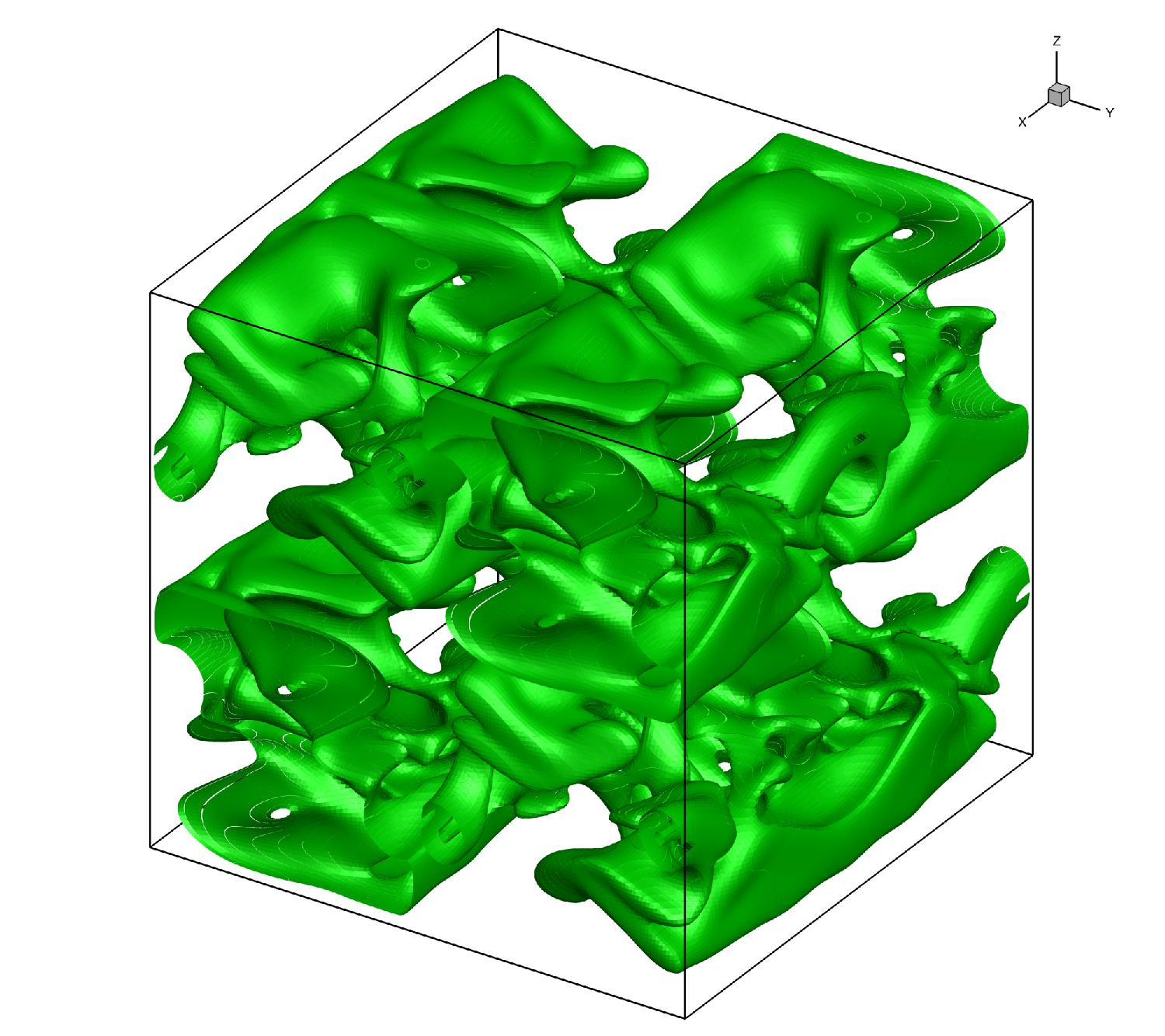}   \\ 
\includegraphics[width=0.32\textwidth]{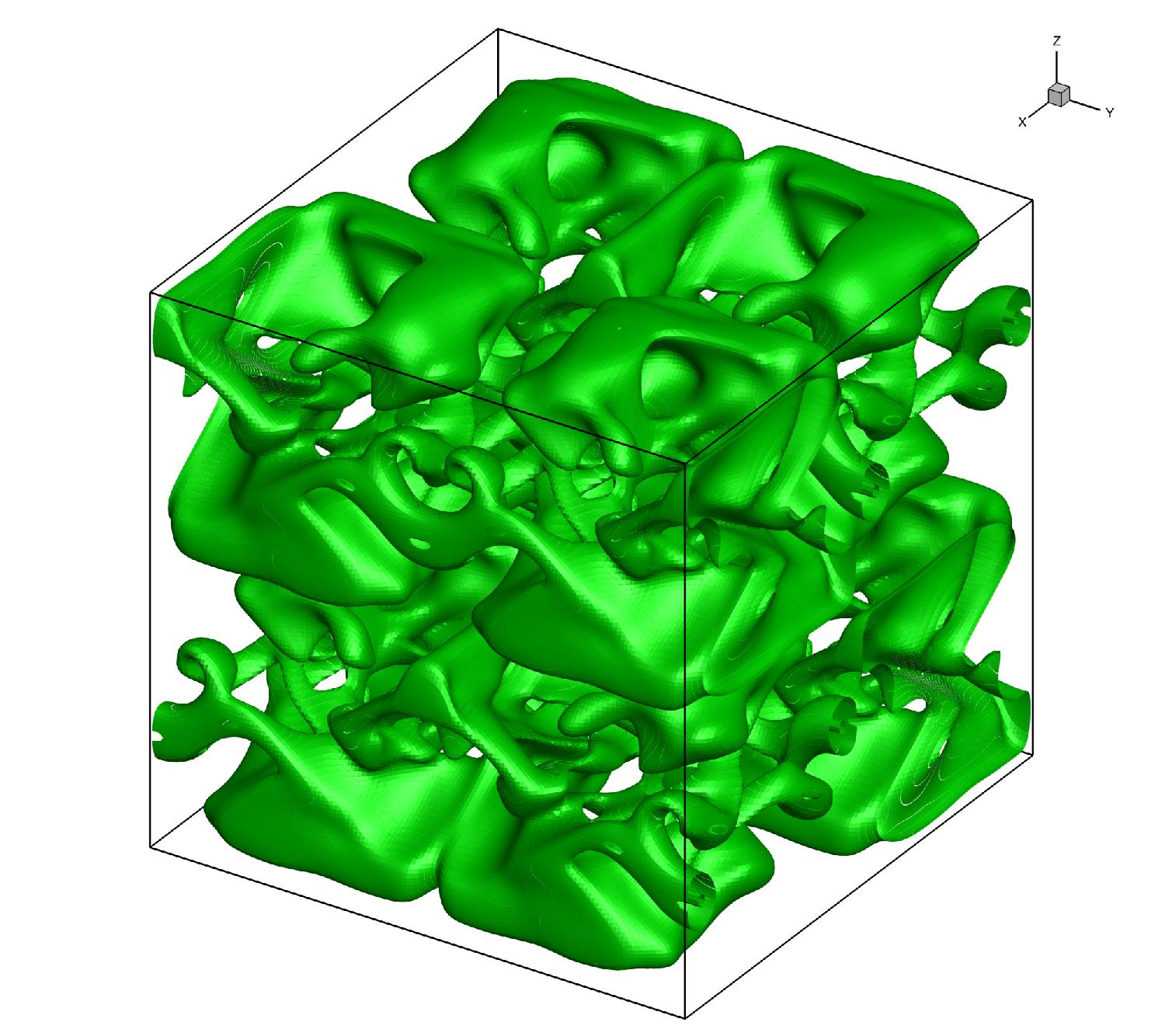}   & 
\includegraphics[width=0.32\textwidth]{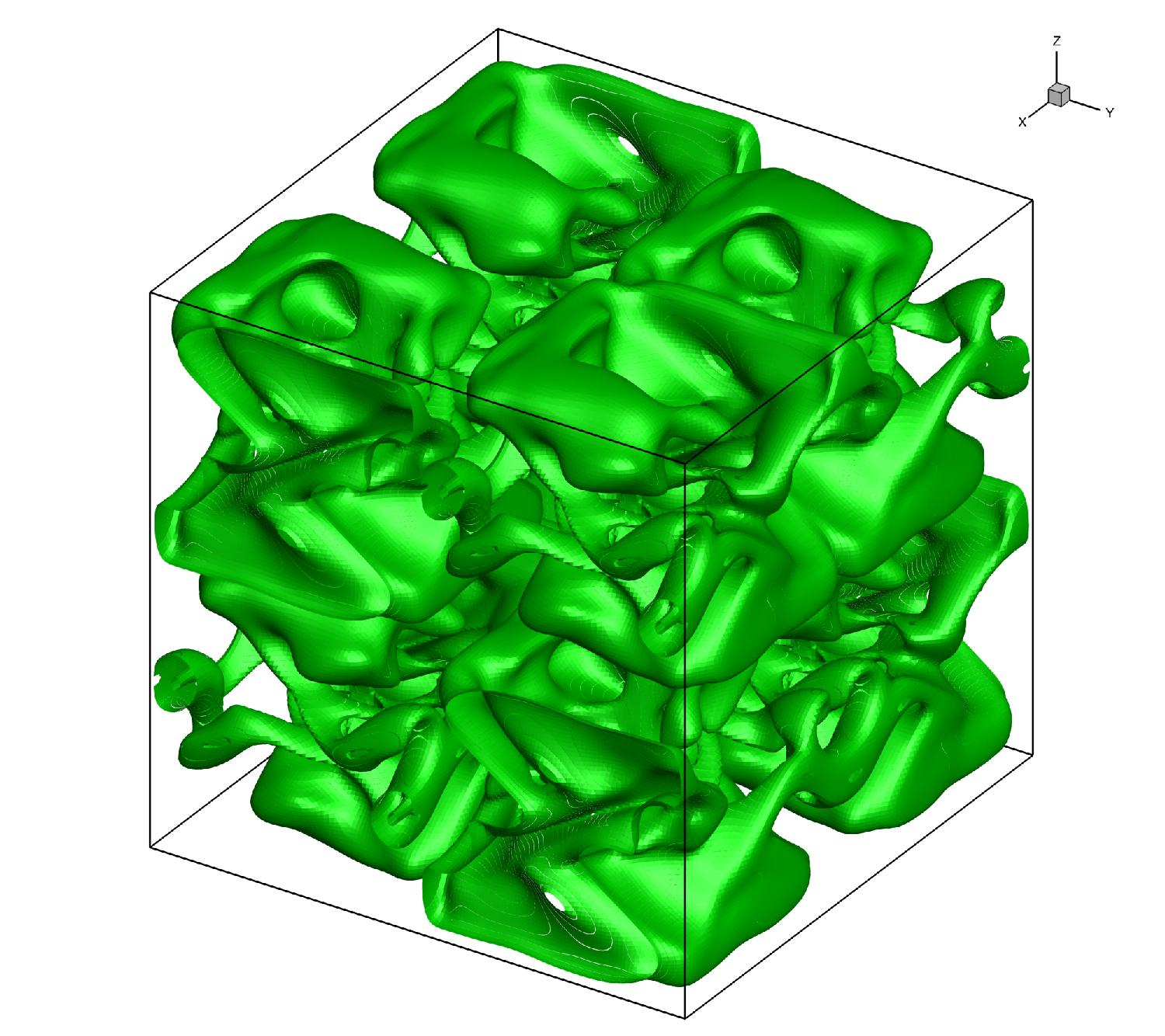}   &   
\includegraphics[width=0.32\textwidth]{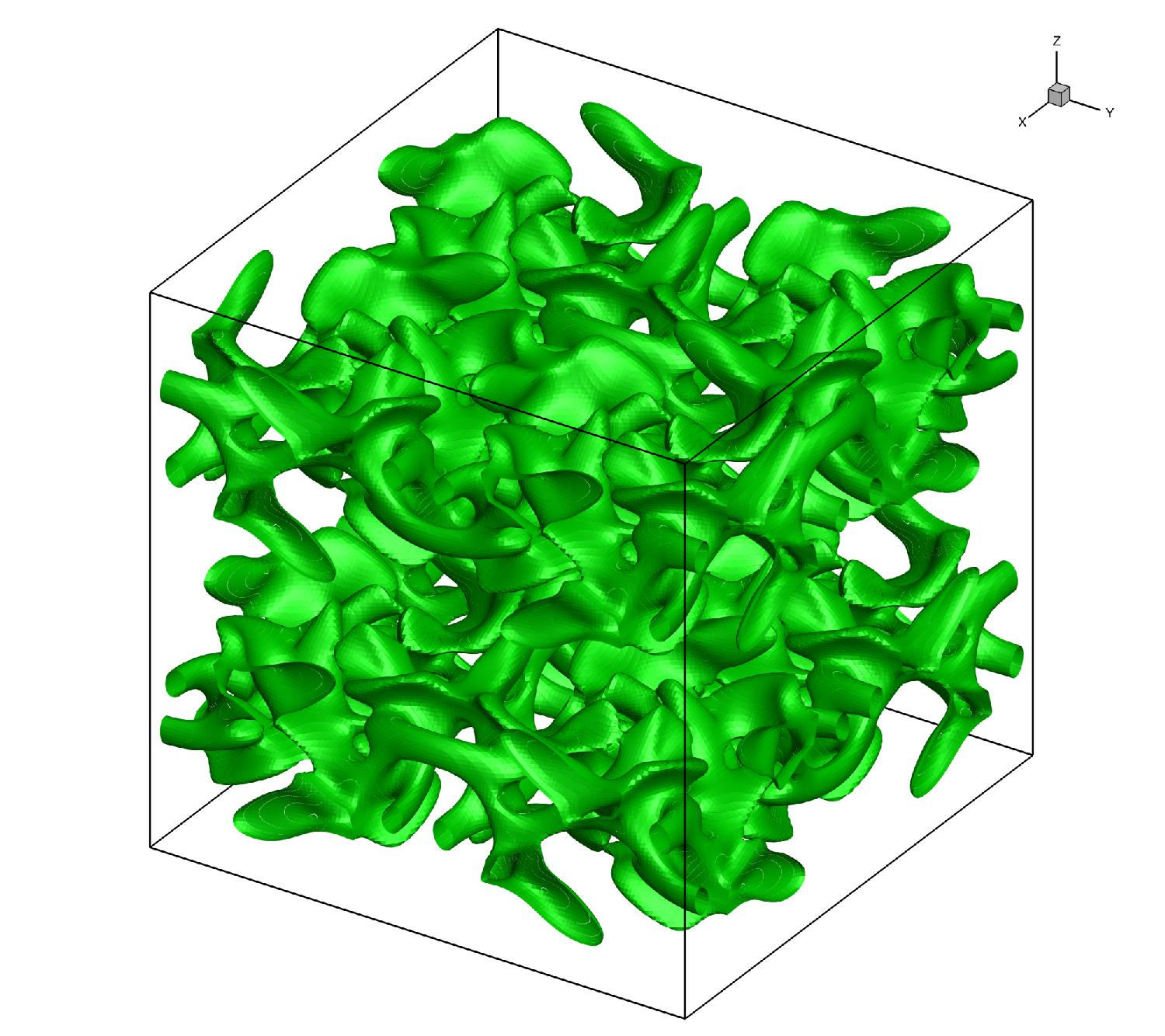}      
\end{tabular} 
\caption{Isocontour surfaces of value $-0.5$ of all components of the tensor $\AAA$ obtained with the first order hyperbolic model of Peshkov \& Romenski 
for the 3D Taylor-Green vortex (Re$=200$) at the final time $t=10$. Component $A_{ik}$ is shown in row $i$ and column $j$ of this figure.  } 
\label{fig.tgv3d.all}
\end{center}
\end{figure}

We now solve the Taylor-Green vortex problem again, but this time in three space dimensions. The computational domain is the box $\Omega = [-\pi,\pi]^3$, with six periodic boundary conditions. 
For large times and large Reynolds numbers, the 3D Taylor-Green vortex is a classical example for the development of flow structures with smaller and smaller spatial scales, up to the 
onset of turbulence. The problem has been widely studied in literature and a reference solution is available via a direct numerical simulation (DNS) provided in the paper of Brachet et al. 
\cite{Brachet1983}. To obtain a low Mach number compressible flow, the following initial condition is chosen, see also \cite{ShuTGV}: $\rho=\rho_0=1$, $p_0=10^2/\gamma$, $\mathbf{A}=\mathbf{I}$, 
$\mathbf{J}=0$, and  
\begin{eqnarray}
      u(\mathbf{x},0) &=&   \sin(x) \cos(y) \cos(z),  \nonumber \\  
      v(\mathbf{x},0) &=& - \cos(x) \sin(y) \cos(z),  \nonumber \\ 
      z(\mathbf{x},0) &=& 0, \nonumber \\
      p(\mathbf{x},0) &=& p_0 + \frac{\rho_0}{16} \left( (\cos(2z) + 2) (\cos(2x) + \cos(2 y)) - 2.0 \right).  
\label{eq:TGV3D}
\end{eqnarray}
With the choice of $p_0$, the maximum Mach number of the flow at the initial time is M$=0.1$. 
The numerical simulations are carried out for two different Reynolds numbers with the full HPR model in three space dimensions up to a final time of $t=10$, using a third order 
ADER-WENO finite volume scheme ($N=0$, $M=2$), on a computational grid composed of $224^3$ elements. For the HPR model, the following set of parameters has been chosen: 
$\gamma=1.4$, $\rho_0=1$, $\mu=1/\textnormal{Re}$, $c_v=2.5$, $c_s=8$, $\alpha=1$, $T_0=1$, $\kappa=\gamma c_v \mu$, hence the resulting Prandtl number is Pr=$1$. An important quantity 
for the comparison with existing DNS data is the kinetic energy dissipation rate 
\begin{equation}
  \frac{d k}{d t} = \frac{d}{d t} \left( \frac{1}{\left| \Omega \right|} \int \limits_\Omega \frac{1}{2} \rho \mathbf{v}^2 d\x \right). 
	\label{eqn.dkdt} 
\end{equation} 
The computational results obtained for the kinetic energy dissipation rate are depicted in Fig. \ref{fig.dkdt} for two Reynolds numbers, $\Re=100$ and $\Re=200$, where we can note 
an good agreement with the DNS reference data of Brachet et al. for Re$=100$, while the employed third order ADER-WENO finite volume scheme seems to be too dissipative for 
$t>6$ in the case $\Re=200$. Further systematic studies of this important test problem will be carried out in the future, using substantially refined grids and higher order 
polynomial degrees in the numerical scheme in order to identify the cause for the observed deviation of the HPR model from the Navier-Stokes reference solution for higher 
Reynolds numbers. 
A 3D view of the time evolution of the developing small-scale flow structures is shown in Fig. \ref{fig.tgv3d} at the aid of the component $A_{11}$ of the distortion tensor 
$\AAA$, while all elements of $\AAA$ are depicted at the final time $t=10$ in Fig. \ref{fig.tgv3d.all}. 

\subsection{Heat conduction in a gas} 
\label{sec:heat}

\begin{figure}[!htbp]
\begin{center}
\begin{tabular}{cc} 
\includegraphics[width=0.4\textwidth]{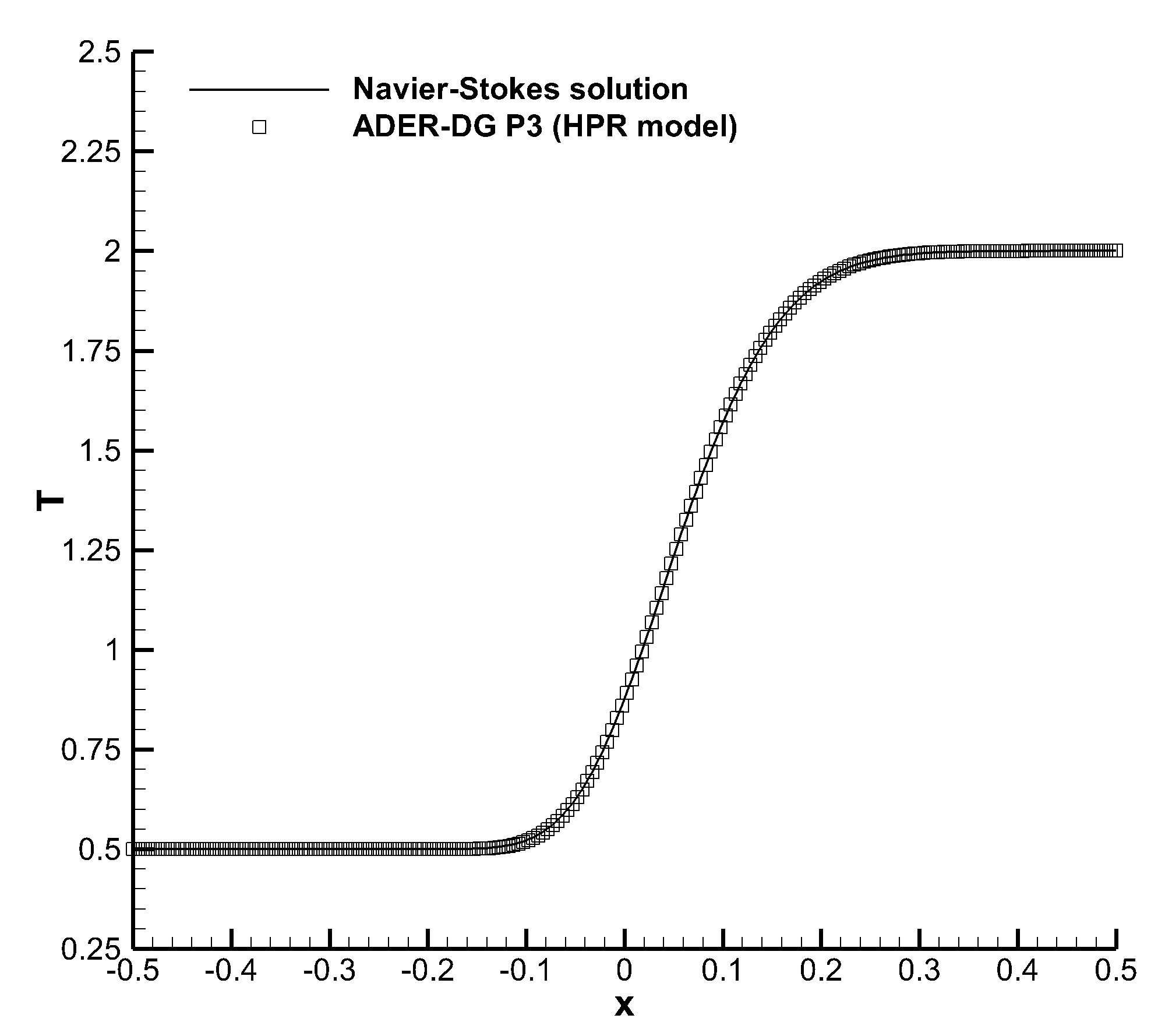}  &  
\includegraphics[width=0.4\textwidth]{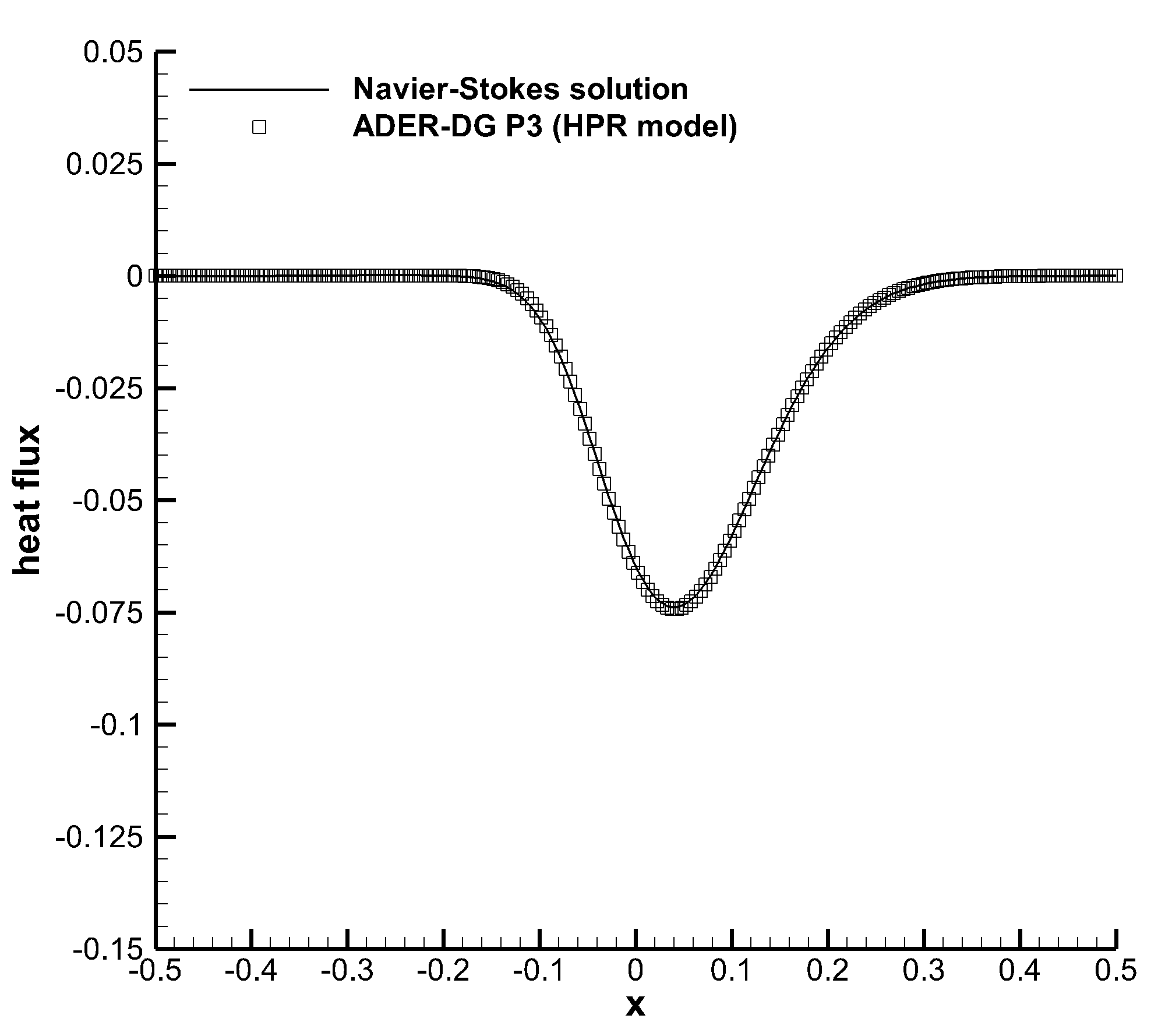}   
\end{tabular} 
\caption{Heat conduction test problem at a final time of $t=1.0$. Temperature distribution (left) and heat flux (right). For the Navier-Stokes solution, 
the classical Fourier heat flux $q_1 = -\kappa T_x$ is shown, while for the first order HPR model, we plot $q_1 = \alpha^2 J_1 T$. }  
\label{fig.heat}
\end{center}
\end{figure}

Here, we solve a simple test problem dominated by the effect of heat conduction. The initial condition for the flow variables is $\rho=2$ for $x<0$ and
$\rho=0.5$ for $x\geq 0$, while $u=v=0$ and $p=1$ everywhere. Furthermore, $\mathbf{A}=\mathbf{I}$ and $\mathbf{J}=0$ at the initial
time. The parameters of the HPR model are defined by $\gamma=1.4$, $\rho_0=1$, $c_v=2.5$, $c_s=1$, $\mu=10^{-2}$, $\alpha=2$, $T_0=1$ and $\kappa = 10^{-2}$. 
The computational domain is $\Omega = [-0.5, 0.5] \times [-0.1, 0.1]$ and simulations are carried out with an ADER-DG $P_3P_3$ scheme ($N=M=3$) 
until $t=1.0$ on a grid composed of $100 \times 5$ elements. In Fig. \ref{fig.heat} a 1D cut through the computational results at $y=0$ is shown for 
the first order HPR model, together with a Navier-Stokes reference solution computed with an ADER-DG scheme \cite{ADERNSE} on the same grid. 
We note an excellent agreement between the two models. 
Furthermore, in Table \ref{tab.heat.cpu} we list the computational time (wallclock time) and the number of time steps needed to reach the final simulation 
time in the context of ADER finite volume and ADER-DG schemes using two different grid resolutions for both, the HPR model and for the 
compressible Navier-Stokes equations. With the above parameters, the explicit discretization of the compressible Navier-Stokes equations already runs into
the parabolic time step restriction, which leads to a quadratic decrease of the time step size with mesh refinement and thus to a significant increase
in computational time. Compared to the finite volume case, the situation is even worse for DG methods, where the discretization of the HPR model does not 
only require much less time steps, but also less CPU time than the discretization of the compressible Navier-Stokes equations. We think that these results 
might be relevant for those readers who are interested in the discretization of viscous compressible flows with heat conduction using explicit time 
integration schemes. 

\begin{table}[!t]
 \caption{CPU time comparison for the heat conduction problem using the Navier-Stokes equations and the first order HPR model.}
\begin{center}
 \begin{tabular}{ccccc}
 \hline
 & \multicolumn{2}{c}{\textbf{Navier-Stokes equations}} &  \multicolumn{2}{c}{\textbf{HPR model}}\\
 \hline
 Mesh & time steps & CPU time & time steps & CPU time \\ 
 \hline
 \multicolumn{5}{c}{ADER-WENO finite volume scheme ($P_0P_2$)} \\
 \hline
 100 & 1587 & 18.7  & 461 & 42.7  \\ 
 200 & 5535 & 112.2 & 922 & 144.8 \\ 
 \hline
 \multicolumn{5}{c}{ADER-DG scheme ($P_3P_3$)} \\
 \hline
 100 & 87080  & 2317.2  & 4554 & 651  \\ 
 200 & 340646 & 18476   & 9104 & 2437 \\ 
 \hline
 \end{tabular}
\end{center}
\label{tab.heat.cpu}
\end{table}

\subsection{Viscous shock profile} 
\label{sec:shock}

The numerical test problems solved so far considered only either low Mach number flows, or at most weakly compressible flows. 
However, the first order HPR model is also valid in the case of supersonic viscous flows. Therefore, in this section we 
solve the problem of an isolated viscous shock wave propagating into a medium at rest with a shock Mach number of $M_s>1$. 
In the case of a Prandtl number of Pr$=0.75$, there exists an exact traveling wave solution of the compressible Navier 
-Stokes equation that was first found by Becker \cite{Becker1923} in 1923. In the following, we briefly recall the exact
solution of Becker, where the indices "0" and "1" denote the upstream and the post-shock states, respectively. 

\begin{figure}[!htbp]
\begin{center}
\begin{tabular}{ccc} 
\includegraphics[width=0.3\textwidth]{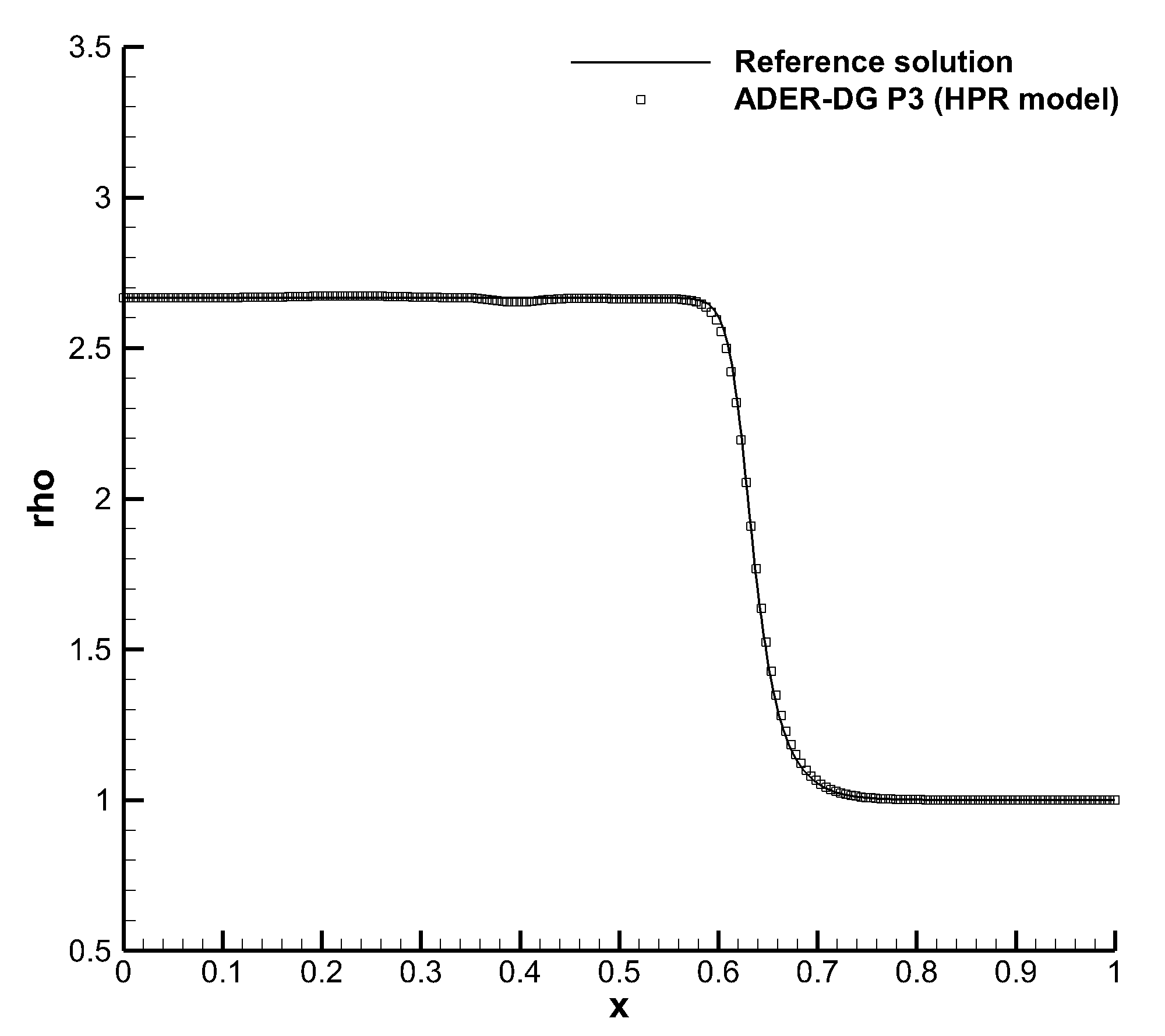}  &  
\includegraphics[width=0.3\textwidth]{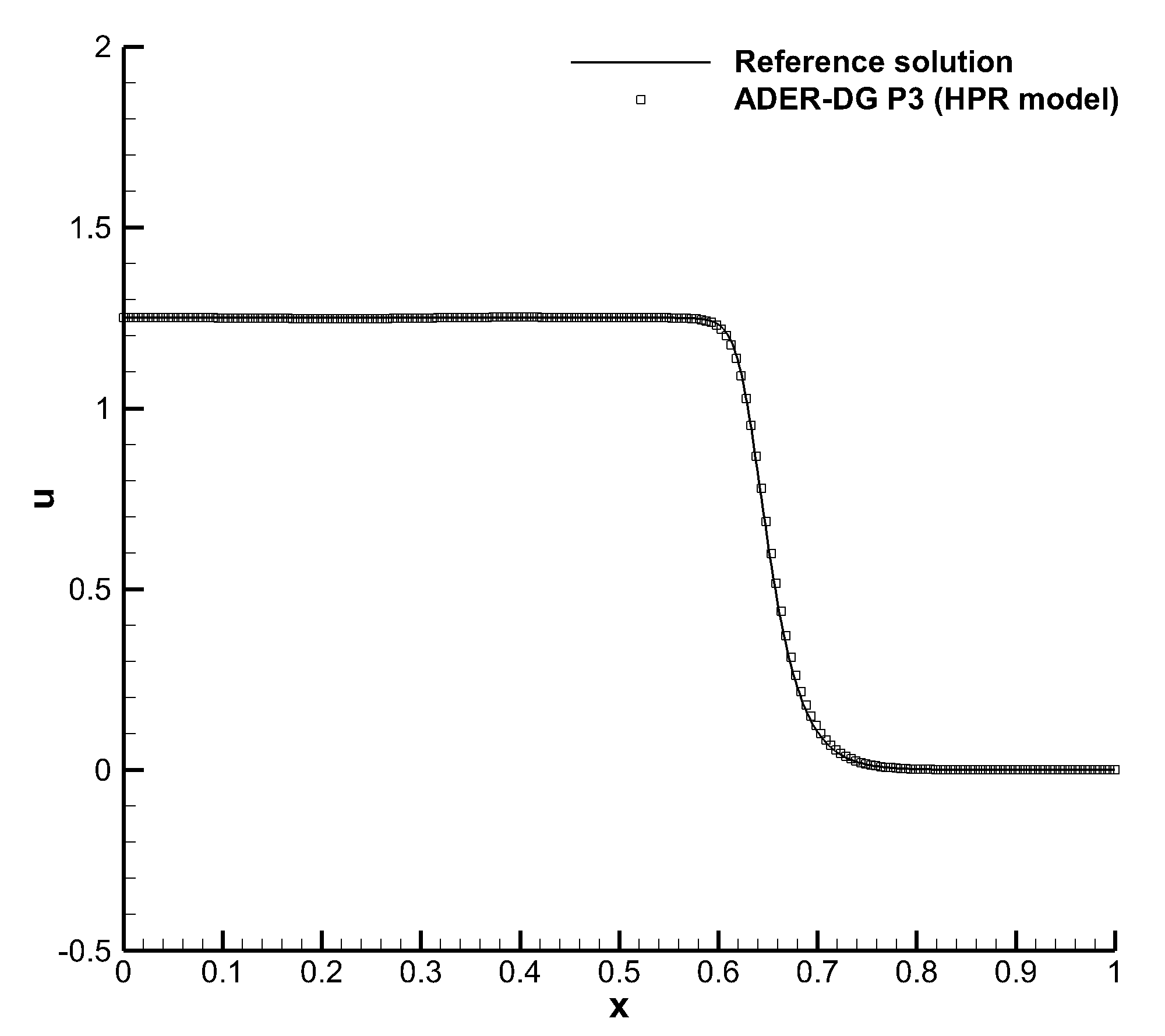}    &  
\includegraphics[width=0.3\textwidth]{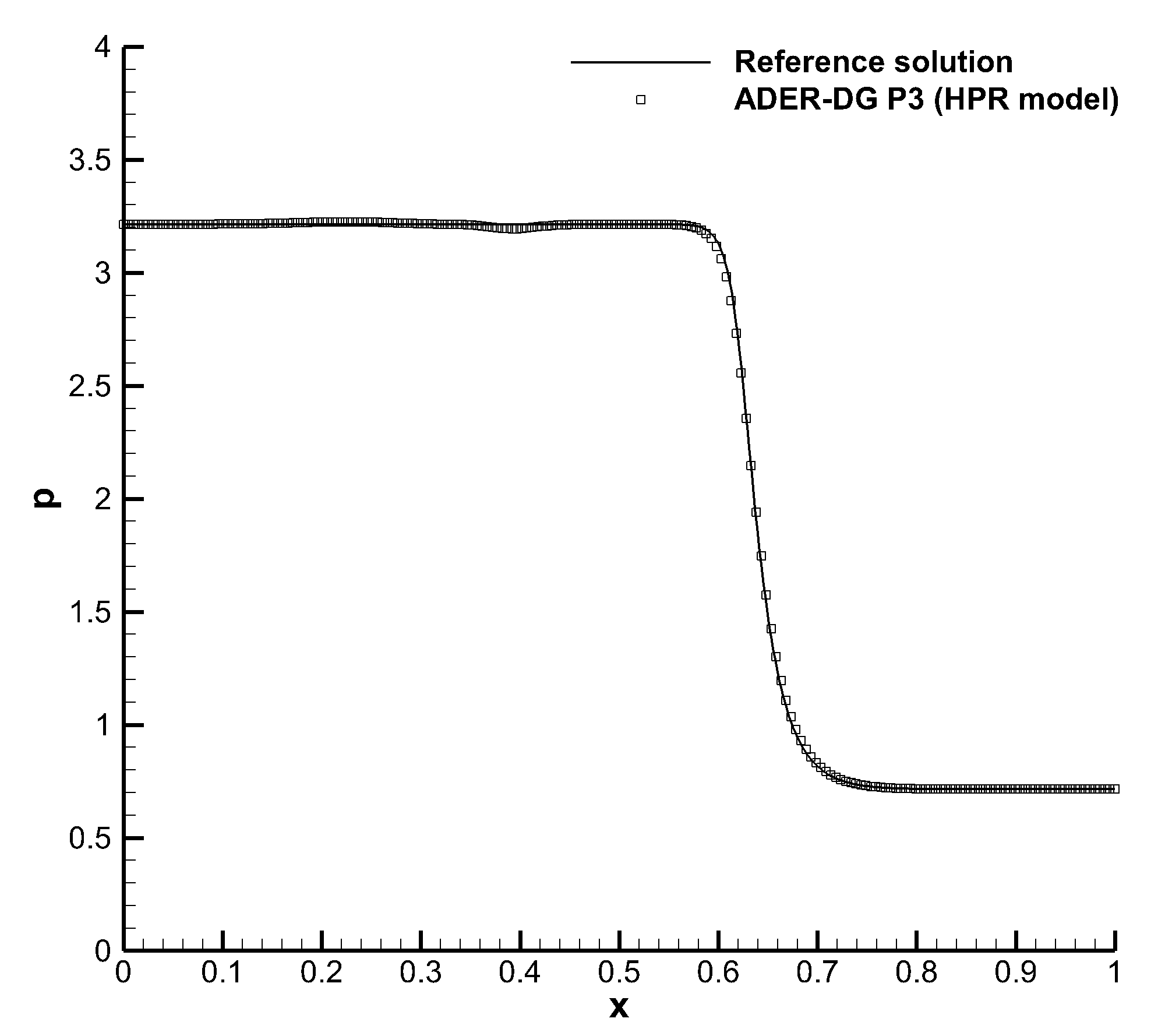}     
\end{tabular} 
\caption{Viscous shock with shock Mach number $M_s=2$ and Prandtl number $Pr=0.75$ at a final time of $t=0.2$. Comparison of the exact solution of the compressible Navier-Stokes
equations according to Becker \cite{Becker1923} with the HPR model. Density profile (left) velocity profile (middle) and pressure profile (right). }  
\label{fig.vshockprim}
\end{center}
\end{figure}

\begin{figure}[!htbp]
\begin{center}
\begin{tabular}{cc} 
\includegraphics[width=0.45\textwidth]{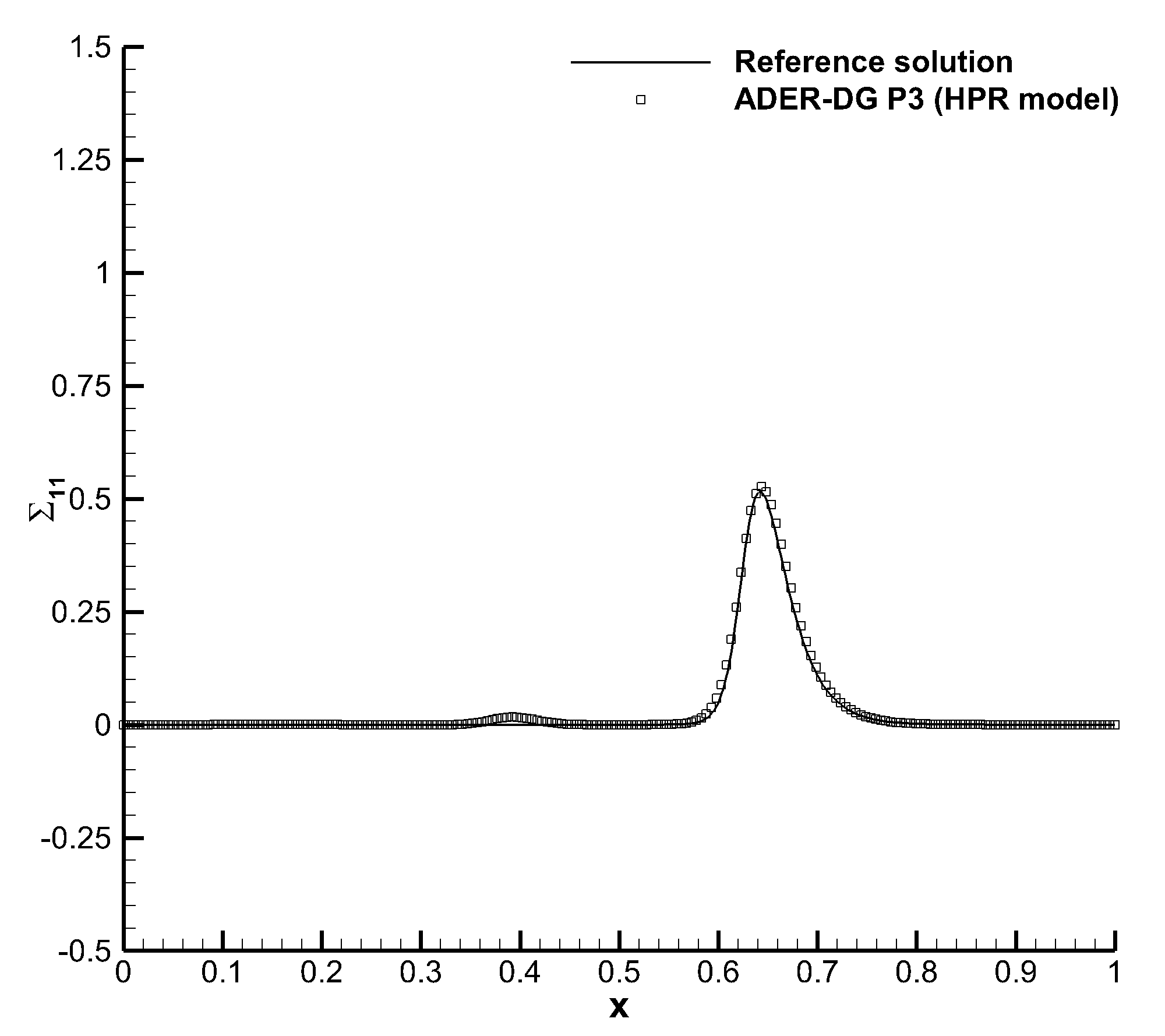}  &  
\includegraphics[width=0.45\textwidth]{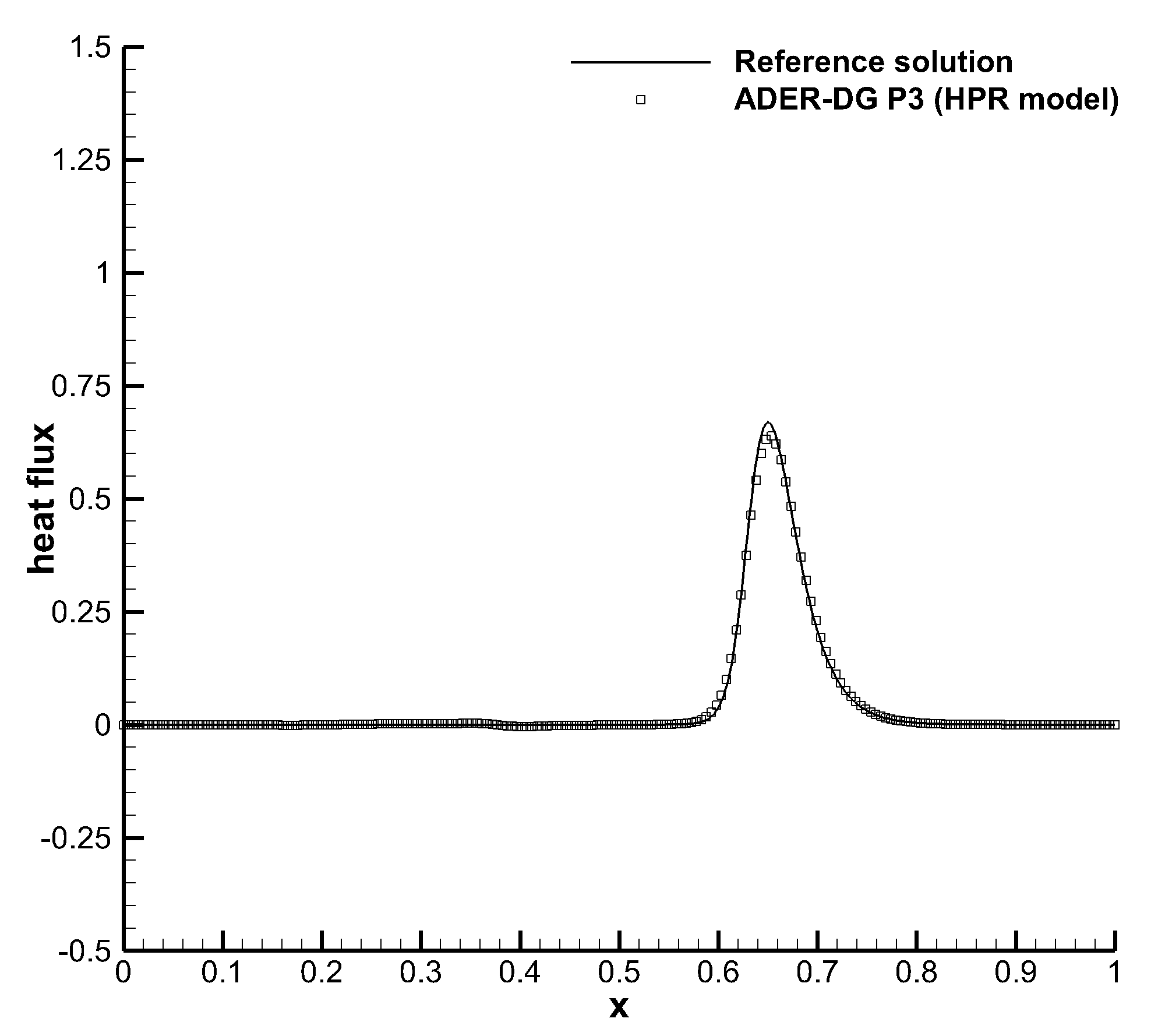}   
\end{tabular} 
\caption{Viscous shock with shock Mach number $M_s=2$ and Prandtl number $Pr=0.75$ at a final time of $t=0.2$. Viscous stress tensor component $\sigma_{11}$ (left) and heat flux (right). 
For the Navier-Stokes solution, the classical Fourier heat flux $q_1 = -\kappa T_x$ is shown, while for the HPR model, we plot $q_1 = \alpha^2 T J_1$. }  
\label{fig.SJ}
\end{center}
\end{figure}

For the special case of a \textit{stationary shock wave} at Prandtl number $Pr=0.75$ and constant viscosity, 
the compressible  Navier-Stokes equations can be reduced to one single ordinary differential 
equation (ODE) that can be solved analytically. The exact solution for the 
dimensionless velocity $\bar u = \frac{u}{M_s \, c_0}$ of this stationary shock wave is then given by the root of 
the following equation, see \cite{Becker1923,BonnetLuneau}:
\begin{equation} 
 \label{eqn.alg.u} 
  \frac{|\bar u - 1|}{|\bar u - \lambda^2|^{\lambda^2}} = \left| \frac{1-\lambda^2}{2} \right|^{(1-\lambda^2)} 
  \exp{\left( \frac{3}{4} \textnormal{Re}_s \frac{M_s^2 - 1}{\gamma M_s^2} x \right)},
\end{equation}
with
\begin{equation}
  \lambda^2 = \frac{1+ \frac{\gamma-1}{2}M_s^2}{\frac{\gamma+1}{2}M_s^2}.
\end{equation}
From eqn. \eqref{eqn.alg.u} one obtains the dimensionless velocity $\bar u$ as a function of $x$. 
The form of the viscous profile of the dimensionless pressure $\bar p = \frac{p-p_0}{\rho_0 c_0^2 M_s^2}$ is given by
the relation 
\begin{equation}
 \label{eqn.alg.p} 
  \bar p = 1 - \bar u +  \frac{1}{2 \gamma}
                         \frac{\gamma+1}{\gamma-1} \frac{(\bar u - 1 )}{\bar u} (\bar u - \lambda^2).  
\end{equation}
Finally, the profile of the dimensionless density $\bar \rho = \frac{\rho}{\rho_0}$ is found from the 
integrated continuity equation: $\bar \rho \bar u = 1$. In order to obtain an unsteady shock wave traveling 
into a medium at rest, it is sufficient to superimpose a constant velocity field $u = M_s c_0$ to the solution of the 
stationary shock wave found in the previous steps. We setup our computation with the exact solution of 
a shock wave (initially centered at $x=0.25$), traveling at $M_s=2.0$ to the right into a medium at rest.  
The values of the unperturbed fluid in front of the shock wave are chosen as $\rho_0=1$ (which also serves 
as reference density for the HPR model), $u_0=v_0=0$ and $p_0=1/\gamma$, hence $c_0 = 1$. The Reynolds number based 
on the shock speed and a unitary reference length ($L=1$) is defined as $\Re_s=\frac{\rho_0 \, c_0 \, M_s \, L }{\mu}$.   
The parameters of the first order HPR model are chosen as $\gamma = 1.4$, $c_v = 2.5$, $c_s=50$, $\mu=2 \cdot 10^{-2}$, 
$\alpha = 50$, $T_0=1$ and $\kappa = 9 \frac{1}{3} \cdot 10^{-2}$.  The resulting 
shock Reynolds number is $\Re_s=100$. The distortion tensor is initialized with $\mathbf{A}=\sqrt[3]{\rho} \, \mathbf{I}$
and the initial heat flux is set to $\mathbf{J}=0$, so that the system is started out of equilibrium. Simulations are carried out 
with an ADER-DG $P_3P_3$ scheme ($N=M=3$) up to a final time of $t=0.2$, by which the shock wave has traveled a distance of 0.4 to the right. 
The computational domain is given by $\Omega = [-0.5, 0.5] \times [-0.1, 0.1]$ and the mesh contains $100 \times 5$ elements. 
The comparison between the numerical solution of the first order HPR model and the exact solution of the 
compressible Navier-Stokes equations \eqref{eqn.alg.u} and \eqref{eqn.alg.p} is presented in Figs. \ref{fig.vshockprim} 
and \ref{fig.SJ}. One observes an excellent agreement of the viscous shock profile, apart from a small spurious wave at the 
left of the shock, which could be due to a small start-up error resulting from the non-equilibrium initial condition in the
variables $\mathbf{A}$ and $\mathbf{J}$.

\subsection{Viscous double Mach reflection problem} 
\label{sec:dmr} 

\begin{figure}[!htbp]
\begin{center}
\begin{tabular}{ccc} 
\includegraphics[width=0.32\textwidth]{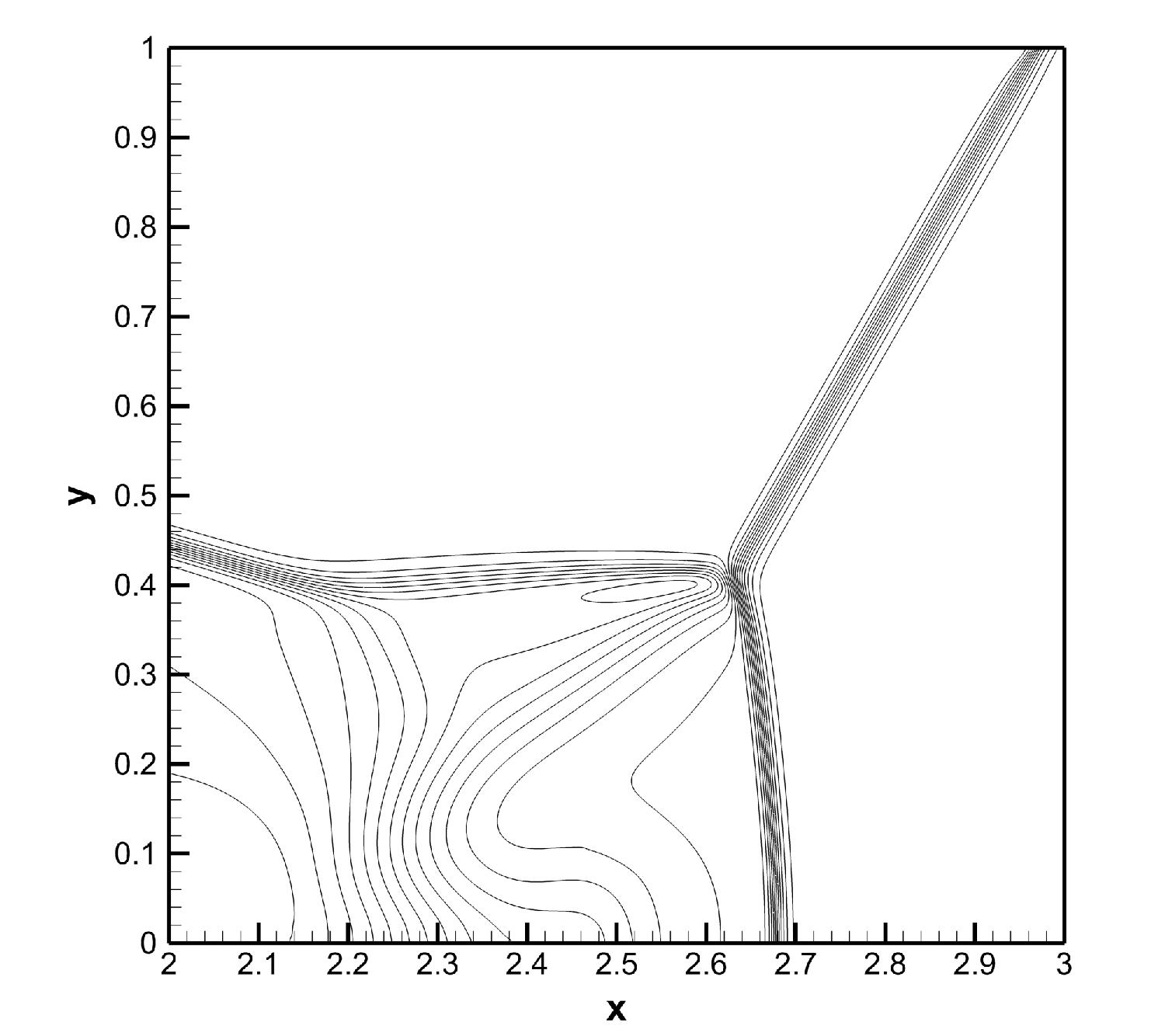}  &   
\includegraphics[width=0.32\textwidth]{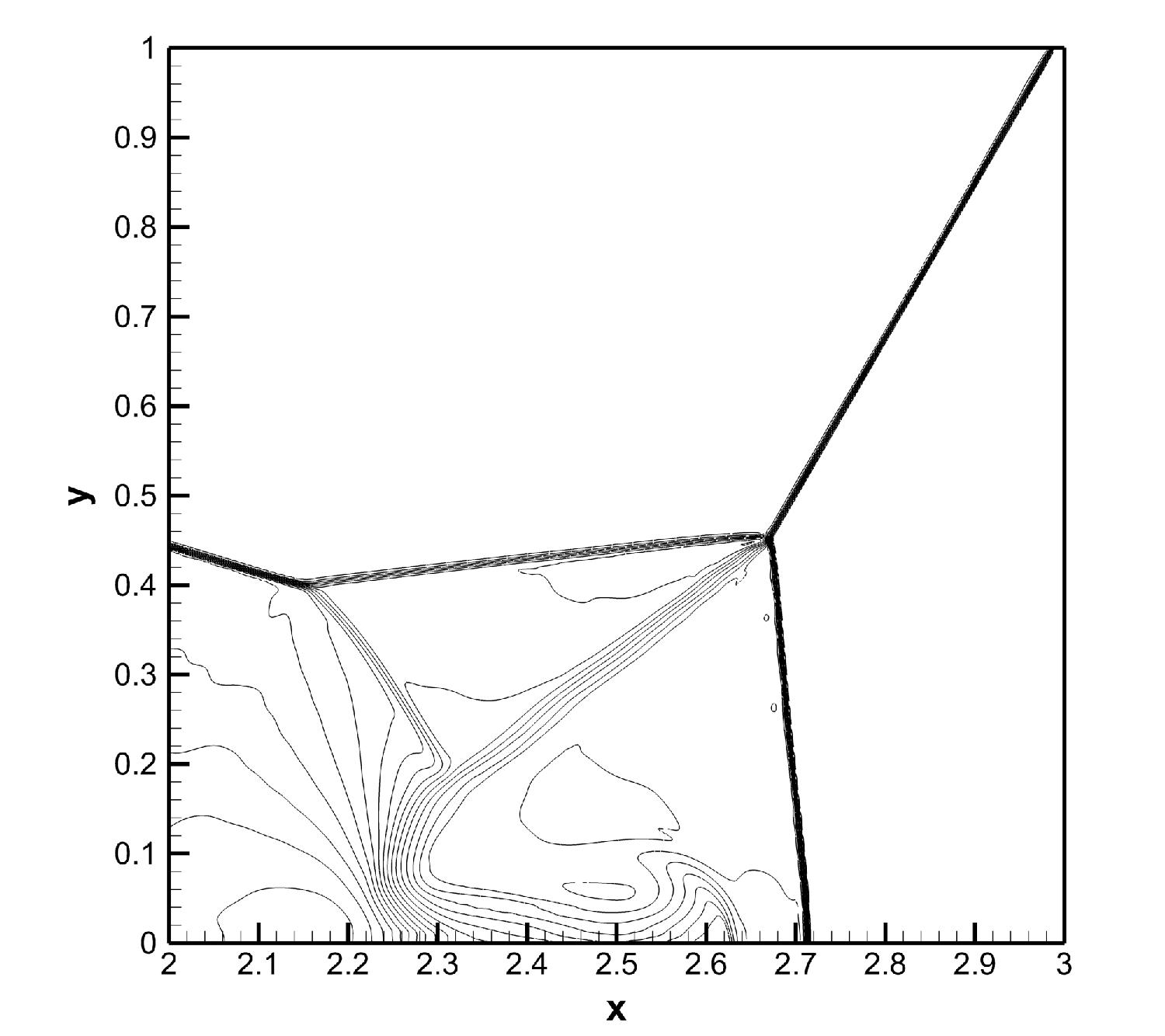}  & 
\includegraphics[width=0.32\textwidth]{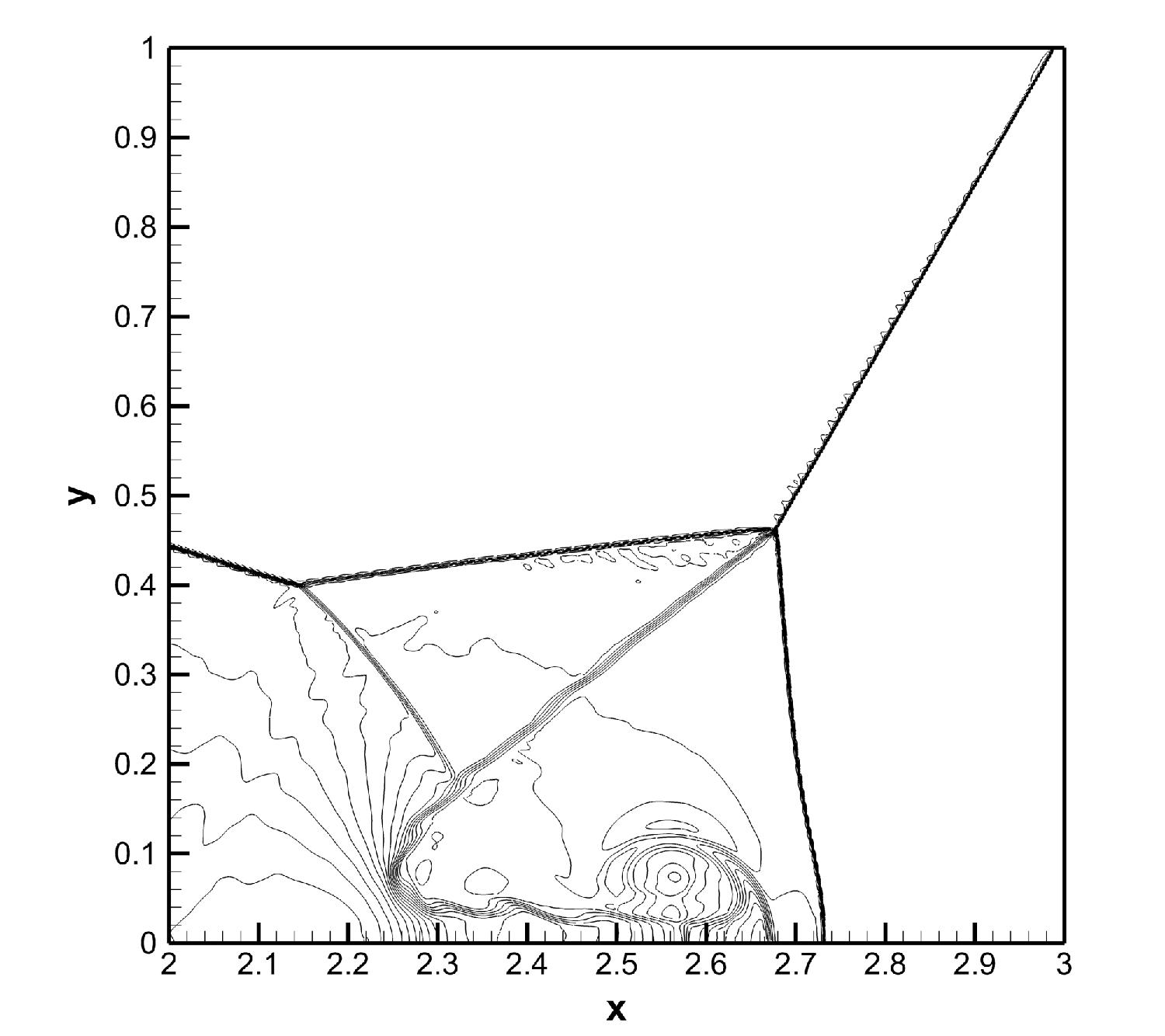}   
\end{tabular} 
\caption{Zoom into the viscous double Mach reflection problem at various shock Reynolds numbers at a final time of $t=0.2$ computed with a third order 
ADER-WENO finite volume method solving the HPR model. Left: $\mu = 10^{-1}$ ($\Re_s=100$). Middle: $\mu = 10^{-2}$, ($\Re_s=1000$). 
Right: Inviscid Euler reference calculation (Re$_s \to \infty$). 41 density contour levels in the interval $[1.5,17.5]$. }  
\label{fig.dmr}
\end{center}
\end{figure}

\begin{figure}[!htbp]
\begin{center}
\begin{tabular}{cc} 
\includegraphics[width=0.4\textwidth]{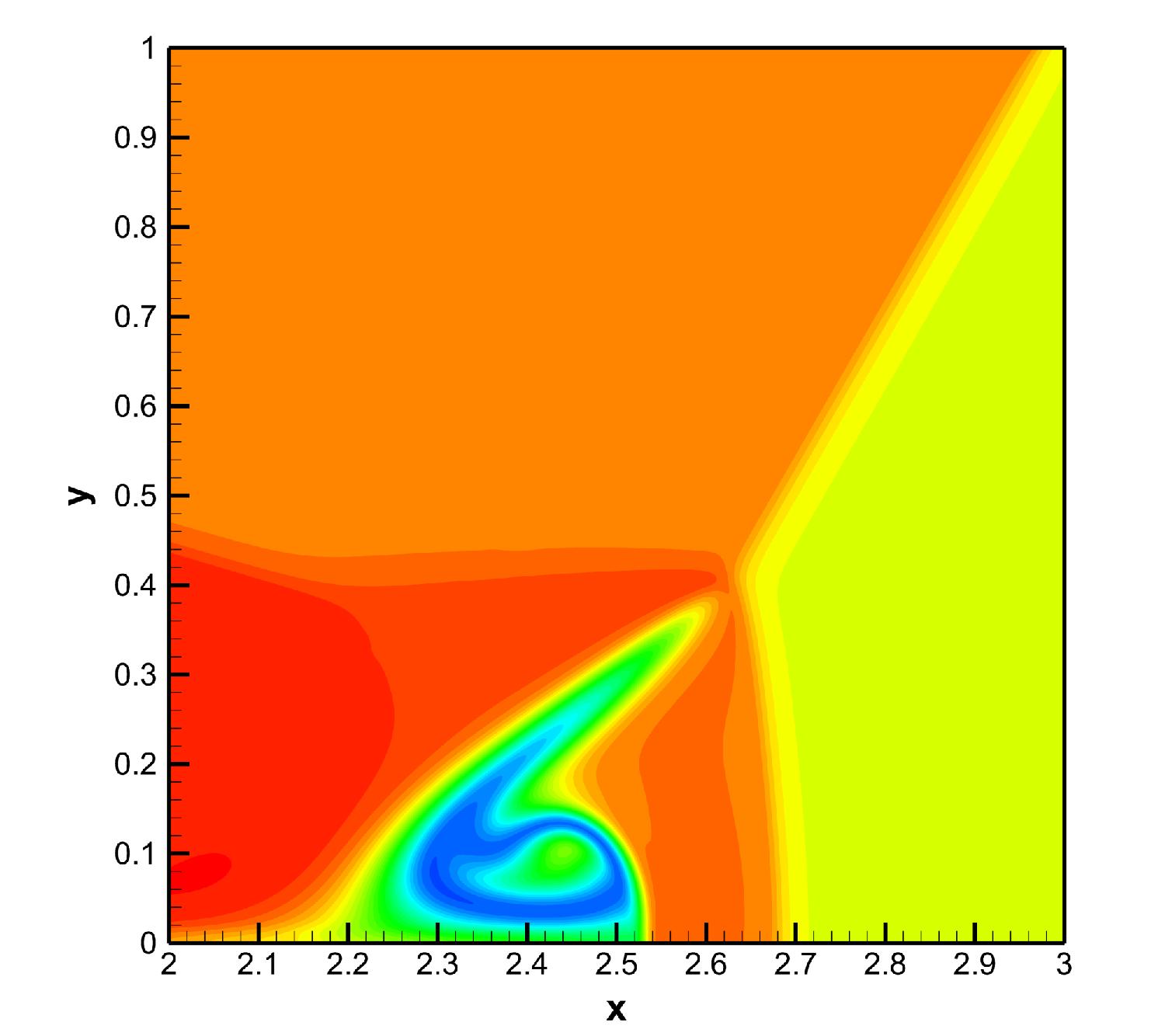}  &   
\includegraphics[width=0.4\textwidth]{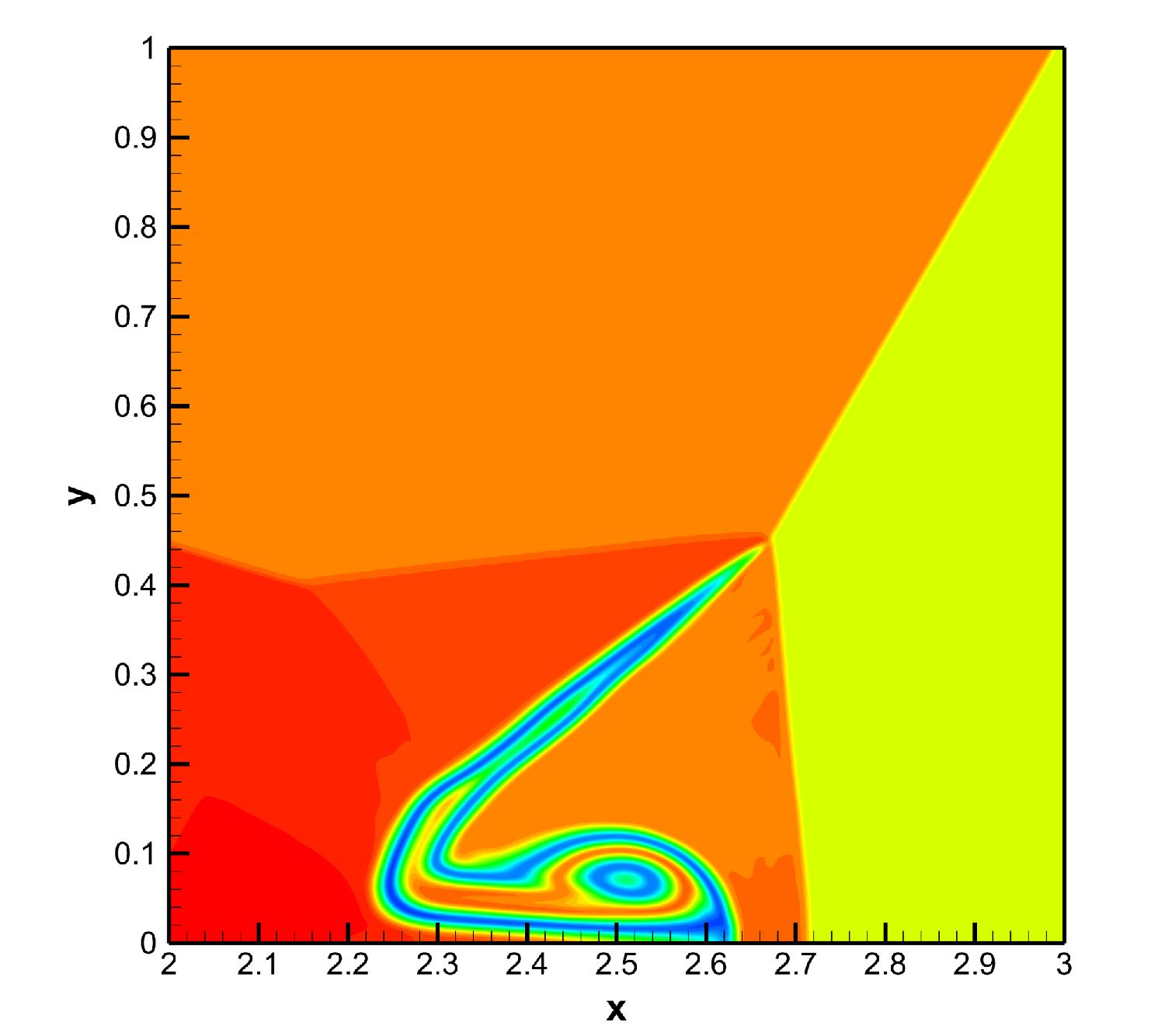}   
\end{tabular} 
\caption{Distortion tensor component $A_{11}$ of the HPR model for the viscous double Mach reflection problem at a final time of $t=0.2$, 
computed with a third order ADER-WENO finite volume method. Left: $\mu = 10^{-1}$ (Re$_s=100$). Right: $\mu = 10^{-2}$, (Re$_s=1000$). 31 contour levels 
are shown in the interval $[-2.25,+2.25]$.}  
\label{fig.dmra11}
\end{center}
\end{figure}

In this section we run a viscous version of the 2D double Mach reflection problem of a strong shock, which has been 
originally proposed for the compressible Euler equations by Woodward and Colella in \cite{woodwardcol84}. This test problem involves 
a Mach $10$ shock that hits a $30^\circ$ ramp. Using the Rankine-Hugoniot conditions of the compressible Euler equations we can deduce 
the initial conditions for the flow variables in front of and behind the shock wave as   
\begin{eqnarray}
(\rho, u, v, p)( \x,t=0) =
\left\{
\begin{array}{cll}
  \frac{1}{\gamma}(8.0, 8.25, 0.0, 116.5), \quad & \text{ if } & \quad x'<0.1, \\
  (1.0, 0.0, 0.0, \frac{1}{\gamma}),       \quad & \text{ if } & \quad x'\geq 0.1, 
\end{array}
\right.
\end{eqnarray}
where $x'$ is the coordinate in a rotated coordinate system. Reflecting slip wall boundary conditions are prescribed  
on the bottom and the exact solution of an isolated moving oblique shock wave with shock Mach 
number $M_s=10$ is imposed on the upper boundary. Inflow and outflow boundary conditions are prescribed on the left side 
and the right side, respectively. \\ 
The computational domain is given by $\Omega = [0;3.5] \times [0;1]$ and the computational grid uses a characteristic  
length of $h=1/400$, leading to $1400 \times 400$ computational cells. We solve this problem with a third order $P_0P_2$ ADER-WENO finite 
volume scheme ($N=0$, $M=2$). The parameters of the HPR model are: $\gamma=1.4$, $c_v=2.5$, $\rho_0=1$, $c_s=20$, 
$\alpha^2=200$, $T_0=1$ and $\kappa= \gamma c_v \mu / \textnormal{Pr}$ with a Prandtl number of Pr=$0.75$. 
The initial condition for the distortion tensor is chosen as $\mathbf{A} = \sqrt[3]{\rho} \, \mathbf{I}$ and the heat flux is
initialized with $\mathbf{J}=0$. The computational results are depicted in Fig. \ref{fig.dmr} for the density variable at a final time of $t=0.2$ 
using two different values for the viscosity coefficient: $\mu=10^{-1}$, which corresponds to a shock Reynolds number of $\Re_s=100$, and $\mu=10^{-2}$, 
corresponding to $\Re_s=1000$. An inviscid reference calculation of the compressible Euler equations with the same ADER-WENO scheme on the same 
grid is also provided in Fig. \ref{fig.dmr}, to show that the missing flow features in the HPR model are actually due to the presence of physical
viscosity and not due to the effect of numerical diffusion. Overall one can observe that the typical flow structures like the incident shock wave, 
the reflected shock wave and the Mach stem are well reproduced. Furthermore, the typical mushroom-type flow structure close to the $x$-axis is 
also present in the viscous computations carried out with the HPR model. However, in this test problem, a rather large physical viscosity has been 
added, hence preventing the development of any unstable small-scale flow structures as observed in \cite{balsarashu,Dumbser2014,Zanotti2015a} 
and as obtained in the inviscid reference calculation. For other numerical results concerning the viscous double Mach reflection problem, 
see \cite{ADERNSE}. In Fig. \ref{fig.dmra11} we also provide a visualization of the distortion tensor component $A_{11}$, which clearly
indicates the presence of the shear layers and the mushroom-type flow feature, underpinning just once more the value of $\mathbf{A}$ for the 
purpose of flow visualization. 


\subsection{Application to solid mechanics} 
\label{sec:solid}

The main key advantage of the HPR model \eqref{eqn.conti}-\eqref{eqn.energy} is its capability to describe in \textit{one single PDE system} the two main branches of 
\textbf{continuum mechanics}, namely \textbf{fluid mechanics} \underline{\textit{and}} \textbf{solid mechanics}. We explicitly stress here that the test problem 
proposed in this section \textbf{cannot} be solved with the conventional Navier-Stokes equations, since it is a typical test problem of \textit{solid mechanics}!  
Here we consider a typical benchmark problem used in computational seismology \cite{komatitsch1998,komatitsch1999,gij1}, 
consisting in the propagation of a wave in a linear elastic solid with free surface. The problem was first discussed and solved by Lamb in 1904, see \cite{Lamb1904} and is therefore often called \textit{Lamb's problem} in the literature. The problem consists in a point force acting on an elastic solid, perpendicular to a free surface. 
In our particular setup, the computational domain is given by 
$\Omega = [-2000, 2000] \times [-2000, 0]$ and we add the following point source to the right hand side of the momentum equation \eqref{eqn.momentum}: 
\begin{equation}
  \mathbf{S}(\mathbf{x},t) = \rho_0 a_1 \left( \halb + a_2 (t - t_D)^2 \right) \exp \left( a_2 (t - t_D)^2 \right) \delta(\mathbf{x}-\mathbf{x}_s) \mathbf{e}_y,   
 \label{eqn.src} 	
\end{equation}  
with the Dirac delta distribution $\delta(\mathbf{x})$, the unit vector pointing in $y$-direction $\mathbf{e}_y=(0,1)$, and the following source parameters: 
$t_D = 0.08$, $a_1 = -2000$, $a_2 = -(\pi f_c)^2$, $fc = 14.5$. The source is located in 
$\mathbf{x}_s = (0,-1)$, hence slightly below the free surface, which is located at $y=0$. At the free surface, the shear stresses $\sigma_{12} = \sigma_{21}$ and the normal stress $\sigma_{22}$
vanish. The parameters of the HPR model with stiffened gas EOS are $\rho_0 = 2200$, $c_v=1$, $\gamma=2$, $c_0 = 2385.160721$, $c_s = 1847.5$, $\alpha=0$, $\tau_1 \to \infty$ and 
$\tau_2 \to \infty$. In practice, we set $\tau_1 = \tau_2 = 10^{20}$. The initial condition is chosen as $\rho = \rho_0$, $u=v=0$, $p=0$, $\mathbf{A}=\mathbf{I}$ and $\mathbf{J}=0$. 
With the previous parameters, the resulting sound speed of longitudinal pressure waves is $c_L = \sqrt{c_0^2 + 4/3 c_s^2} = 3200$.  
Simulations are performed with a $P_4P_4$ ADER-DG scheme ($N=M=4$) up to a final time of $t=1.3$ using a grid composed of $200 \times 100$ elements. 

\begin{figure}[!htbp]
\begin{center}
\begin{tabular}{c} 
\includegraphics[width=0.7\textwidth]{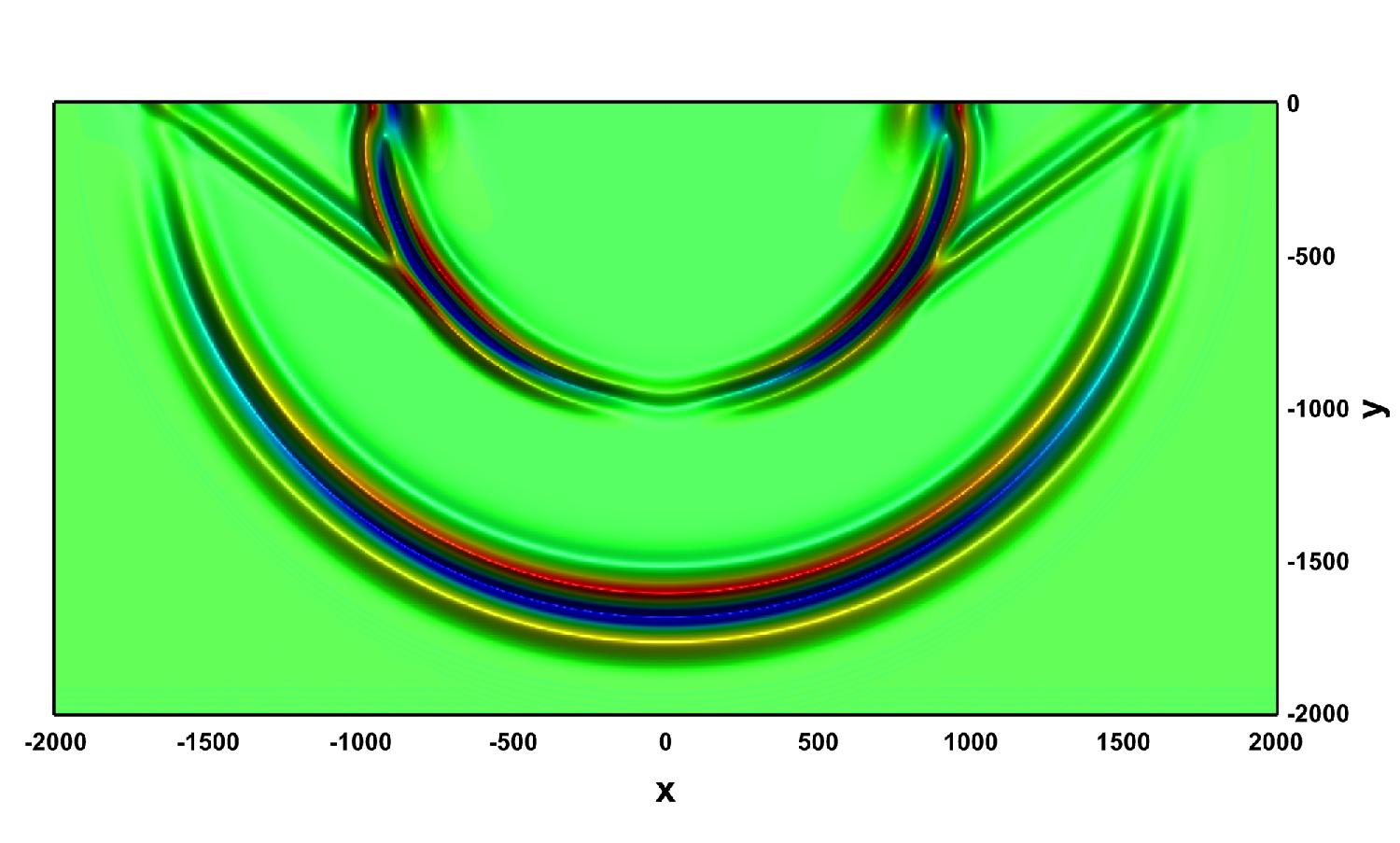}  \\  
\includegraphics[width=0.7\textwidth]{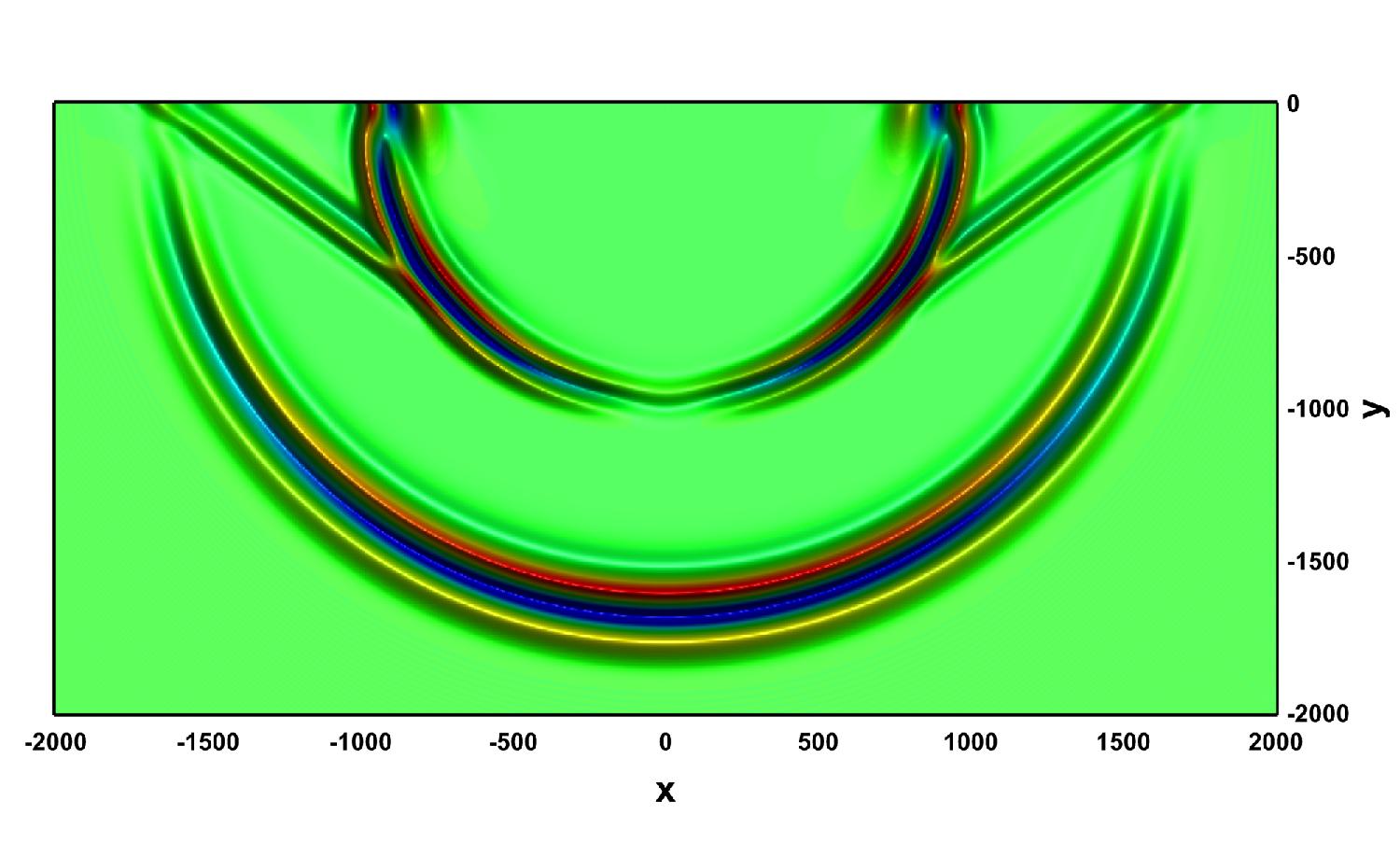}   
\end{tabular} 
\caption{Lamb's problem in 2D at a final time of $t=0.6$ computed with an ADER-DG P4 method. Contour colors of $v$ for the HPR model (top) and the reference solution 
based on the equations of linear elasticity (bottom). }  
\label{fig.Lamb2D}
\end{center}
\end{figure}

\begin{figure}[!htbp]
\begin{center}
\begin{tabular}{c} 
\includegraphics[width=0.95\textwidth]{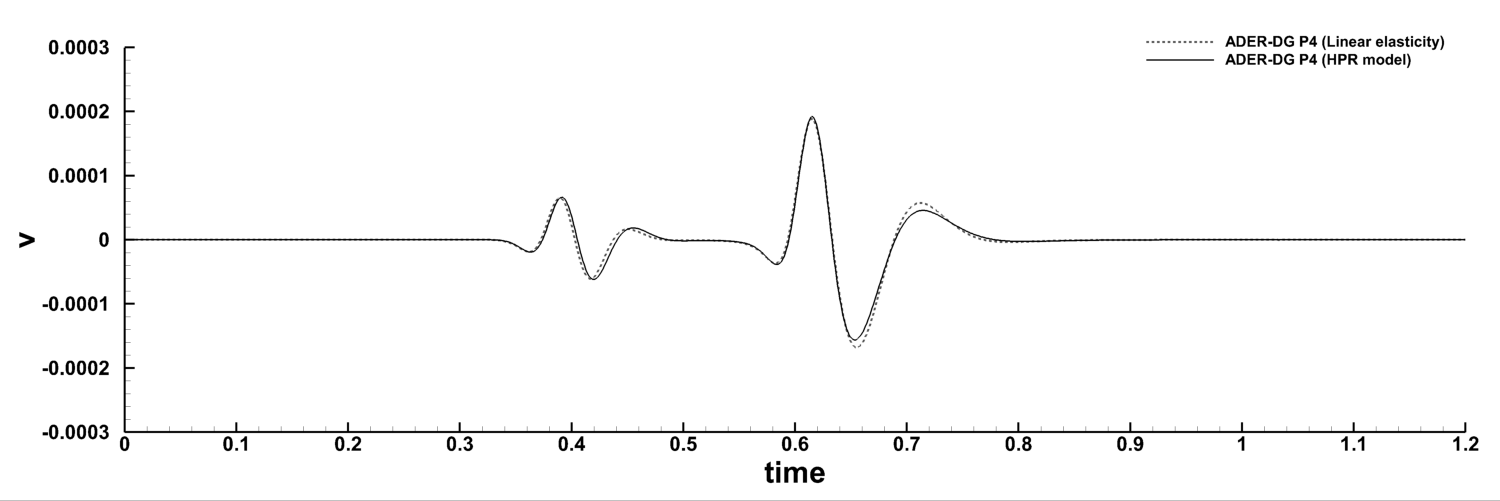}  
\end{tabular} 
\caption{Time signal of the velocity component $v$ recorded at $(990,0)$ for Lamb's problem in 2D. Comparison between the HPR model and the reference solution based on 
the equations of linear elasticity. }  
\label{fig.Lamb.seis}
\end{center}
\end{figure}

In order to have a direct comparison with classical theory of linear elasticity, we solve the problem again with the same ADER-DG scheme on the same grid, but in the second run
we directly solve the classical equations of linear elasticity in velocity-stress formulation \cite{komatitsch1998,komatitsch1999,gij1},  
\begin{eqnarray}
  \frac{\partial}{\partial t} \sigma_{xx} - (\lambda + 2\mu) \frac{\partial}{\partial x} u - \lambda \frac{\partial}{\partial y} v & = &  0, \nonumber \\  
  \frac{\partial}{\partial t} \sigma_{yy} - \lambda \frac{\partial}{\partial x} u - (\lambda + 2\mu) \frac{\partial}{\partial y} v & = &  0, \nonumber \\  
  \frac{\partial}{\partial t} \sigma_{xy} - \mu \frac{\partial}{\partial x} v - \mu \frac{\partial}{\partial y} u & = &  0, \nonumber \\  
 \rho \frac{\partial}{\partial t} u - \frac{\partial}{\partial x} \sigma_{xx} -  \frac{\partial}{\partial y} \sigma_{xy} & = &  0, \nonumber \\  
 \rho \frac{\partial}{\partial t} v - \frac{\partial}{\partial x} \sigma_{xy} -  \frac{\partial}{\partial y} \sigma_{yy} & = &  0. 
\label{eqn.linwave} 
\end{eqnarray} 
The two Lam\'e constants in \eqref{eqn.linwave} $\lambda = 7.509672500 \cdot 10^{9}$ and $\mu = 7.509163750 \cdot 10^{9}$ are chosen in order to obtain the same wave  propagation speeds as in the HPR model, i.e. $c_p = \sqrt{(\lambda+2\mu)/\rho}=3200$ for longitudinal pressure waves and $c_s = \sqrt{\mu/\rho} = 1847.5$ for shear waves. The density is, of course, also set 
to $\rho=2200$. The same point source \eqref{eqn.src} as in the HPR model is added to the right hand side of the last equation of \eqref{eqn.linwave}. The 
computational results obtained with both models are compared against each other in Fig. \ref{fig.Lamb2D}, where we note an excellent agreement of the two
computed wave fields. In Fig. \ref{fig.Lamb.seis} we compare the velocity signal $v$ at an observation point located at $x=(990,0)$, where also a very good agreement
between the nonlinear HPR model and the reference solution based on the linear elastic wave equations can be observed.

\section{Conclusion} \label{sec:conclusion}

To the knowledge of the authors, this is the first time that a numerical method has been applied to the full first order hyperbolic Peshkov-Romenski  
model \cite{PeshRom2014} with heat conduction, see Eqns. \eqref{eqn.conti}-\eqref{eqn.energy}. The proposed family of high order one-step ADER finite volume and 
ADER discontinuous Galerkin finite element schemes is able to discretize rather general hyperbolic systems of partial differential equations with non-conservative 
products and stiff source terms and has been applied in the frame of the HPR model to a large set of different test problems, ranging from viscous compressible 
fluids to elastic solids. We have shown numerical convergence results, as well as detailed comparisons with different analytical and numerical reference solutions. 
It is very important to stress again that the first order HPR model is able to represent the basic equations of continuum mechanics in a unified manner, including 
\textit{fluid mechanics} and \textit{solid mechanics} as two special limiting cases of the \textit{same} mathematical model. 
The nonlinear material behavior is entirely governed by the equation of state and by the strain relaxation mechanism $\boldsymbol{\psi}$. 
To the knowledge of the authors, such a universal formulation of continuum mechanics in a first order hyperbolic system is \textit{unique} and has never been 
tackled before with high order shock capturing methods for hyperbolic conservation laws. 

Future applications will concern the extension of the numerical method to moving unstructured meshes in the frame of ADER-WENO-ALE schemes 
\cite{LagrangeNC,Lagrange3D}, as well as to non-Newtonian fluids and complex visco-plastic solids. The use of high order schemes for hyperbolic 
PDE on space-time adaptive meshes, as outlined in \cite{AMR3DCL,AMR3DNC,Zanotti2015a}, might also become useful in near future in combination 
with the HPR model in the context of crack generation and crack propagation in nonlinear solid mechanics. 

Since the HPR model has already \textit{finite wave speeds} for all involved physical processes, i.e. heat and mass transport, 
as well as viscous momentum transport, future research will be carried out in order to extend it also to the \textit{relativistic regime}, following the promising 
investigations by \cite{Gundlach2012}. 

\section*{Acknowledgments}
\vspace{-7mm} 
\noindent
\hspace{-3mm} 
\begin{tabular}{lr} 
\begin{minipage}[c]{0.8\textwidth}
The research presented in this paper has been financed by the European Research Council (ERC) under the
European Union's Seventh Framework Programme (FP7/2007-2013) with the research project \textit{STiMulUs}, 
ERC Grant agreement no. 278267. 
M.D. and O.Z. have further received funding from the European Union's Horizon 2020 Research and 
Innovation Programme under the project \textit{ExaHyPE}, grant agreement number no. 671698 (call
FETHPC-1-2014). 
\end{minipage}
& 
\begin{minipage}[c]{0.2\textwidth}
\includegraphics[angle=0,width=0.95\textwidth]{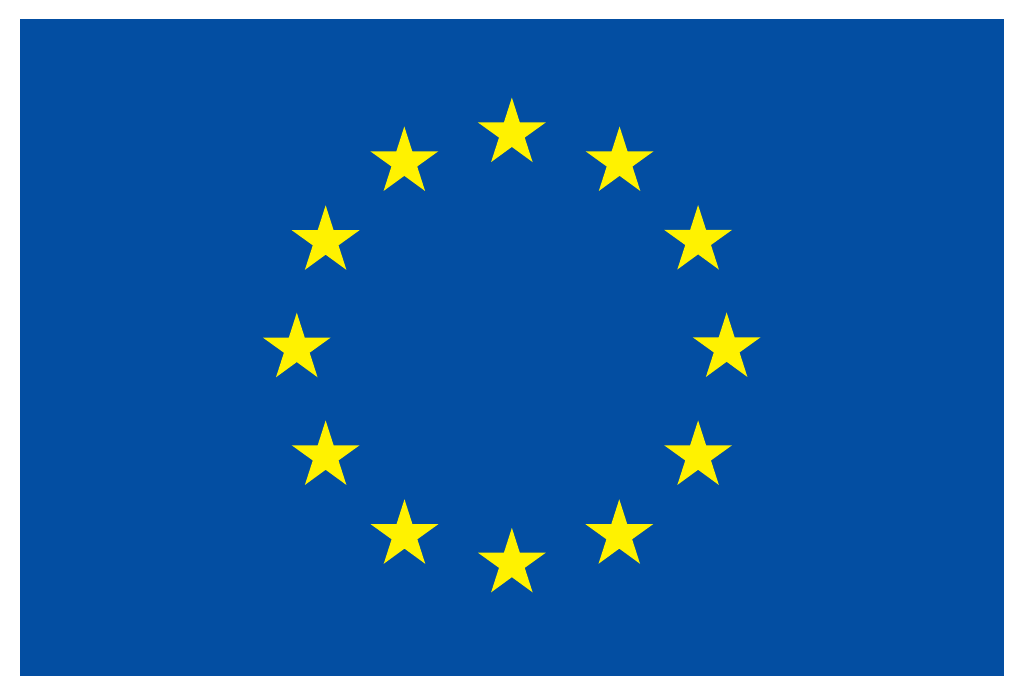}
\end{minipage}  
\end{tabular} 
The authors would like to acknowledge PRACE for awarding access to the SuperMUC 
supercomputer based in Munich, Germany, at the Leibniz Rechenzentrum (LRZ). 

\vspace{3mm} 

Last, but not least, the authors are very grateful to Sergei Konstantinovich Godunov, for his groundbreaking and 
very inspiring seminal ideas that are at the basis of both, the theoretical as well as the numerical framework employed in this paper.



\clearpage

\appendix
\section{Eigenvalues of matrices $ \mathsf{A}_k $}\label{ap.matrix}
In this section, we give the formulas for the eigenvalues of matrices $ \mathsf{A}_k $ of the viscous subsystem of \eqref{eqn.HPR}. Thus, if the heat conducting effect is ignored, then matrices $ \mathsf{A}_k $ look as follows
\begin{equation}
\mathsf{A}_k(\QQ)=\left (\begin{array}{cccccccc}
	v_k &  0  &  z_1   &  z_2   &  z_3   &   0    &   0    &   0    \\
	 0  & v_k &  Z_1   &  Z_2   &  Z_3   &   0    &   0    &   0    \\
	R_1 & P_1 &  v_k   &   0    &   0    & X_{11} & X_{12} & X_{13} \\
	R_2 & P_2 &   0    &  v_k   &   0    & X_{21} & X_{32} & X_{23} \\
	R_3 & P_3 &   0    &   0    &  v_k   & X_{31} & X_{32} & X_{33} \\
	 0  &  0  & A_{11} & A_{12} & A_{13} &  v_k   &   0    &   0    \\
	 0  &  0  & A_{21} & A_{22} & A_{23} &   0    &  v_k   &   0    \\
	 0  &  0  & A_{31} & A_{32} & A_{33} &   0    &   0    &  v_k
\end{array}\right ),
\end{equation}
where (no summation over repeated index $ k $)
\[\rho R_i=\dfrac{\partial \sigma_{ik}}{\partial \rho},\ \ \ \ \rho P_i=\delta_{ik}+\dfrac{\partial \sigma_{ik}}{\partial p},\ \ \ \ \ \rho X_{ij}=\dfrac{\partial \sigma_{ik}}{\partial A_{jk}}\]
\[z_k=\rho,\ \ \ \ z_m=z_n=0,\ \ \ \ m\neq k, n\neq k,\]
\[Z_k\dfrac{\partial \rho E_1}{\partial p}=\rho E_1 + \rho E_2  - \rho \dfrac{\partial \rho E_1}{\partial \rho} + p + \sigma_{kk},\ \ \ \ Z_m\dfrac{\partial \rho E_1}{\partial p}=\sigma_{mk},\ \  \ \ Z_n\dfrac{\partial \rho E_1}{\partial p}=\sigma_{nk},\ \ m\neq k, n\neq k.\]

The eigenvalues of the matrix $ \mathsf{A}_k(\QQ) $ are given by the formulas
\begin{equation}
v_k-\lambda_3,\ \ v_k-\lambda_2,\ \ v_k-\lambda_1,\ \ v_k, \ \ v_k, \ \ v_k+\lambda_1,\ \ v_k+\lambda_2,\ \ v_k+\lambda_3,
\end{equation}
where $ \lambda_1\leq\lambda_2<\lambda_3 $ are three eigenvalues of the 3-by-3 matrix 
\[\WW =\WW_1\WW_2,\]
where
\[
\WW_1=\left(
\begin{array}{ccccc}
	R_1 & P_1 & X_{11} & X_{12} & X_{13} \\
	R_2 & P_2 & X_{21} & X_{22} & X_{23} \\
	R_3 & P_3 & X_{31} & X_{32} & X_{33}
\end{array}
\right),\ \ \WW_2=\left(
\begin{array}{ccc}
	 z_1   &  z_2   &  z_3   \\
	 Z_1   &  Z_2   &  Z_3   \\
	A_{11} & A_{12} & A_{13} \\
	A_{21} & A_{22} & A_{23} \\
	A_{31} & A_{32} & A_{33}
\end{array}
\right).
\]

Because of that the EOS is assumed to be a convex function of the state variables $ \QQ $ \cite{PeshRom2014} and hence the HPR model is symmetric hyperbolic, the eigenvalues are thus assumed to be real, and can be found by analytical formulas as the roots of the cubic polynomial $ \det(\WW-\lambda\II)=0 $. They are
\[\lambda_k=\sqrt{\beta_k+{\rm tr}(\WW)/3},\]
where $ \beta_k=2\sqrt{-a/3}\cos((\phi+2(k-1)\pi))/3) $ and 
\[\phi=\left\{\begin{array}{l}
{\rm acos(-\sqrt{-27b^2/(4a^3)})},\ \ {\rm if}\ b>0\\[2mm]
{\rm acos(\sqrt{-27b^2/(4a^3)})},\ \ {\rm if}\ b<0
\end{array}\right.\]
Here, $ a=(I_1^2-3 I_2)/6 $ and $  b= \left (5 I_1^3/9-I_1I_2-6 \det(\WW)\right )/6 $, $I_1={\rm tr}(\WW)$ and $ I_2={\rm tr}(\WW^2) $.

\clearpage

\end{document}